\documentclass[11pt]{article}

\usepackage{smile}
\usepackage{enumerate,mathrsfs}
\usepackage{algorithm,algorithmic}
\usepackage{extarrows}
\usepackage[OT1]{fontenc}
\usepackage[colorlinks,linkcolor=red,anchorcolor=blue,citecolor=blue]{hyperref}
\usepackage{fullpage}
\usepackage[protrusion=false,expansion=true]{microtype}
\usepackage{graphicx}
\usepackage{comment}
\usepackage{stmaryrd}

\newcommand{\mT}{\mathscr{T}}
\newcommand{\mX}{\mathscr{X}}
\newcommand{\mR}{\mathbb{R}}

\newcommand{\mP}{\mathbb{P}}

\newcommand{\mI}{\mathbb{I}}
\newcommand{\tv}{\text{vec}}

\newcommand{\supp}{\text{supp}}
\newcommand{\bdm}{\boldsymbol{m}}

\newtheorem{Definition}{Definition}
\newtheorem{Condition}{Condition}
\newtheorem{Corollary}{Corollary}
\newtheorem{Theorem}{Theorem}

\newtheorem{Remark}{Remark}
\newtheorem{Lemma}{Lemma}

\newcommand{\normmm}[1]{{\left\vert\kern-0.25ex\left\vert\kern-0.25ex\left\vert #1
		\right\vert\kern-0.25ex\right\vert\kern-0.25ex\right\vert}}

\begin{document}

\title{ \bf Sparse and Low-rank Tensor Estimation via Cubic Sketchings}

\author
{
	Botao Hao\thanks{Postdoctoral Researcher, Department of Electrical Engineering, Princeton University. E-mail: haobotao000@gmail.com.},~
	Anru Zhang\thanks{Assistant Professor, Department of Statistics, University of Wisconsin-Madison. E-mail: anruzhang@stat.wisc.edu.},~
	Guang Cheng\thanks{Professor, Department of Statistics, Purdue University, West Lafayette. E-mail: chengg@purdue.edu. Guang Cheng would like to acknowledge support by NSF DMS-1712907, DMS-1811812, DMS-1821183, and Office of Naval Research (ONR N00014-18-2759). In addition, Guang Cheng is a member of Institute for Advanced Study, Princeton and visiting Fellow of SAMSI for the Deep Learning Program in the Fall of 2019; he would like to thank both Institutes for their hospitality.}
}

\date{}
\maketitle

\maketitle
\begin{abstract}
In this paper, we propose a general framework for sparse and low-rank tensor estimation from cubic sketchings. A two-stage non-convex implementation is developed based on sparse tensor decomposition and thresholded gradient descent, which ensures exact recovery in the noiseless case and stable recovery in the noisy case with high probability. The non-asymptotic analysis sheds light on an interplay between optimization error and statistical error. The proposed procedure is shown to be rate-optimal under certain conditions. As a technical by-product, novel high-order concentration inequalities are derived for studying high-moment sub-Gaussian tensors. An interesting tensor formulation illustrates the potential application to high-order interaction pursuit in high-dimensional linear regression.

\end{abstract}


{
\hypersetup{linkcolor=black}
\tableofcontents
}

\section{Introduction}\label{sec:intro}

The rapid advance in modern scientific technology gives rise to a wide range of high-dimensional tensor data \cite{K08, kolda2009}. Accurate estimation and fast communication/processing of tensor-valued parameters are crucially important in practice. For example, a tensor-valued predictor which characterizes the association between brain diseases and scientific measurements becomes the point of interest \cite{zhou2013,li2013tucker, SL17}. Another example is the tensor-valued image acquisition algorithm that can considerably reduce the number of required samples by exploiting the compressibility property of signals \cite{CC13, friedland2014compressive}. 

The following tensor estimation model is widely considered in recent literatures,
\begin{equation}\label{eq:general-model}
y_i = \langle \mT^\ast, \mX_i\rangle + \epsilon_i, \quad i=1,\ldots, n.
\end{equation}
Here, $\mX_i$ and $\epsilon_i$ are the measurement tensor and the noise, respectively. The goal is to estimate the unknown tensor $\mT^\ast$ from measurements $\{y_i, \mX_i\}_{i=1}^n$. A number of specific settings with varying forms of $\mX_i$ have been studied, e.g., tensor completion \cite{liu2013,yuan2016tensor,yuan2017incoherent,zhang2016cross,montanari2016spectral, ghadermarzy2018near, zhang2016exact, bengua2017efficient}, tensor regression \cite{SL17,zhou2013, li2013tucker,raskutti2015convex,chen2016non,li2017parsimonious,zhang2019islet}, multi-task learning \cite{romera2013multilinear}, etc.

\begin{figure}[h!]
	\centering
		\vspace{0.1cm}
		\includegraphics[scale=0.5]{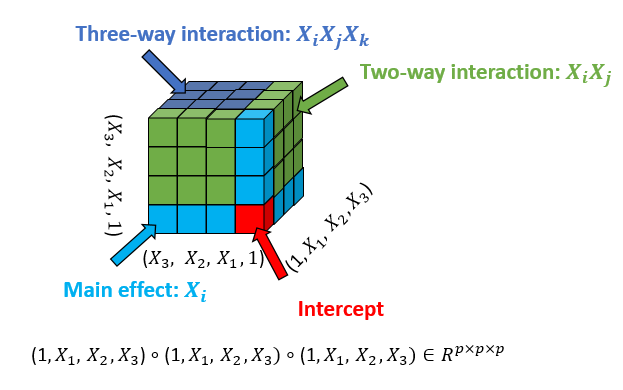}\;\;\;\;\;\;\;\;\;
		\includegraphics[scale=0.3]{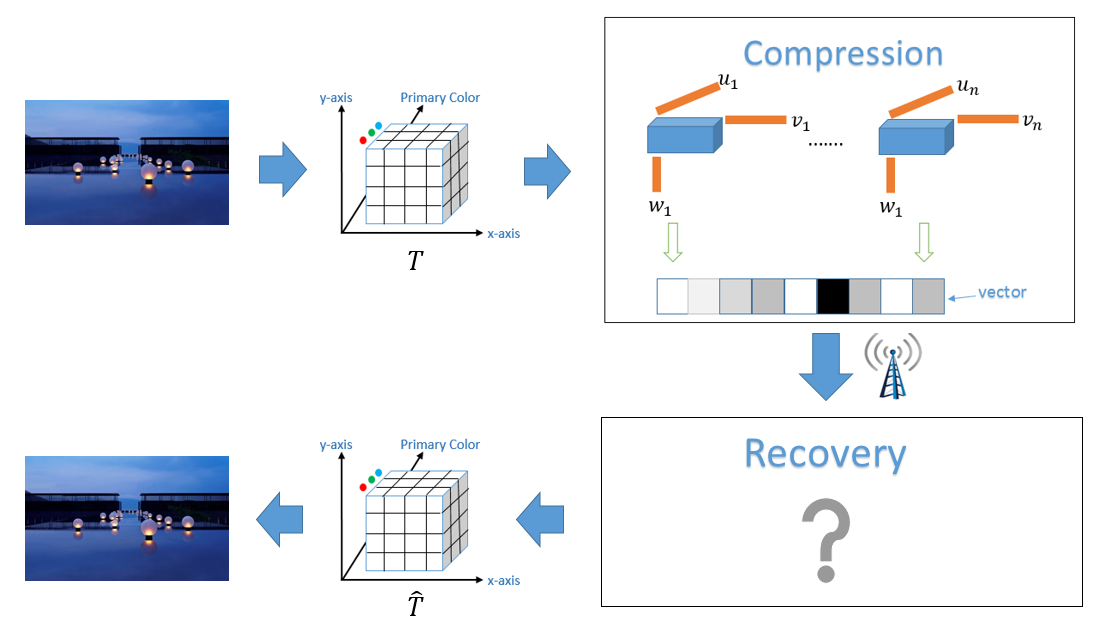}
		\vspace{-0.1in}
	\caption{Illustration for interaction reformulation and tensor image/video compression}
	\label{fig:illurstration}
\end{figure}

In this paper, we focus on the case that the measurement tensor can be written in a cubic sketching form. For example, $\mX_i = \bx_i\circ\bx_i\circ\bx_i$ or $\mX_i=\bu_i\circ \bv_i\circ \bw_i$, depending on whether $\mT^*$ is symmetric or not. The cubic sketching form of $\mX_i$ is motivated by a number of applications.
\begin{itemize}
	\item \emph{Interaction effect estimation:} High-dimensional high-order interaction models have been considered under a variety of settings \cite{bien2013lasso, hao2014interaction,fan2016interaction, basu2018iterative}. By writing $\mX_i = \bx_i\circ\bx_i\circ\bx_i$, we find that the interaction model has an interesting tensor representation (see left panel of Figure \ref{fig:illurstration}) which allows us to estimate high-order interaction terms using tensor techniques. This is in contrast with the existing literature that mostly focused on pair-wise interactions due to the model complexity and computational difficulties. More detailed discussions will be provided in Section \ref{sec:interaction}.
	\item \emph{High-order imaging/video compression:} High-order imaging/video compression is an important task in modern digital imaging with various applications (see right panel of Figure \ref{fig:illurstration}), such as hyper-spectral imaging analysis \cite{li2010tensor} and facial imaging recognition \cite{vasilescu2003multilinear}. 
One could use Gaussian ensembles for compression such that each entry of $\mX_i$ is i.i.d. randomly generated \cite{zhou2013, raskutti2015convex,chen2016non}. In contrast, the non-symmetric cubic sketchings, i.e., $\mX_i=\bu_i\circ \bv_i\circ \bw_i$, reduce the memory storage from $O(np_1p_2p_3)$ to $O(n(p_1+p_2+p_3))$ ($n$ is the sample size and $(p_1, p_2, p_3)$ is the tensor dimension), but still preserve the optimal statistical rate. More detailed discussions will be provided in Section \ref{sec:non-symmetrc}. 
\end{itemize}
In practice, the total number of measurements $n$ is considerably smaller than the number of parameters in the unknown tensor $\mT^*$, due to all kinds of restrictions such as time and storage. Fortunately, a variety of high-dimensional tensor data possess intrinsic structures, such as low-rankness \cite{kolda2009} and sparsity \cite{Sun2016}. This could highly reduce the effective dimension of the parameter and make the accurate estimation possible. Please refer to \eqref{eq:mT_low-rank} and \eqref{eq:low_rank_nonsym} for low-rankness and sparsity assumptions.

In this paper, we propose a computationally efficient non-convex optimization approach for sparse and low-rank tensor estimation via cubic-sketchings. Our procedure is two-stage: 
\begin{itemize}
	\item[(i)] obtain an initial estimate via the method of tensor moment (motivated by high-order Stein's identity), and then apply sparse tensor decomposition to the initial estimate to output a warm start;
	\item[(ii)] use a thresholded gradient descent to iteratively refine the warm start in each tensor mode until convergence.
\end{itemize}

Theoretically, we carefully characterize the optimization and statistical errors at each iteration step. The output estimate is shown to converge in a geometric rate to an estimation with minimax optimal rate in statistical error (in terms of tensor Frobenius norm). In particular, after a logarithmic number of iterations, whenever $n \gtrsim K^2(s\log(ep/s))^{\tfrac{3}{2}}$, the proposed estimator $\hat{\mT}$ achieves
\begin{equation}\label{eq:rintro-recovery-rate}
	\left\|\hat{\mT} - \mT^*\right\|_F^2 \leq C\sigma^2\frac{Ks\log (p/s)}{n}
\end{equation}
with high probability, where $s$, $K$, $p$, and $\sigma^2$ are the sparsity, rank, dimension, and noise level, respectively. 
We further establish the matching minimax lower bound to show that \eqref{eq:rintro-recovery-rate} is indeed optimal over a large class of sparse low-rank tensors. Our optimality result can be further extended to the non-sparse case (such as tensor regression \cite{zhou2013, chen2016non,rauhut2017low, li2017near}) -- to the best of our knowledge, this is the first statistical rate optimality result in both sparse and non-sparse low-rank tensor regressions.

The above theoretical analyses are non-trivial due to the non-convexity of the empirical risk function, and the need to develop some new high-order sub-Gaussian concentration inequalities. Specifically, the empirical risk function in consideration satisfies neither restricted strong convexity (RSC) condition nor sparse eigenvalue (SE) condition in general. Thus, many previous results, such as the one based on local optima analysis \cite{WLZ14, LW15,chen2016non}, are not directly applicable. Moreover, the structure of cubic-sketching tensor leads to high-order products of sub-Gaussian random variables. Thus, the matrix analysis based on Hoeffding-type or Bernstein-type concentration inequality \cite{Cai2015,chen2015exact} will lead to sub-optimal statistical rate and sample complexity. This motivates us to develop new high-order concentration inequalities and sparse tensor-spectral-type bound, i.e., Lemmas \ref{cor:concentr_high_order} and \ref{lemma:tensor_spectral} in Section \ref{sec:high-order-concentration-inequality}. These new technical results are obtained based on the careful partial truncation of high-order products of sub-Gaussian random variables and the argument of bounded $\psi_{\alpha}$-norm \cite{ALPT2011}, and may be of independent interest.

The literature on low-rank matrix estimation methods, e.g., the spectral method and nuclear norm minimization \cite{candes2009exact,keshavan2009matrix,koltchinskii2011nuclear}, is also related to this work. However, our cubic sketching model is by-no-means a simple extension from matrix estimation problems. In general, many related concepts or methods for matrix data, such as singular value decomposition, are problematic to apply in the tensor framework \cite{richard2014statistical,zhang2017tensor}. It is also found that simple unfolding or matricizing of tensors may lead to suboptimal results due to the loss of structural information \cite{Mu14}. Technically, the tensor nuclear norm is NP-hard to even approximate \cite{yuan2016tensor,yuan2017incoherent,friedland2017nuclear}, and thus the method to handle tensor low-rankness is distinct from the matrix.

The rest of the paper is organized as follows. Section \ref{sec:prelim} provides preliminaries on notation and basic knowledge of tensor. A two-stage method for symmetric tensor estimation is proposed in Section \ref{sec:symmetric}, with the corresponding theoretical analysis given in Section \ref{sec:theory-symmetric}. A concrete application to high-order interaction effect models is described in Section \ref{sec:interaction}. The non-symmetric tensor estimation model is introduced and discussed in Section \ref{sec:non-symmetrc}. Numerical analysis is provided in Section \ref{sec:simu} to support the proposed procedure and theoretical results of this paper. Section \ref{sec:discussion} discusses extensions to higher-order tensors. The proofs of technical results are given in supplementary materials.

\section{Preliminary}\label{sec:prelim}
Throughout the paper, vector, matrix, and tensor are denoted by boldface lower-case letters (e.g., $\bx, \by$), boldface upper-case letters (e.g., $\bX, \bY$), and script letters (e.g., $\cX, \cY$), respectively. For any set $A$, let $|A|$ be the cardinality. The $\diag(\bx)$ is a diagonal matrix generated by $\bx$. For two vectors $\bx$ and $\by$, $\bx\circ\by$ is the outer product. Define $\|\bx\|_q:=(|x_1|^q+\cdots+|x_p|^q)^{1/q}$. We also define the $l_0$ quasi-norm by $\|\bx\|_0 = \#\{j:x_j\neq 0\}$ and $l_{\infty}$ norm by $\max_{1\leq j \leq p}|x_j|$. Denote the set $\{1,2,\ldots, n\}$ by $[n]$. Let $\be_j$ be the canonical vectors, whose $j$-th entry equals to 1 and all other entries equal to zero. For any two sequences $\{a_n\}_{n=1}^\infty, \{b_n\}_{n=1}^\infty$, we say $a_n = \cO(b_n)$ if there exists some positive constant $C_0$ and sufficiently large $n_0$ such that $|a_n|\leq C_0b_n$ for all $n\geq n_0$. We also write $a_n\asymp b_n$ if there exists $C, c>0$ such that $ca_n\leq b_n \leq Ca_n$ for all $n\geq 1$. Additionally, $C_1,C_2,\ldots,c_1,c_2,\ldots$ are generic constants, whose actual values may be different from line to line. 

We next introduce notations and operations on the matrix. For matrices $\bA = [\ba_1,\ldots, \ba_J]\in \mathbb R^{I\times J}$ and $\bB = [\bb_1,\ldots, \bb_L]\in \mathbb R^{K\times L}$, their \emph{Kronecker product} is defined as a $(IK)$-by-$(JL)$ matrix $\bA \otimes \bB = [\ba_1\otimes \bB ~ \cdots \ba_J \otimes \bB]$, where $\ba_j\otimes \bB=(a_{j1}\bB^{\top}, \ldots, a_{jI}\bB^{\top})^{\top}$. If $\bA$ and $\bB$ have the same number of columns $J=L$, the \emph{Khatri-Rao product} is defined as $\bA \odot \bB = [\ba_1\circ \bb_1, \ba_2\circ \bb_2, \cdots, \ba_J\circ \bb_J] \in \mathbb{R}^{IK \times J}$. If the matrices $\bA$ and $\bB$ are of the same dimension, the \emph{Hadamard product} is their element-wise matrix product, such that $(\bA * \bB)_{ij} = \bA_{ij}\cdot \bB_{ij}$. For matrix $\bX = [\bx_1 \cdots \bx_n]\in \mathbb R^{m\times n}$, we also denote the vectorization $\text{vec}(\bX) = (\bx_1^{\top}, \ldots, \bx_n^{\top})\in\mathbb R^{1\times mn}$ and column-wise $\ell_2$ norms as $\text{Norm}(\bX)=(\|\bx_1\|_2,\ldots, \|\bx_n\|_2)\in\mathbb R^{1 \times n}$.

In the end, we focus on tensor notation and relevant operations. Interested readers are referred to \cite{kolda2009} for more details. Suppose $\cX\in\mathbb R^{p_1\times p_2\times p_3}$ is an order-3 tensor. Then the $(i, j, k)$-th element of $\cX$ is denoted by $[\cX]_{ijk}$. The successive tensor multiplication with vectors $\bu\in\mathbb R^{p_2}$, $\bv\in\mathbb R^{p_3}$ is denoted by $\cX\times_2 \bu \times_3\bv = \sum_{j\in[p_2], l\in[p_3]}u_j v_l \cX_{[:,j,l]}\in \mathbb R^{p_1}$. We say $\cX \in \mathbb{R}^{p_1\times p_2\times p_3}$ is \emph{rank-one} if it can be written as the outer product of three vectors, i.e., $\cX = \bx_1 \circ \bx_2 \circ \bx_3$ or $[\cX]_{ijk} = x_{1i}x_{2j}x_{3k}$ for all $i, j, k$. Here ``$\circ$" represents the vector outer product. We say $\cX$ is symmetric if $[\cX]_{ijk} = [\cX]_{ikj} = [\cX]_{jik} = [\cX]_{jki} = [\cX]_{kij} = [\cX]_{kji}$ for all $i, j ,k$. Then, $\cX$ is rank-one and symmetric if and only if it can be decomposed as $\cX = \bx\circ \bx\circ \bx$ for some vector $\bx$. 

More generally, we may decompose a tensor as the sum of rank one tensors as follows,
\begin{equation}\label{eq:CP-decomposition}
\cX = \sum_{k=1}^K \eta_k\bx_{1k}\circ \bx_{2k} \circ \bx_{3k},
\end{equation} 
where $\eta_k\in \mathbb{R}, \bx_{1k}\in \mathbb{S}^{p_1-1}, \bx_{2k}\in \mathbb{S}^{p_2-1}, \bx_{3k}\in \mathbb{S}^{p_3-1}$.
This is the so-called CANDECOMP/PARAFAC, or CP decomposition \cite{kolda2009} with CP-rank being defined as the minimum number $K$ such that \eqref{eq:CP-decomposition} holds. Then, $\{\bx_{1k}\}_{k=1}^K, \{\bx_{2k}\}_{k=1}^K, \{\bx_{3k}\}_{k=1}^K$ are called \emph{factors} along first, second and third mode. Note that factors are normalized as unit vectors to guarantee the uniqueness of decomposition, and $\boldsymbol{\eta} = \{\eta_1, \ldots, \eta_K\}$ plays an analogous role of singular values in matrix value decomposition here. 
Several tensor norms also need to be introduced. The tensor Frobenius norm and tensor spectral norm are defined respectively as
\begin{equation}\label{def:tensor_spectral_norm}
\|\cX\|_{F} = \sqrt{\sum_{i=1}^{p_1}\sum_{j=1}^{p_2}\sum_{k=1}^{p_3}\cX_{ijk}^2}, \ \|\cX\|_{op}:=\sup_{\bu\in\mathbb R^{p_1},\bv\in\mathbb R^{p_2},\bw\in\mathbb R^{p_3}}\frac{|\langle\cX, \bu\circ\bv\circ\bw\rangle|}{\|\bu\|_2\|\bv\|_2\|\bw\|_2},
\end{equation}
where $\langle\cX,\cY \rangle = \sum_{i,j,k} \cX_{ijk}\cY_{ijk}$. Clearly, $\|\cX\|_F^2 = \langle\cX,\cX\rangle$. We also consider the following sparse tensor spectral norm, 
\begin{equation}\label{def:sparse_spectral}
\|\cX\|_{s} := \sup_{\substack{\|\ba\|=\|\bb\|=\|\bc\|=1 \\ \max\{\|\ba\|_0,\|\bb\|_0, \|\bc\|_0\}\leq s}}\big|\langle\cX, \ba\circ\bb\circ\bc\rangle\big|.
\end{equation}
By definition, $\|\cX\|_s\leq \|\cX\|_{op}$. Suppose $\cX = \bx_1\circ \bx_2\circ \bx_3$ and $\cY = \by_1\circ \by_2 \circ \by_3$ are two rank-one tensors. Then it is easy to check that
$\|\cX\|_{F} = \|\bx_1\|_2\|\bx_2\|_2\|\bx_3\|_2$ and $\langle\cX, \cY\rangle = (\bx_1^{\top}\by_1)(\bx_2^{\top}\by_2)(\bx_3^{\top}\by_3)$. 

\section{Symmetric Tensor Estimation via Cubic Sketchings}\label{sec:symmetric}

In this section, we focus on the estimation of sparse and low-rank symmetric tensors,
\begin{equation}\label{eq:model}
y_i = \big\langle\mT^*, \mX_i \big\rangle+\epsilon_i, \quad \mX_i = \bx_i\circ \bx_i\circ \bx_i\in \mathbb R^{p\times p\times p}, \quad i=1,\ldots,n,
\end{equation}
where $\bx_i$ are random vectors with i.i.d. standard normal entries. As previously discussed, the tensor parameter $\mT^*$ often satisfies certain low-dimensional structures in practice, among which the factor-wise sparsity and low-rankness \cite{raskutti2015convex} commonly appear. We thus assume $\mT^*$ is CP rank-$K$ for $K \ll p$ and the corresponding factors are sparse,
\begin{equation}\label{eq:mT_low-rank}
	\begin{split}
	\mT^* = \sum_{k=1}^K \eta_k^*\bbeta_{k}^*\circ\bbeta_{k}^*\circ\bbeta_{k}^*, \quad\text{with} \quad \|\bbeta_k^*\|_2 = 1, \|\bbeta_k^*\|_0\leq s, \forall k\in [K].
	\end{split}
\end{equation}
The CP low-rankness has been widely assumed in literature for its nice scalability and simple formulation \cite{SL17,li2010tensor,li2017parsimonious}. Different from the matrix factor analysis, we do not assume the tensor factors $\bbeta_k^*$ here are orthogonal. On the other hand, since the low-rank tensor estimation is NP-hard in general \cite{hillar2013most}, we will introduce an incoherence condition in the forthcoming Condition \ref{con:incoherence_sym} to ensure that the correlation among different factors $\bbeta_k^*$ is not too strong. Such a condition has been used in recent literature on tensor data analysis \cite{anandkumar2014tensor}, compressed sensing \cite{donoho2006compressed}, matrix decomposition \cite{chandrasekaran2011rank}, and dictionary learning \cite{arora2014new}. 

Based on observations $\{y_i, \mX_i\}_{i=1}^n$, we propose to estimate $\mT^*$ via minimizing the empirical squared loss since the close-form gradient provides computational convenience,
\begin{equation}\label{eq:empirical risk minimization}
\hat{\mT} = \argmin_{\mT}\mathcal{L}(\mT)\quad \text{subject to $\mT$ is sparse and low-rank},
\end{equation} 
where
\begin{equation}\label{eq:L_s-T}
\begin{split}
\cL(\mT) & = \mathcal{L}\left(\eta_k, \bbeta_1, \ldots, \bbeta_K\right) = \frac{1}{n}\sum_{i=1}^n \left(y_i - \langle\mT, \mX_i\rangle\right)^2\\
&=  \frac{1}{n}\sum_{i=1}^n\left(y_i - \sum_{k=1}^K\eta_k\left(\bx_i^\top \bbeta_k\right)^3\right)^2. 
\end{split}
\end{equation}
Equivalently, \eqref{eq:empirical risk minimization} can be written as,
\begin{equation}\label{eqn:risk_function_sym}
\begin{split}
&\min_{\eta_k, \bbeta_k} \frac{1}{n}\sum_{i=1}^n\Big(y_i-\sum_{k=1}^K\eta_k(\bx_i^{\top}\bbeta_k)^3\Big)^2, \\
&\text{s.t.} \  \|\bbeta_k\|_2 = 1, \|\bbeta_k\|_0\leq s, \ \text{for} \ k\in[K].
\end{split}
\end{equation} 
Clearly, \eqref{eqn:risk_function_sym} is a non-convex optimization problem. To solve it, we propose a two-stage method as described in the next two subsections.

\subsection{Initialization}\label{sec:procedure-symmetric}

Due to the non-convexity of \eqref{eqn:risk_function_sym}, a straightforward implementation of many local search algorithms, such as gradient descent and alternating minimization, may easily get trapped into local optimums and result in sub-optimal statistical performance. Inspired by recent advances of spectral method (e.g., EM algorithm \cite{zhang2016spectral}, phase retrieval \cite{CLM16}, and tensor SVD \cite{zhang2017tensor}), we propose to evaluate an initial estimate $\{\eta_k^{(0)}, \bbeta_k^{(0)}\}$ via the method of moment and sparse tensor decomposition (a variant of high-order spectral method) in the following Steps 1 and 2, respectively. The pseudo-code is given in Algorithm \ref{alg:overall_symmetric}.

\textbf{Step 1: Unbiased Empirical Moment Estimator.} Construct the empirical moment-based estimator $\cT_s$,
	\begin{equation}
	\label{eqn:moment_tensor_sym}
	\begin{split}
	& \cT_s:= \frac{1}{6}\Big[\frac{1}{n}\sum_{i=1}^n y_i \bx_i\circ \bx_i\circ \bx_i \\
	&- \sum_{j=1}^p\Big(\bdm_1\circ \be_j\circ \be_j+ \be_j\circ \bdm_1 \circ \be_j+\be_j\circ \be_j\circ \bdm_1\Big)\Big],\\
	& \text{where}\ \bdm_1 := \frac{1}{n}\sum_{i=1}^n y_i\bx_i,\ \be_j \ \text{is the canonical vector}.
	\end{split}
	\end{equation}
	~\\
	\noindent  Based on Lemma \ref{lemma:moment_symmetric}, $\cT_s$ is an unbiased estimator of $\mT^\ast$. The construction of \eqref{eqn:moment_tensor_sym} is motivated by the high-order Stein's identity (\cite{janzamin2014score}; also see Theorem \ref{thm:steins_higher} for a complete statement). Intuitively speaking, based on the third-order score function of a Gaussian random vector $\bx$: $\cS_3(\bx) = \bx\circ \bx\circ \bx - \sum_{j=1}^p(\bx\circ \be_j\circ \be_j+ \be_j\circ \bx \circ \be_j+\be_j\circ \be_j\circ \bx)$, we can construct the unbiased estimator of $\mT^*$ by properly choosing a continuously differentiable function in high-order Stein's identity. See the proof of Lemma \ref{lemma:moment_symmetric} for details.

\textbf{Step 2: Sparse Tensor Decomposition.}  Based on the method of moment estimator obtained in Step 1, we further obtain good initialization for the factors $\{\eta_k^{(0)}, \bbeta_k^{(0)}\}$ via truncation and alternating rank-1 power iterations \cite{Sun2016,AGJ14},
	$$\cT_s \approx \sum_{k=1}^K \eta_k^{(0)}\bbeta_{k}^{(0)}\circ \bbeta_{k}^{(0)}\circ \bbeta_{k}^{(0)}.$$ 
Note that the tensor power iterations recover one rank-1 component per time. To identify all rank-1 components, we generate a large number of different initialization vectors, implement a clustering step, and choose the centroids as the estimates in the initialization stage. This scheme originally appears in tensor decomposition literature \cite{anandkumar2014tensor,AGJ14}, although our problem setting and proof techniques are very different. This procedure is also very different from the matrix setting since the rank-1 component in singular value decomposition is mutually orthogonal, but we do not enforce the exact orthogonality here for $\mT^*$.
	
More specifically, we first choose a large integer $M \gg K$ and generate $M$ starting vectors $\{\bb_m^{(0)}\}_{m=1}^M \in\mathbb{R}^p$ through sparse SVD as described in Algorithm \ref{alg:sparse_svd_sym}. Then for each $\bb_m^{(0)}$, we apply the following truncated power updates:
$$\tilde{\bb}_{m}^{(l+1)} = \frac{\cT_s\times_2\bb_{m}^{(l)}\times_3 \bb_{m}^{(l)}}{\|\cT_s\times_2\bb_{m}^{(l)}\times_3 \bb_{m}^{(l)}\|_2}, \quad \bb_{m}^{(l+1)} =\frac{T_d(\tilde{\bb}_m^{(l+1)})}{\|T_d(\tilde{\bb}_m^{(l+1)})\|_2},\quad l=0,\ldots,$$
where $\times_2, \times_3$ are tensor multiplication operators defined in Section \ref{sec:prelim} and $T_d(\bx)$ is a truncation operator that sets all but the largest $d$ entries in absolute values to zero for any vector $\bx$. It is noteworthy that the symmetry of $\cT_s$ implies
\begin{equation*}
    \cT_s\times_2\bb_{m}^{(l)}\times_3 \bb_{m}^{(l)}= \cT_s\times_1\bb_{m}^{(l)}\times_3 \bb_{m}^{(l)}= \cT_s\times_1\bb_{m}^{(l)}\times_2 \bb_{m}^{(l)}.
\end{equation*}
This means the multiplications along different modes are the same. We run power iterations till its convergence, and denote $\bb_m$ as the outcome. Finally, we apply $K$-means to partition $\{\bb_{m}\}_{m=1}^M$ into $K$ clusters, let the centroids of the output clusters be $\{\bbeta_k^{(0)}\}_{k=1}^K$, and calculate $\eta_k^{(0)} = \cT_s\times_1\bbeta_{k}^{(0)}\times_2 \bbeta_{k}^{(0)}\times_3 \bbeta_{k}^{(0)}$ for $k\in [K]$.
\begin{algorithm}
	\small
	\caption{Initialization in cubic sketchings}
	\begin{algorithmic}[1]
		\REQUIRE response $\{y_i\}_{i=1}^n$, sketching vector $\{\bx_i\}_{i=1}^n$, truncation level $d$, rank $K$, stopping error $\epsilon = 10^{-4}$.
		\STATE  \textbf{Step 1:} Calculate the moment-based tensor $\cT_s$ as \eqref{eqn:moment_tensor_sym}.\\
		\STATE \textbf{Step 2:} 
		\STATE \quad\textbf{For} $m = 1$ \textbf{to} $M$ \\
		\quad\quad Generate $\bb_m^{(0)}$ through Algorithm \ref{alg:sparse_svd_sym}.
		\STATE \quad\quad \textbf{Repeat}  power update:
		\begin{eqnarray*}
			&&\tilde{\bb}_{m}^{(l+1)} = \frac{\cT_s\times_2\bb_{m}^{(l)}\times_3 \bb_{m}^{(l)}}{\|\cT_s\times_2\bb_{m}^{(l)}\times_3 \bb_{m}^{(l)}\|_2}, \quad \bb_{m}^{(l+1)} =\frac{T_d(\tilde{\bb}_m^{(l+1)})}{\|T_d(\tilde{\bb}_m^{(l+1)})\|_2}, \quad l = l+1.
		\end{eqnarray*}
		\STATE \quad\quad  \textbf{Until} $\|\bb_m^{(l+1)} - \bb_m^{(l)}\|_2\leq \epsilon$.
		\STATE \quad\textbf{End for}
		\STATE \quad Perform $K$-means for $\{ \bb_{m}^{(l)}\}_{m=1}^M$. Denote the centroids of $K$ clusters by $\{\bbeta_{k}^{(0)}\}_{k=1}^K$.\\ \STATE \quad Calculate $\eta_k^{(0)} = \cT_s\times_1\bbeta_{k}^{(0)}\times_2 \bbeta_{k}^{(0)}\times_3 \bbeta_{k}^{(0)}, k\in[K].$\\
		\RETURN symmetric tensor estimator $\{\eta_k^{(0)}, \bbeta_k^{(0)}\}_{k=1}^K$
	\end{algorithmic}
	\label{alg:overall_symmetric}
\end{algorithm}

\subsection{Thresholded Gradient Descent}

After obtaining a warm start in the first stage, we propose to apply the thresholded gradient descent to iteratively refine the solution to the non-convex optimization problem \eqref{eqn:risk_function_sym}. Specifically, denote $\bX = (\bx_1,\ldots, \bx_n)\in\mathbb{R}^{p\times n}$, $\by=(y_1,\ldots, y_n)^{\top}\in \mathbb R^{n}$, $\boldsymbol{\eta} = (\eta_1, \ldots, \eta_K)^{\top}\in \mathbb R^{K}$, and $\bB = (\bbeta_1, \ldots, \bbeta_K)\in \mathbb{R}^{p\times K}$. Since $\cL(\bB, \boldsymbol{\eta}) = \cL(\mT)$, we let 
$$
\nabla_{\bB} \cL(\bB, \boldsymbol{\eta}) = (\nabla_{\bbeta_1} \cL(\bB, \boldsymbol{\eta})^{\top}, \ldots, \nabla_{\bbeta_K} \cL(\bB, \boldsymbol{\eta})^{\top}) \in \mathbb R^{1\times pK},
$$ 
be the gradient function with respect to $\bB$. Based on the detailed calculation in Lemma \ref{lemma:gradient_equivalence_sym}, $\nabla_{\bB} \cL(\bB, \boldsymbol{\eta})$ can be written as
\begin{equation}
\label{eqn:symmetric_gradient}
\nabla_{\bB} \cL(\bB, \boldsymbol{\eta}) = \frac{6}{n}[\{(\bB^{\top}\bX)^{\top}\}^3\boldsymbol{\eta} - \by]^{\top}[(\{(\bB^{\top}\bX)^{\top}\}^2\odot\boldsymbol{\eta}^{\top})^{\top} \odot \bX]^{\top},
\end{equation}
where $\{(\bB^\top \bX)^\top\}^3$ and $\{(\bB^\top \bX)^\top\}^2$ are entry-wise cubic and squared matrices of $(\bB^\top \bX)^\top$. Define $\varphi_h(x)$ as the thresholding function with a level $h$ that satisfies the following minimal assumptions:
\begin{equation}
\label{eqn:thresholding}
\left|\varphi_h(x)-x\right|\leq h, \forall x\in\mR, \ \text{and} \ \varphi_{h}(x) = 0 , \text{ when }|x|\leq h.
\end{equation}
Many widely used thresholding schemes, such as hard thresholding $H_h(x)=xI_{(|x|>h)}$, soft-thresholding $S_h(x)=\sign(x)\max(|x|-h, x)$, satisfy \eqref{eqn:thresholding}. With a slight abuse of notation, we further define the vector thresholding function as $\varphi_h(\bx) =(\varphi_h(x_1), \ldots, \varphi_h(x_p))$ for $\bx \in \mathbb{R}^p$. 

The initial estimates $\boldsymbol{\eta}^{(0)}$ and $\bB^{(0)}$ will be updated by thresholded gradient descent in two steps summarized in Algorithm	\ref{alg:overall_symmetric-gradient}. It is noteworthy that only $\bB$ is updated in Step 3, while $\boldsymbol{\eta}$ will be updated in Step 4 after finishing the update of $\bB$. 

\textbf{Step 3: Updating $\bB$ via Thresholded Gradient descent.} We update $\bB^{(t)}$ via thresholded gradient descent, 
\begin{equation}
\label{eqn:thresholed_update_sym}
\tv(\bB^{(t+1)}) = \varphi_{\frac{\mu\bh(\bB^{(t)})}{\phi}}(\tv(\bB^{(t)}) - \frac{\mu}{\phi}\nabla_{\bB} \cL(\bB^{(t)}, \boldsymbol{\eta}^{(0)})).
\end{equation}
Here, 
\begin{itemize}
	\item $\mu$ is the step size and $\phi=\sum_{i=1}^ny_i^2/n$ serves as an approximation for $(\sum_{k=1}^K\eta_k^*)^2$ (see Lemma \ref{lemma:marginal_effect});
	\item $\bh(\bB)\in \mR^{1\times K}$ is the thresholding level defined as
	\begin{equation*}
	\label{eqn:threshold_level_symmetric}
	\bh(\bB) = \sqrt{\frac{4\log np}{n^2}}[\{\{(\bB^{\top}\bX)^{\top}\}^3\boldsymbol{\eta}^{(0)} - \by\}^2]^{\top}\{\{(\bB^{\top}\bX)^{\top}\}^2\odot\boldsymbol{\eta}^{(0)\top}\}^2.
	\end{equation*}
\end{itemize}

\textbf{Step 4: Updating $\boldsymbol{\eta}$ via Normalization.} We normalize each column of $\bB^{(T)}$ and estimate the weight parameter as
	\begin{equation}\label{eq:weight-update}
	\begin{split}
	&\hat{\bB} =(\hat{\bbeta}_1, \ldots, \hat{\bbeta}_K)^{\top} = \Big(\frac{\bbeta_1^{(T)}}{\|\bbeta_1^{(T)}\|_2}, \ldots, \frac{\bbeta_K^{(T)}}{\|\bbeta_K^{(T)}\|_2}\Big),\\
	&\hat{\boldsymbol{\eta}}=  (\hat{\eta}_1, \ldots, \hat{\eta}_K)^{\top} = \Big(\eta_1^{(0)}\|\bbeta_1^{(T)}\|_2^{3}, \ldots, \eta_K^{(0)}\|\bbeta_K^{(T)}\|_2^{3}\Big)^{\top}.
	\end{split}
	\end{equation}
   The final estimator for $\mT^\ast$ is
    $$
	\hat{\mT} = \sum_{k=1}^{K}\hat{\eta}_k\hat{\bbeta}_{k}\circ \hat{\bbeta}_{k}\circ \hat{\bbeta}_{k}.
	$$
	\begin{algorithm}
	\small
	\caption{Thresholded gradient descent in cubic sketchings}
	\begin{algorithmic}[1]
		\REQUIRE response $\{y_i\}_{i=1}^n$, sketching vector $\{\bx_i\}_{i=1}^n$, step size $\mu$, rank $K$, stopping error $\epsilon = 10^{-4}$, warm-start $\{\eta_k^{(0)}, \bbeta_k^{(0)}\}_{k=1}^K$.
		\STATE \textbf{Step 3:} Let $t=0$. \\
		\STATE \quad \textbf{Repeat} thresholded gradient descent
		\STATE \begin{itemize}
			\item Compute thresholding level $\bh(\bB)$.
			\item Calculate the thresholded gradient descent update 
			$$\tv(\bB^{(t+1)}) = \varphi_{\frac{\mu\bh(\bB)}{\phi}}\Big(\tv(\bB^{(t)}) -\frac{\mu}{\phi}\nabla_{\bB}\cL(\bB^{(t)}, \boldsymbol{\eta}^{(0)})\Big),$$ 
			where $\phi=\frac{1}{n}\sum_{i=1}^ny_i^2$. The detailed form of $\nabla_{\bB}\cL(\bB, \boldsymbol{\eta}^{(0)})$ refers to  \eqref{eqn:symmetric_gradient}.
		\end{itemize}
		\STATE \quad \textbf{Until} $\|\bB^{(T+1)} - \bB^{(T)}\|_F\leq \epsilon$. \\
		
		\STATE \textbf{Step 4:} Perform column-wise normalization and update the weight as \eqref{eq:weight-update}. Construct the final estimator $\hat{\mT} = \sum_{k=1}^{K}\hat{\eta}_k\hat{\bbeta}_{k}\circ \hat{\bbeta}_{k}\circ \hat{\bbeta}_{k}$.\\
		\RETURN symmetric tensor estimator $\hat{\mT}$
	\end{algorithmic}
	\label{alg:overall_symmetric-gradient}
\end{algorithm}

\begin{algorithm}
	\caption{Sparse SVD}
	\begin{algorithmic}[1]
		\REQUIRE tensor $\cT_s$, cardinality parameter $d$.
		\STATE Compute $\tilde{\btheta}  =T_d(\btheta)$, where $\btheta\sim \cN(0,I_d)$.
		\STATE Calculate $\bu$ as the leading singular vector of $\cT_s\times_1\tilde{\btheta}$.
		\RETURN $T_d(\bu)/\|\bu\|_2$
	\end{algorithmic}
	\label{alg:sparse_svd_sym}
\end{algorithm}

\begin{Remark}[(Stochastic Thresholded Gradient descent)]
The evaluation of the gradient \eqref{eqn:symmetric_gradient} requires $\cO(npK^2)$ operations at each iteration and can be computationally intense for large $n$ or $p$. To economize the computational cost, a stochastic version of thresholded gradient descent algorithm can be easily carried out by sampling a subset of summand functions \eqref{eqn:symmetric_gradient} at each iteration. This will accelerate the procedure especially in the case of large-scale settings. See Section \ref{subsec:stochas} for details. 
\end{Remark}

\section{Theoretical Analysis}\label{sec:theory-symmetric}

In this section, we establish the geometric convergence rate in optimization error and minimax optimal rate in statistical error of the proposed symmetric tensor estimator.
\subsection{Assumptions}
We first introduce the assumptions for theoretical analysis. Conditions \ref{con:ident}-\ref{con:incoherence_sym} are on the true tensor parameter $\mT^*$ and Conditions \ref{con:noise}-\ref{con:sample} are on the measurement scheme. Specifically, the first condition ensures the model identifiability for CP-decomposition. 
\begin{Condition}[(Uniqueness of CP-decomposition)]\label{con:ident}
	The CP-decomposition in \eqref{eq:mT_low-rank} is unique in the sense that if there exists another CP-decomposition $\mT^* = \sum_{k=1}^{K'} \eta_k^{*'}\bbeta_{k}^{*'}\circ\bbeta_{k}^{*'}\circ\bbeta_{k}^{*'}$, it must have $K=K'$ and be invariant up to a permutation of $\{1,\ldots, K\}$.
\end{Condition}
For technical purposes, we introduce the following conditions to regularize the CP-decomposition of $\mT^\ast$. 
Similar assumptions were imposed in recent tensor literature, e.g., \cite{zhou2013, Sun2016} and Assumption 1.1 (A4) \citep{cai2019nonconvex}.
\begin{Condition}[(Parameter space)]\label{con:parameter}
	The CP-decomposition $\mT^* = \sum_{k=1}^K \eta_k^*\bbeta_{k}^*\circ\bbeta_{k}^*\circ\bbeta_{k}^*$ satisfies
	\begin{equation}\label{eq:condition-parameter-space}
	\|\mT^*\|_{op}\leq C\eta_{\max}^*, \quad K=\cO(s), \quad \text{and}\quad R =\eta_{\max}^*/\eta_{\min}^*\leq C'
	\end{equation}
	for some absolute constants $C,C'$, where $\eta_{\min}^*=\min_{k}\eta_k^*$ and $\eta_{\max}^*=\max_{k}\eta_k^*$. Recall that $s$ is the sparsity of $\bbeta_k^*$.
\end{Condition}
\begin{Remark}
In Condition \ref{con:parameter}, $R$ plays a similar role as a ``condition number.'' This assumption means that the tensor is ``well-conditioned," i.e., each rank-1 component is roughly of the same size. 
\end{Remark} 

As shown in the seminal work of \cite{hillar2013most}, the estimation of low-rank tensors can be NP-hard in general. Hence, we impose the following incoherence condition.
\begin{Condition}[(Parameter incoherence)]
	\label{con:incoherence_sym}
	The true tensor components are incoherent such that 
	\begin{eqnarray*}
		\Gamma:=\max_{1\leq k_1\neq k_2\leq K}|\langle\bbeta_{k_1}^*, \bbeta_{k_2}^*\rangle|\leq \min\{C^{''}K^{-\tfrac{3}{4}}R^{-1}, s^{-\tfrac{1}{2}}\},
	\end{eqnarray*}
where $R$ is the singular value ratio defined in \eqref{eq:condition-parameter-space} and $C^{''}$ is some small constant. 
\end{Condition}

\begin{Remark}
The preceding incoherence condition has been widely used in different scenarios in recent high-dimensional research, such as tensor decomposition \cite{Sun2016,AGJ14}, compressed sensing \cite{donoho2006compressed}, matrix decomposition \cite{chandrasekaran2011rank}, and dictionary learning \cite{arora2014new}. It can be also viewed as a relaxation of orthogonality: if $\{\bbeta_1^*, \ldots, \bbeta_K^*\}$ are mutually orthogonal, $\Gamma$ equals zero. We can show from both theory (Lemma \ref{lemma:proof_incoherence} in the supplementary materials) and simulation (Section \ref{sec:simu}) that the low-rank tensor $\mT^*$ induced by \eqref{eq:mT_low-rank} satisfies the incoherence condition with high probability, if the component vectors $\bbeta_k^*$ are randomly generated, say from Gaussian distribution.
\end{Remark}

We also introduce the following conditions on noise distribution. 
\begin{Condition}[(Sub-exponential noise)]\label{con:noise}
	The noise $\{\epsilon_i\}_{i=1}^n$ are i.i.d. randomly generated with mean 0 and variance $\sigma^2$ satisfying $0<\sigma<C\sum_{k=1}^K\eta_k^*$. $(\epsilon_i/\sigma)$ is sub-exponential distributed, i.e., there exists constant $C_{\epsilon}>0$ such that $\|(\epsilon_i/\sigma)\|_{\psi_1} := \sup_{p\geq 1} p^{-1}(\mathbb{E}|\epsilon_i/\sigma|^p)^{1/p}\leq C_{\epsilon}$, and is independent of $\{\mX_i\}_{i=1}^n$.
\end{Condition}

The sample complexity condition	is crucial for our algorithm especially in the initialization stage. Ignoring any polylog factors, Condition \ref{con:sample} is even weaker than the sparse matrix estimation case $(n\gtrsim s^2)$ in \cite{CLM16}.
\begin{Condition}[(Sample complexity)] \label{con:sample}
	\begin{equation*}
	n \geq C^{'''}K^2(s\log(ep/s))^{\tfrac{3}{2}}\log^4n.
	\end{equation*}
\end{Condition}

\subsection{Main Theoretical Results}\label{sec:main-theoretical}

Our main Theorem \ref{thm:symmetric_main} shows that based on a proper initializer, the output of the proposed procedure can achieve optimal estimation error rate after a sufficient number of iterations. Here, we define the contraction parameter 
$$0<\kappa = 1-32\mu K^{-2}R^{-\tfrac{8}{3}}<1$$ 
and also denote $\cE_1 = 4K\eta_{\max}^{*\tfrac{2}{3}}\varepsilon_0^2$ and $\cE_2= C_0\eta_{\min}^{*-\tfrac{4}{3}}/16$ for some $C_0>0.$
\begin{Theorem}[(Statistical and Optimization Errors)]
	\label{thm:symmetric_main}
    Suppose Conditions \ref{con:incoherence_sym}-\ref{con:sample} hold, $|\supp(\bbeta_k^{(0)})|\lesssim s$, and the initial estimator $\{\bbeta_k^{(0)}, \eta_k^{(0)}\}_{k=1}^K$ satisfy
    \begin{equation}\label{ineq:initial}
    \max_{1\leq k\leq K}\Big\{\big\|\bbeta_k^{(0)} - \bbeta_k^*\big\|_2, |\eta_k^{(0)} - \eta_k^*|\Big\}\lesssim K^{-1}
\end{equation}
with probability at least $1-\cO(1/n)$. Assume the step size $\mu\leq \mu_0$, where $\mu_0$ is defined in \eqref{eqn:mu_0}. Then, the output of the thresholded gradient descent update in \eqref{eqn:thresholed_update_sym} satisfies:
	\begin{itemize}
		\item For any $t=0,1,2,\ldots$, the factor-wise estimator satisfies
		\begin{eqnarray}\label{ineqn:error_bound}
			\sum_{k=1}^{K}\Big\|\sqrt[3]{\eta_k^{\scriptscriptstyle{(0)}}}\bbeta_k^{(t+1)} -\sqrt[3]{\eta_k^*}\bbeta_k^{*} \Big\|_2^2
			\leq \cE_1\kappa^t+\cE_2\frac{\sigma^2s\log p}{n}
		\end{eqnarray}
		with probability at least $1- \cO(tKs/n)$.
		\item When the total number of iterations is no smaller than
		\begin{eqnarray}\label{eqn:T^*}
			T^*=\Big(\log (\frac{n}{\sigma^2s \log p} \vee 1) + \log \frac{\cE_1}{\cE_2}\Big)/\log \kappa^{-1},
\end{eqnarray}
there exists a constant $C_1$ (independent of $K, s, p, n, \sigma^2$) such that the final estimator $\hat{\mT} = \sum_{k=1}^K\eta_k^{(0)}\bbeta_k^{(T^*)}\circ\bbeta_k^{(T^*)}\circ\bbeta_k^{(T^*)}$ satisfies
\begin{eqnarray}\label{ineqn:final_rate}
\Big\|\hat{\mT}-\mT^*\Big\|_F^2\leq \frac{C_1\sigma^2Ks\log p}{n}
\end{eqnarray}
with probability at least $1-\cO(T^*Ks/n)$.
	\end{itemize}
\end{Theorem}
\begin{Remark}
The error bound \eqref{ineqn:error_bound} can be decomposed into an optimization error $\cE_1\kappa^t$ (which decays with a geometric rate as iterations) and a statistical error $\cE_2\frac{\sigma^2s\log p}{n}$ (which does not decay as iterations). 
In the special case that $\sigma =0$, $\hat{\mT}$ exactly recover $\mT^\ast$ with high probability.
\end{Remark}


The next theorem shows that Steps 1 and 2 of Algorithm \ref{alg:overall_symmetric} provides a good initializer required in Theorem \ref{thm:symmetric_main}. 
\begin{Theorem}[(Initialization Error)]
	\label{thm:initial_symmetric}
    Recall $\Gamma = \max_{1\leq k_1\neq k_2\leq K}|\langle\bbeta_{k_1}^*, \bbeta_{k_2}^*\rangle|$. Suppose the number of initializations $L \geq K^{C_3\gamma^{-4}}$, where $\gamma$ is a constant defined in \eqref{def:gamma}. Given that Conditions \ref{con:ident}-\ref{con:noise} hold, the initial estimator obtained from Steps 1-2 with a truncation level $s \leq d\leq Cs$ satisfies
	\begin{equation}\label{eqn:initial_rate}
	\max_{1\leq k \leq K}\Big\{\|\bbeta_{k}^{(0)}-\bbeta_{k}^*\|_2, |\eta_k^{(0)}-\eta_k^*|\Big\} \leq C_2KR\delta_{n,p,s} + \sqrt{K}\Gamma^2
	\end{equation}
	and 
	$$|\supp(\bbeta_k^{(0)})|\lesssim s
	$$ 
	with probability at least $1-5/n$, where 
	\begin{equation}\label{eqn:initial_stat_rate}
	\delta_{n,p,s} = (\log n)^3\Big(\sqrt{\frac{s^3\log ^3(ep/s)}{n^2}} + \sqrt{\frac{s\log (ep/s)}{n}}\Big).
	\end{equation}
	Moreover, if the sample complexity condition \ref{con:sample} holds, then the above bound satisfies \eqref{ineq:initial}.
\end{Theorem}
\begin{Remark}[(Interpretation of initialization error)]\label{rm:initialization error}
The upper bound of \eqref{eqn:initial_rate} consists of two terms that correspond to the approximation error of $\cT_s$ to $\mT^\ast$ and the incoherence among $\bbeta_{k}^\ast$'s, respectively. Especially, the former converges to zero as $n$ grows while the latter does not. 

\end{Remark}
The proof of Theorems \ref{thm:symmetric_main} and \ref{thm:initial_symmetric} are postponed to Section \ref{subsec:proof_ini}-\ref{subsec:gradient-update-proof} in the supplementary materials. The combination of Theorems \ref{thm:symmetric_main} and \ref{thm:initial_symmetric} immediately yields the following upper bound for the final estimator, which is one main result of this paper.
\begin{Theorem}[(Upper Bound)]
	\label{thm:symmetric_upper_bound}
Suppose Conditions \ref{con:ident} -- \ref{con:sample} hold, $s\leq d\leq Cs$. After $T^*$ iterations, there exists a constant $C_1$ not depending on $K, s, p, n, \sigma^2$, such that the proposed procedure yields
\begin{eqnarray}\label{ineqn:final_rate_main}
\Big\|\hat{\mT}-\mT^*\Big\|_F^2\leq \frac{C_1\sigma^2Ks\log p}{n}
\end{eqnarray}
with probability at least $1-\cO(T^*Ks/n)$, where $T^*$ is defined in \eqref{eqn:T^*}.
\end{Theorem}


The above upper bound turns out to match the minimax lower bound for a large class of sparse and low-rank tensors.
\begin{Theorem}[(Lower Bound)]
	\label{th:lower_bound}
	Consider the following class of sparse and low-rank tensors,
	\begin{equation}\label{eq:mathcal_F_symmetric}
	\mathcal{F}_{p, K, s} = \left\{\mT: \begin{array}{l}
	\mT = \sum_{k=1}^K\eta_k \bbeta_k\circ \bbeta_k\circ \bbeta_k, \|\bbeta_k\|_0 \leq s, \text{ for } k \in[K],\\ 
	\mT \ \text{satisfies Conditions \ref{con:ident},  \ref{con:parameter}, and \ref{con:incoherence_sym}}
	\end{array}\right\}.
	\end{equation}
	Suppose that $\{\mX_i\}_{i=1}^n$ are i.i.d standard normal cubic sketchings with i.i.d. $N(0, \sigma^2)$ noise in  \eqref{eq:model}, $p\geq 20s$, and $s\geq 4$. We have the following lower bound result,
	\begin{equation*}
		\inf_{\tilde{\mT}}\sup_{\mT \in \mathcal{F}_{p, K, s}} \mathbb{E} \left\|\tilde{\mT} - \mT\right\|_F^2\geq c\sigma^2\frac{Ks \log(ep/s)}{n}.
	\end{equation*}
\end{Theorem}
The proof of Theorem \ref{th:lower_bound} is deferred to Section \ref{proof:thm_minimax} in the supplementary materials. Combining Theorems \ref{thm:symmetric_upper_bound} and \ref{th:lower_bound}, we immediately obtain the following minimax-optimal rate for sparse and low-rank tensor estimation with cubic sketchings when $\log p \asymp \log (p /s)$:
	\begin{equation}\label{eqn:optimal_rate}
	\inf_{\tilde{\mT}}\sup_{\mT^* \in \mathcal{F}_{p, K, s}}\mathbb{E}\left\|\tilde{\mT} - \mT^*\right\|_F^2 \asymp \sigma^2\frac{Ks\log (ep/s)}{n}.
	\end{equation}
The rate in \eqref{eqn:optimal_rate} sheds light upon the effect of dimension $p$, noise level $\sigma^2$, sparsity $s$, sample size $n$ and rank $K$ to the estimation performance.

\begin{Remark}
Recently, Li, Haupt, and Woodruff \cite{li2017near} studied the optimal sketching for the low-rank tensor regression and gave an near-optimal sketching complexity with a sharp $(1+\varepsilon)$-worse-case error bound. Different from the framework of \cite{li2017near} that focuses on a deterministic setting, we study a probabilistic model with random observation noises, propose a new algorithm, and studied the minimax optimal rate of estimation errors. In addition, \cite{SL17,raskutti2015convex,chen2016non} considered different types of convex/non-convex algorithms for low-rank tensor regression with statistical assumptions. To our best knowledge, we are the first to achieve an optimal rate in estimation error based on polynomial-time algorithms for the tensor regression problem.

\end{Remark}

\begin{Remark}[(Non-sparse low-rank tensor estimation via cubic-sketchings)]
	When the low-rank tensor $\mT^\ast$ is not necessarily sparse, i.e.,
	\begin{equation*}
	\mT^* \in \mathcal{F}_{p, K} = \left\{\mT: \begin{array}{l}
	\mT = \sum_{k=1}^K\eta_k \bbeta_k\circ \bbeta_k\circ \bbeta_k, \text{ for } k \in[K],\\ 
	\mT \ \text{satisfies Conditions \ref{con:ident}, \ref{con:parameter}, and \ref{con:incoherence_sym}}
	\end{array}\right\},
	\end{equation*}
	we can apply the proposed procedure with all the truncation/thresholding steps removed. If $n \geq \cO (p^{3/2})$, we can use similar arguments of Theorems \ref{thm:symmetric_main}-\ref{thm:symmetric_upper_bound} to show that the estimator $\hat{\mT}'$ satisfies
	\begin{equation}\label{eq:non-sparse-rate}
		\left\|\hat{\mT}' - \mT^*\right\|_F^2\lesssim \frac{\sigma^2 Kp}{n}
	\end{equation}
	for any $\mT^* \in \mathcal{F}_{p, K}$ with high probability. Furthermore, similar arguments of Theorem \ref{th:lower_bound} imply that the rate in \eqref{eq:non-sparse-rate} is minimax optimal. 
\end{Remark}
	
\begin{Remark}[(Comparison with existing matrix results)]\label{rm:comparison-matrix-setting}
Our cubic sketching tensor results are far more than extensions of the existing matrix ones. For example, \cite{Cai2015,chen2015exact} studied the low-rank matrix recovery via rank-1 projections: $y_i = \bx_i^\top \bT \bx_i + \epsilon_i$ and proposed the convex nuclear norm minimization methods. The theoretical properties of their estimate are analyzed under a $\ell_1/\ell_2$-RIP or Restricted Uniform Boundedness condition (RUB). However, the tensor nuclear norm is computationally infeasible and one can check that our cubic sketching framework does not satisfy RIP or RUB conditions in general following the arguments in \cite{CLM16,CLS15}. Thus, these previous results cannot be directly applied.
	
In addition, the analysis of gradient updates for the tensor case is significantly more complicated than the matrix case. First, it requires high-order concentration inequalities for the tensor case since the cubic-sketching tensor leads to high-order products of sub-Gaussian random variables (see Section \ref{sec:high-order-concentration-inequality} for details). The necessity of high-order expansions in the analysis of gradient updates for the tensor case also significantly increases the hardness of the problem. To ensure the geometric convergence, we need much more subtle analysis comparing to the ones in the matrix case \cite{CLS15}.
\end{Remark}

\subsection{Key Lemmas: High-order Concentration Inequalities}\label{sec:high-order-concentration-inequality}

As mentioned earlier, one major challenge for theoretical analysis of cubic sketching is to handle heavy tails of high-order Gaussian moments. One can only handle up-to second moments of sub-Gaussian random variables by directly applying the Hoeffding's or Bernstein's concentration inequalities. Therefore, we need to develop the following high-order concentration inequalities as technical tools: Lemma \ref{cor:concentr_high_order} characterizes the tail bounds for the sum of sub-Gaussian products, and Lemma \ref{lemma:tensor_spectral} provides the concentration inequalities for Gaussian cubic sketchings. The proofs of Lemmas \ref{cor:concentr_high_order} and \ref{lemma:tensor_spectral} are given in Section \ref{subsec:concentration}.

\begin{Lemma}[(Concentration inequality for sum of sub-Gaussian products)]\label{cor:concentr_high_order}
	Suppose $\bX_i= (\bx_{1i}^{\top}, \ldots, \bx_{mi}^{\top})^{\top}\in \mathbb R^{m\times p}$, $i\in[n]$ are $n$ i.i.d  random matrices. Here, suppose $\bx_{ij}$, the $j$-th row of $\bX_i$, is an isotropic sub-Gaussian vector, i.e., $\mathbb{E} \bx_{ij}=0$ and $\Cov(\bx_{ij}) = I$. Then for any vectors $\ba = (a_1\ldots,a_n)\in\mathbb R^n$, $\{\bbeta_j\}_{j=1}^m\subseteq \mathbb R^{p}$, and $0<\delta<1$, we have
	\begin{equation*}
	\begin{split}
		&\Big|\sum_{i=1}^na_i\prod_{j=1}^m(\bx_{ij}^{\top}\bbeta_j)-\mathbb E\Big(\sum_{i=1}^na_i\prod_{j=1}^m(\bx_{ij}^{\top}\bbeta_j)\Big)\Big|\\
		&\leq C\prod_{j=1}^m\|\bbeta_j\|_2\Big(\|\ba\|_{\infty}(\log \delta^{-1})^{m/2}+\|\ba\|_2(\log \delta^{-1})^{1/2}\Big)
		\end{split}
	\end{equation*}
	with probability at least $1-\delta$ for some constant $C$. 
\end{Lemma}
	Note that in Lemma \ref{cor:concentr_high_order}, each $\bX_i$ does not necessarily have independent entries, even though $\{\bX_i\}_{i=1}^n$ are independent matrices. Building on Lemma \ref{cor:concentr_high_order}, Lemma \ref{lemma:tensor_spectral} provides a generic spectral-type concentration inequality that can be used to quantify the approximation error of $\cT_s$ introduced in Step 1 of the proposed procedure. 
\begin{Lemma}[(Concentration inequality for Gaussian cubic sketchings)]
	\label{lemma:tensor_spectral}
	Suppose $\{\bx_{1i}\}_{i=1}^n \overset{iid}{\sim} \cN(0, \bI_{p_1})$, $\{\bx_{2i}\}_{i=1}^n\overset{iid}{\sim} \cN(0, \bI_{p_2})$, $\{\bx_{3i}\}_{i=1}^n \overset{iid}{\sim} \cN(0, \bI_{p_3})$, $\bbeta_1\in \mathbb R^{p_1}$, $\bbeta_2\in \mathbb R^{p_2}$, $\bbeta_3\in \mathbb R^{p_3}$ are fixed vectors. 
	\begin{itemize}
		\item Define $M_{\text{nsy}} = \frac{1}{n}\sum_{i=1}^n\langle \bx_{1i}\circ \bx_{2i}\circ \bx_{3i}, \bbeta_1\circ \bbeta_2\circ \bbeta_3\big\rangle\bx_{1i}\circ\bx_{2i}\circ\bx_{3i}$. Then $\mathbb{E}(M_{\text{nsy}}) = \bbeta_1\circ\bbeta_2 \circ\bbeta_3$ and
		\begin{equation*}
		\begin{split}
			&\Big\|M_{\text{nsy}} - \mathbb E(M_{\text{nsy}})\Big\|_s\\
			\leq & C(\log n)^3\Big(\sqrt{\frac{s^3\log^3 (ep/s)}{n^2}} + \sqrt{\frac{s\log (ep/s)}{n}}\Big)\|\bbeta_1\|_2\|\bbeta_2\|_2\|\bbeta_3\|_2
			\end{split}
		\end{equation*}
		with probability at least $1-10/n^3-1/p$.
		\item Define $M_{\text{sym}} = \frac{1}{n}\sum_{i=1}^n\langle \bx_{1i}\circ \bx_{1i}\circ \bx_{1i}, \bbeta_1\circ \bbeta_1\circ \bbeta_1\big\rangle\bx_{1i}\circ\bx_{1i}\circ\bx_{1i}$. Then $\mathbb E(M_{\text{sym}})=6\bbeta_1\circ \bbeta_1 \circ \bbeta_1+3\sum_{m=1}^p(\bbeta_1\circ \be_m\circ \be_m + \be_m\circ \bbeta_1\circ\be_m + \be_m\circ \be_m\circ \bbeta_1)$ and
		\begin{equation*}
		\begin{split}
		&	\Big\|M_{\text{sym}} - \mathbb E(M_{\text{sym}})\Big\|_s\\
		\leq & C(\log n)^3\Big(\sqrt{\frac{s^3\log^3 (ep/s)}{n^2}} + \sqrt{\frac{s\log (ep/s)}{n}}\Big)\|\bbeta_1\|_2^3
			\end{split}
		\end{equation*}
		with probability at least $1-10/n^3-1/p$.
	\end{itemize}	
	Here, $C$ is an absolute constant and $\|\cdot\|_{s}$ is the sparse tensor spectral norm defined in \eqref{def:sparse_spectral}.
\end{Lemma}

\section{Application to High-order Interaction Effect Models}\label{sec:interaction}

In this section, we study the high-order interaction effect model in the cubic sketching framework. Specifically, we consider the following three-way interaction model
\begin{equation}\label{eq:interaction_model}
y_l  = \xi_0 +  \sum_{i=1}^p \xi_i z_{li}  + \sum_{i, j = 1}^{p} \gamma_{ij} z_{li}z_{lj} + \sum_{i,j,k = 1}^p \eta_{ijk} z_{li}z_{lj}z_{lk}+\epsilon_l, \quad l=1,\ldots, n.
\end{equation}
Here $\bxi$, $\bgamma$, and $\boldsymbol{\eta}$ are coefficients for the main effect, pairwise interaction, and triple-wise interaction, respectively. More importantly, \eqref{eq:interaction_model} can be reformulated into the following tensor form (also see the left panel of Figure \ref{fig:illurstration})
\begin{equation}\label{eq:skeching}
y_l = \langle \cB, \bx_l \circ \bx_l \circ \bx_l\rangle + \epsilon_l,\quad l=1,\ldots, n,
\end{equation}
where $\bx_l = (1, \bz_l^\top)^\top \in \mathbb{R}^{p+1}$ and $\cB\in \mathbb{R}^{(p+1)\times (p+1)\times (p+1)}$ is a tensor parameter corresponding to coefficients in the following way:
\begin{equation}
\left\{\begin{array}{l}
\cB_{[0, 0, 0]} = \xi_0,\\ 
\cB_{[1:p, 1:p, 1:p]} = (\eta_{ijk})_{1\leq i,j,k\leq p}, \\
\cB_{[0, 1:p, 1:p]} = \cB_{[1:p, 0, 1:p]} = \cB_{[1:p, 1:p, 0]} = (\gamma_{ij}/3)_{1\leq i,j\leq p}, \\
\cB_{[0, 0, 1:p]} = \cB_{[0, 1:p, 0]} = \cB_{[1:p, 0, 0]} = (\xi_i/3)_{1\leq i \leq p}.
\end{array}\right.
\end{equation}

We provide the following justification for assuming the tensorized coefficient $\cB$ is low-rank and sparse. First, in modern applications, such as the biomedical research \cite{Hung2016}, the response is often driven by a small portion of coefficients and a small number of factors, leading to a highly entry-wise sparse and low-rank $\cB$. Second, \cite{SK12} suggested that it is suitable to model entry-wise sparse and low-enough rank tensors as arising from sparse loadings. 
Therefore, we assume $\cB$ is CP rank-$K$ with $s$-sparse factors:
$$\cB = \sum_{k=1}^K \eta_k\bbeta_k\circ\bbeta_k\circ\bbeta_k, \quad \|\bbeta_k\|_0 \leq s,$$
where $K, s\ll p$. Then the number of parameters in \eqref{eqn:tensor_interaction}, $K(p+1)$, is significantly smaller than $(p+1)^3$, the total number of parameters in the original three-way interaction effect model \eqref{eq:interaction_model}, which makes the consistent estimation of $\cB$ possible in the high-dimensional case. In this case, \eqref{eq:skeching} can be written as 
\begin{equation}\label{eqn:tensor_interaction}
\begin{split}
	& y_l=\Big\langle \sum_{k=1}^K\eta_k\bbeta_k\circ\bbeta_k\circ\bbeta_k,\bx_l\circ\bx_l\circ\bx_l\Big \rangle + \epsilon_l,\quad l=1,\ldots, n,\\
	& \text{where}\quad \|\bbeta_k\|_2 =1, \quad \|\bbeta_k\|_0 \leq s, \quad k\in [K].
\end{split}
\end{equation}

By assuming $\bz_l \overset{iid}{\sim}N_{p}(0, \bI_{p})$, the high-order interaction effect model \eqref{eq:skeching} reduces to the symmetric tensor estimation model \eqref{eq:model}, except one slight difference that the first coordinate of $\bx_l$, i.e., the intercept, is always 1. To accommodate this difference, we only need to adjust the initial unbiased estimate in the above two-step procedure. Let 
\begin{equation}\label{eqn:moment_inter}
\begin{split}
\cT_s = & \frac{1}{6n}\sum_{l=1}^n y_l \bx_l \circ \bx_l \circ \bx_l- \frac{1}{6}\sum_{j=1}^p(\ba\circ \be_j\circ \be_j+ \be_j\circ \ba \circ \be_j+\be_j\circ \be_j\circ \ba), \\
& \text{where } \quad \ba = \frac{1}{n}\sum_{l=1}^n y_l \bx_l.
\end{split}
\end{equation}
Then we construct the empirical moment-based initial tensor $\cT_{s'}$ as
\begin{itemize}
	\item For $i,j,k\neq 0$, $\cT_{s'[i,j,k]} = \cT_{s[i,j,k]}$, $\cT_{s'[i,j,0]} = \cT_{s[i,j,0]}, \cT_{s'[0,j,k]} = \cT_{s[0,j,k]}$, and $\cT_{s'[i,0,k]} = \cT_{s[i,0,k]}$.
	\item For $i\neq 0$, $\cT_{s'[0,0,i]} = \cT_{s'[ 0, i, 0]} = \cT_{s'[i,0,0]} = \frac{1}{3}\cT_{s[0,0,i]} - \frac{1}{6}(\sum_{k=1}^p\cT_{s[k,k,i]}-(p+2)a_i)$.
	\item $\cT_{s'[0,0,0]} = \frac{1}{2p-2}(\sum_{k=1}^p\cT_{s[0,k,k]} - (p+2)\cT_{s[0,0,0]})$. 
\end{itemize}
Lemma \ref{lemma:moment_interaction} shows that $\cT_{s'}$ is an unbiased estimator for $\mathcal{B}$.

The theoretical results in Section \ref{sec:theory-symmetric} imply the following upper and lower bounds for the three-way interaction effect estimation.
\begin{Corollary}
	Suppose $\bz_1,\ldots, \bz_n$ are i.i.d. standard Gaussian random vectors and $\mathcal{B}$ satisfies Conditions \ref{con:ident}, \ref{con:parameter} and \ref{con:incoherence_sym}. The output, denoted as $\hat \cB$, from the proposed Algorithms \ref{alg:overall_symmetric} and \ref{alg:overall_symmetric-gradient} based on $\cT_{s'}$ satisfies
	\begin{equation}
	\left\|\hat{\cB} - \cB\right\|_F^2 \leq C\frac{\sigma^2Ks\log p}{n}
	\end{equation}
	with high probability. On the other hand, considering the following class of $\mathcal{B}$,
	\begin{equation*}
	\mathcal{F}_{p+1, K, s} = \left\{\mathcal{B}: \begin{array}{l}
	\mathcal{B} = \sum_{k=1}^K\eta_k \bbeta_k\circ \bbeta_k\circ \bbeta_k, \|\bbeta_k\|_0 \leq s, \text{ for } k \in[K],\\ 
	\mathcal{B} \ \text{satisfies Conditions \ref{con:ident},  \ref{con:parameter}, and \ref{con:incoherence_sym}},
	\end{array}\right\}.
	\end{equation*}
	Then the following lower bound holds,
	\begin{equation*}
	\inf_{\hat{\mathcal{B}}}\sup_{\mathcal{B}\in \mathcal{F}_{p+1, K, s}}\mathbb{E}\left\|\hat{\cB} - \cB\right\|_F^2 \geq C\frac{\sigma^2Ks\log p}{n}.
	\end{equation*}
\end{Corollary}

\section{Non-symmetric Tensor Estimation Model}\label{sec:non-symmetrc}

In this section, we extend the previous results to the non-symmetric tensor case. Specifically, we have $\mT^*\in \mathbb R^{p_1\times p_2\times p_3}$ and
\begin{equation}\label{eq:model_asymmetric}
\begin{split}
y_i = \langle\mT^*, \mX_i \rangle+\epsilon_i, \ \mX_i = \bu_i\circ \bv_i\circ \bw_i, \ i\in[n],
\end{split}
\end{equation}
where $\bu_i\in \mathbb{R}^{p_1}, \bv_i\in \mathbb{R}^{p_2}, \bw_i\in \mathbb{R}^{p_3}$ are random vectors with i.i.d. standard normal entries. Again, we assume $\mT^*$ is sparse and low-rank in a similar sense that
\begin{equation}\label{eq:low_rank_nonsym}
\begin{split}
& \mT^* = \sum_{k=1}^K \eta_k^*\bbeta_{1k}^*\circ\bbeta_{2k}^*\circ\bbeta_{3k}^*,\\  
& \|\bbeta_{1k}^*\|_2 = \|\bbeta_{2k}^*\|_2 = \|\bbeta_{3k}^*\|_2 = 1,\quad \max\{\|\bbeta_{1k}^*\|_0, \|\bbeta_{2k}^*\|_0, \|\bbeta_{3k}^*\|_0\}\leq s.
\end{split}
\end{equation}
Denote
\begin{itemize}
	\item $\bB_1 = (\bbeta_{11}, \cdots,\bbeta_{1K})$, $\bB_2 = (\bbeta_{21}, \cdots,\bbeta_{2K})$, $\bB_3 = (\bbeta_{31}, \cdots,\bbeta_{3K})$,
	\item $\bU=(\bu_1, \ldots,\bu_n),  \ \bV=(\bv_1, \ldots,\bv_n), \ \bW=(\bw_1, \ldots,\bw_n)$, \\
	 $\boldsymbol{\eta} = (\eta_1, \ldots, \eta_k)^{\top}, \by = (y_1, \ldots, y_n)^{\top}.$
\end{itemize} 
Then, the empirical risk function can be written compactly as
\begin{equation}\label{def:risk_function_nym}
\cL(\bB_1, \bB_2, \bB_3, \boldsymbol{\eta} ) 
= \frac{1}{n}\Big\|(\bU^{\top}\bB_1) * (\bV^{\top}\bB_2) * (\bW^{\top}\bB_3)\cdot\boldsymbol{\eta} - \by\Big\|_2^2.
\end{equation}
Since \eqref{def:risk_function_nym} is non-convex but fortunately tri-convex in terms of $\bB_1$, $\bB_2$, and $\bB_3$, we develop a block-wise thresholded gradient descent algorithm as detailed below. The complete algorithm is deferred to Section \ref{alg:overall_nonsymmetric} in the supplementary materials.
\begin{description}
	\item[Step 1: (Method of Tensor Moments)] Construct the empirical moment-based estimator
	\begin{eqnarray}\label{eq:moment_asy}
	\cT:=\frac{1}{n}\sum_{i=1}^n y_i \bu_i\circ \bv_i\circ \bw_i\in \mathbb R^{p_1\times p_2\times p_3}
	\end{eqnarray}	 
	to which sparse tensor decomposition is applied for initialization. 
	\item[Step 2: (Block-wise Gradient Descent)] Lemma \ref{lemma:gradient_equivalence_aym} shows that the gradient function for \eqref{def:risk_function_nym} with respect to $\bB_1$ can be written as
	\begin{equation}
	\label{eqn:gradient_nonsym}
		\nabla_{\bB_1}\cL(\bB_1, \bB_2, \bB_3, \boldsymbol{\eta}) = \bD^{\top}(\bC_1^{\top}\odot\bU)^{\top}\in \mathbb R^{1\times (p_1K)},
	\end{equation}
	where $\bD = (\bB_1^{\top}\bU)^{\top}*(\bB_2^{\top}\bV)^{\top}*(\bB_3^{\top}\bW)^{\top}\boldsymbol{\eta} - \by$ and $\bC_1 = (\bB_2^{\top}\bV)^{\top}*(\bB_3^{\top}\bW)^{\top}\odot \boldsymbol{\eta}^{\top}$.
  	For $t=1,\ldots, T$, we fix $\bB_2^{(t)}, \bB_3^{(t)}$ and update $\bB_1^{(t+1)}$ via block-wise thresholded gradient descent,
	\begin{eqnarray*}
		\tv(\bB_1^{(t+1)}) = \varphi_{\tfrac{\mu h(\bB_1^{(t)})}{\phi}}\Big(\tv(\bB_1^{(t)}) - \frac{\mu}{\phi}\nabla_{\bB_1}\cL(\bB_1^{(t)}, \bB_2^{(t)}, \bB_3^{(t)}, \boldsymbol{\eta})\Big), 
	\end{eqnarray*}
	where $\phi= \sum_{i=1}^{n}y_i^2/n$, $\mu$ is the step size, and $\bh(\bB) = \sqrt{\frac{4\log np}{n^2}\{\bD^{2}\}^{\top}\{\bC^2\}}$. The updates of $\bB_2, \bB_3$ are similar.
\end{description}

The theoretical analysis for the non-symmetric case is different from the symmetric one in two folds. First, the non-symmetric cubic sketching tensor is formed by three Gaussian vectors rather than one, which leads to many differences in the calculation of high-order moments. Second, the CP-decomposition of non-symmetric tensor $\mT^*$ \eqref{eq:low_rank_nonsym} forms a tri-convex optimization. At this point, the standard convex analysis for vanilla gradient descent \cite{Bubeck15} could be applied given a proper initialization. 

With the regularity conditions detailed in Section \ref{alg:overall_nonsymmetric}, we present the theoretical results for non-symmetric tensor estimation as follows.
\begin{Theorem}[(Upper Bound)]
	\label{thm:main_nym}
	Suppose Conditions \ref{con:ident_nym} -- \ref{con:noise_non} hold and $n\gtrsim (s\log (p_0/s))^{3/2}$, where $p_0 = \max\{p_1, p_2, p_3\}$. For any $t=0,1,2,\ldots$, the output of Algorithm \ref{alg:overall_nonsymmetric} satisfies
		\begin{eqnarray*}
		\sum_{k=1}^{K}\sum_{j=1}^3\Big\|\sqrt[3]{\eta_{k}}\bbeta_{jk}^{(t+1)} -\sqrt[3]{\eta_k^*}\bbeta_{jk}^{*} \Big\|_2^2 \leq \cO_p\Big(\kappa^t+\frac{\sigma^2s\log p_0}{n}\Big)
		\end{eqnarray*}
	for some $0<\kappa<1$. When the total number of iterations is no smaller than $\log (\frac{n}{\sigma^2s \log p_0} \vee 1) /\log \kappa^{-1}$, the final estimator $\hat{\mT}$ satisfies
		\begin{eqnarray*}
			\Big\|\hat{\mT}-\mT^*\Big\|_F^2\leq \cO_p\Big(\frac{\sigma^2Ks\log p_0}{n}\Big).
		\end{eqnarray*}
\end{Theorem}
\begin{Theorem}[(Lower Bound)]
	\label{th:lower_bound_asymmetric}
	Consider the class of incoherent sparse and low-rank tensors $\mathcal{F} = \{\mT: \mT = \sum_{k=1}^K \bbeta_{1k}\circ \bbeta_{2k}\circ \bbeta_{3k}, \|\bbeta_{i,k}\|_0 \leq s \text{ for } i = 1, 2, 3, k=1,\ldots, K\}$. If $\{\mX_i\}_{i=1}^n$ are i.i.d standard normal cubic sketchings, $\epsilon \overset{iid}{\sim} N(0,\sigma^2)$, $\min\{p_1,p_2,p_3\}\geq 20s$, and $s\geq 4$, we have
	\begin{equation}
	\inf_{\hat{\mT}}\sup_{\mT \in \mathcal{F}} \mathbb{E} \left\|\hat{\mT} - \mT\right\|_F^2\geq \frac{C\sigma^2sK\log(e\cdot p_0/s)}{n}.
	\end{equation}
\end{Theorem}
Theorems \ref{thm:main_nym} and \ref{th:lower_bound_asymmetric} imply that the proposed algorithm achieves a minimax-optimal rate of estimation error in the class of $\mathcal{F}$ as long as $\log(p_0) \asymp \log(p_0/s)$.

\section{Numerical Results}\label{sec:simu}

In this section, we investigate the effect of noise level, CP-rank, sample size, dimension, and sparsity on the estimation performance by simulation studies. We also investigate the numerical performance of the proposed algorithm when the incoherence assumption required in the theoretical analysis fails to hold.

In each setting, we generate $\mT^* = \sum_{k=1}^K \bbeta_k^*\circ \bbeta_k^* \circ \bbeta_k^*$, where $|\supp(\bbeta_k^*)|=s$, the support of $\bbeta_k^*$ is uniformly selected from $\{1,\ldots, p\}$, and the nonzero entries of $\bbeta_k^*$ are drawn randomly from standard normal distribution. Then, we calculate $\eta_k^* \leftarrow \|\bbeta_k^*\|_2^3$ and normalize $\bbeta_k^* \leftarrow \bbeta_k^*/ \|\bbeta_k^*\|_2$. The cubic sketchings $\{\mX_i\}_{i=1}^n$ are generated as $\mX_i = \bx_i\circ \bx_i\circ \bx_i$ and $\bx_i\overset{iid}{\sim} N(0, 1)$. The noise satisfies $\{\epsilon_i\}_{i=1}^n \overset{iid}{\sim} N(0, \sigma^2)$ or $\text{Laplace}(0, \sigma/\sqrt{2})$. Additionally, we adopt the following stopping rules in iterations: (1) the initialization iteration (Step 2 in Algorithm \ref{alg:overall_symmetric}) is stopped if $\|\bb_m^{(l+1)} - \bb_m^{(l)}\|_2\leq 10^{-6}$; (2) the gradient update iteration (Step 3 in Algorithm \ref{alg:overall_symmetric-gradient}) is stopped if $\|\bB^{(T+1)} - \bB^{(T)}\|_F\leq 10^{-6}$. The numerical results are based on 200 repetitions unless otherwise specified. The code was written in R and implemented on an Intel Xeon-E5 processor with
64 GB of RAM. 

\begin{figure}\label{fig_exact}
		\centering
			\includegraphics[width =0.52\linewidth,height=2.2in]{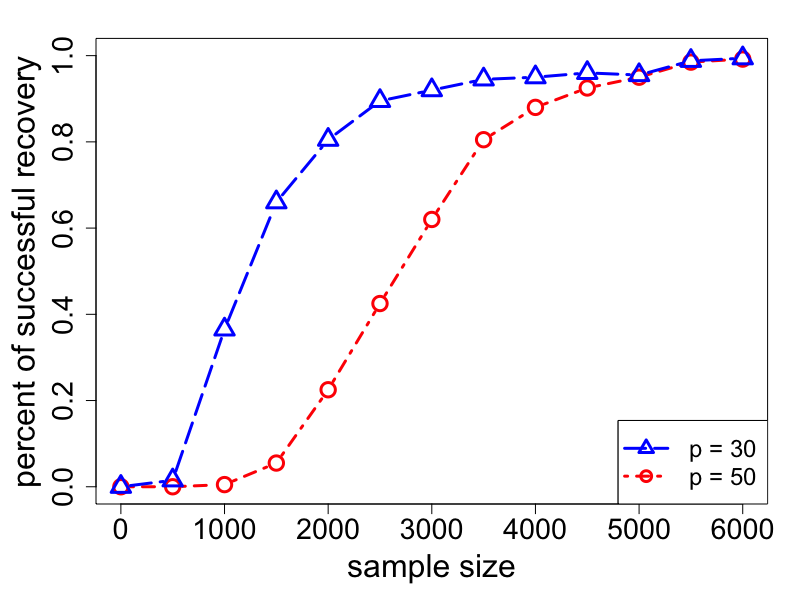}
			\caption{Successful rate of recovery with varying sample size}
	\end{figure}

First, we consider the percentage of successful recovery in the noiseless case. Let $K=3$, $s/p=0.3$, $p = 30$ or $50$, so that the total number of unknown parameters in $\mT^\ast$ is $2.7\times 10^4$ or $1.25\times 10^5$. The sample size $n$ ranges from 500 to 6000. Each recovery is called ``successful" if the relative error $\|\hat{\mT}-\mT^*\|_F/\|\mT^*\|_F < 10^{-4}$. We report the average successful recovery rate in Figure \ref{fig_exact}.
We can see from Figure \ref{fig_exact} that the empirical relation among successful recovery, dimension, and sample size is consistent with the theoretical results in Section \ref{sec:theory-symmetric}.

\begin{figure}
	\centering
	\includegraphics[width =0.48\linewidth,height=2.0in]{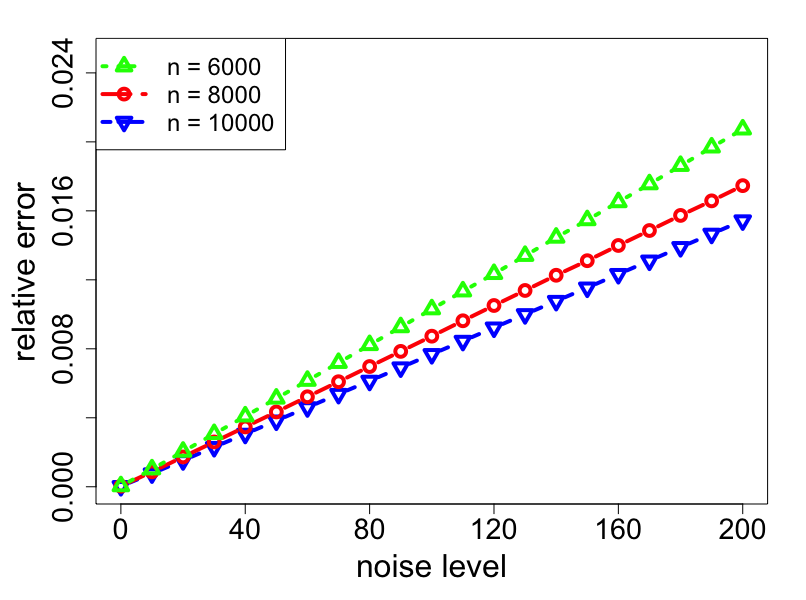}
	\includegraphics[width =0.48\linewidth,height=2.0in]{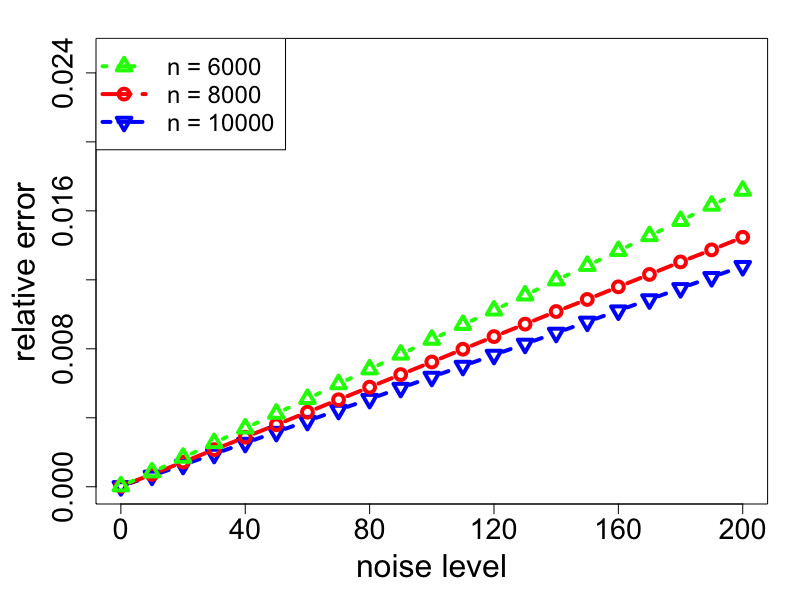}
	\caption{Estimation error under different noise levels. Left panel: $p=30$, right panel: $p=50$}
	\label{fig_noise_level}
	\includegraphics[width =0.48\linewidth,height=2.0in]{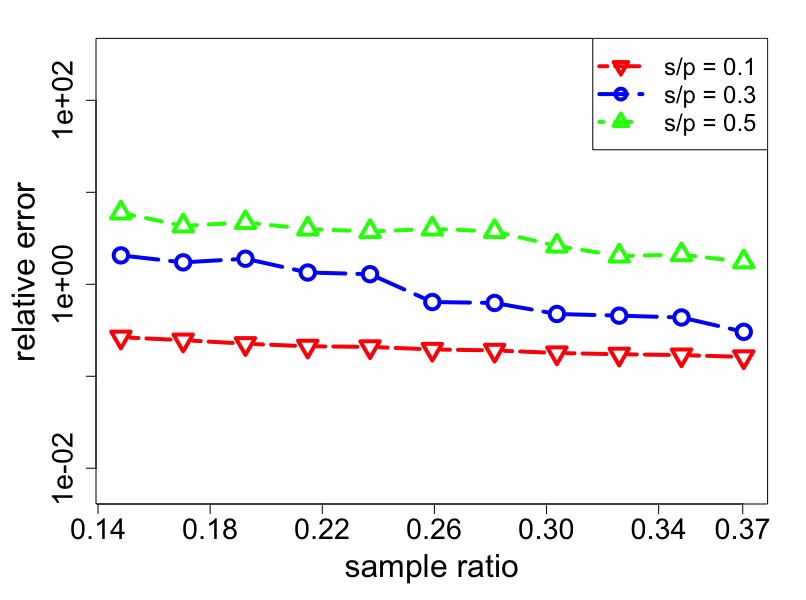}
	\includegraphics[width =0.48\linewidth,height=2.0in]{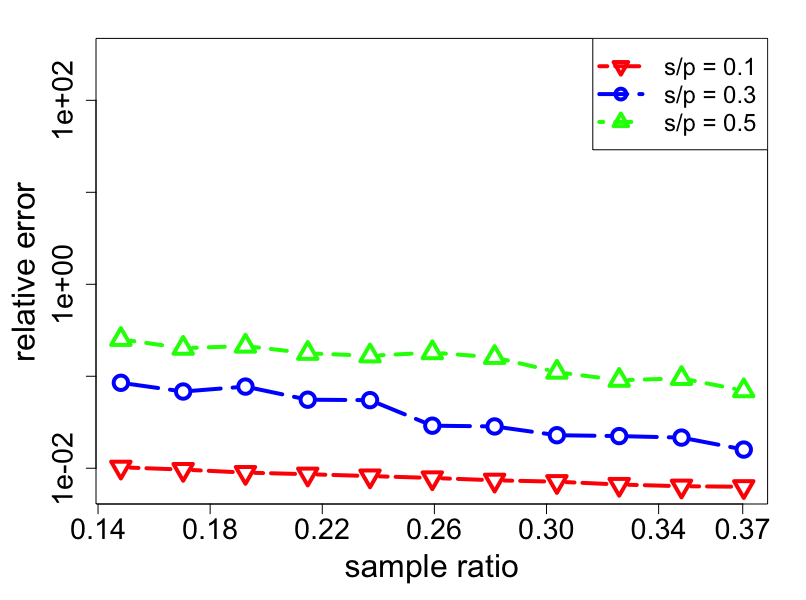}
	\caption{Estimation error under different dimension/sample ratio ($n/p^3$). Left panel: initial estimation error, right panel: final estimation error}
	\label{fig_sample_size}
\end{figure}

We then move to the noisy case. Select $K=3$, $s/p=0.3$, $p\in \{30, 50\}$, $\{\epsilon_i\}_{i=1}^n \overset{iid}{\sim} N(0, \sigma^2)$. We consider two scenarios: (1) sample size $n$ = 6000, 8000, or 10000, $s/p=0.3$, the noise level $\sigma$ varies from 0 to 200; (2) noise level $\sigma = 200$, sample size $n$ varies from 4000 to 10000, $p=30$, $s/p=0.1, 0.3, 0.5$. The estimation errors in terms of $\|\hat{\mT} - \mT^*\|_F / \|\mT^*\|_F$ in these two scenarios are plotted in Figures \ref{fig_noise_level} and \ref{fig_sample_size}, respectively. 
These results show that the proposed procedure achieves a good performance -- Algorithms \ref{alg:overall_symmetric} and \ref{alg:overall_symmetric-gradient} yield more accurate estimation with smaller variance $\sigma^2$ and/or large value of sample size $n$. 

\begin{figure}
\centering
    \includegraphics[width =0.48\linewidth,height=2.0in]{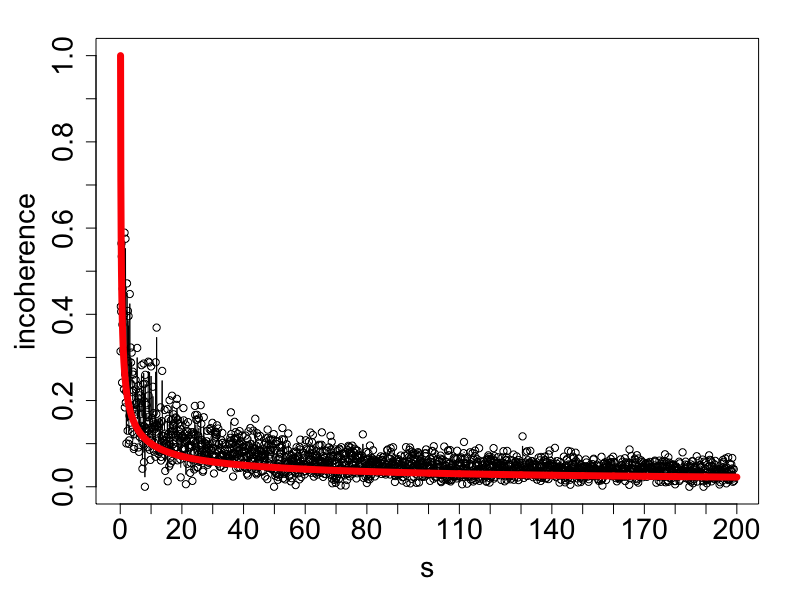}			\includegraphics[width =0.48\linewidth,height=2.0in]{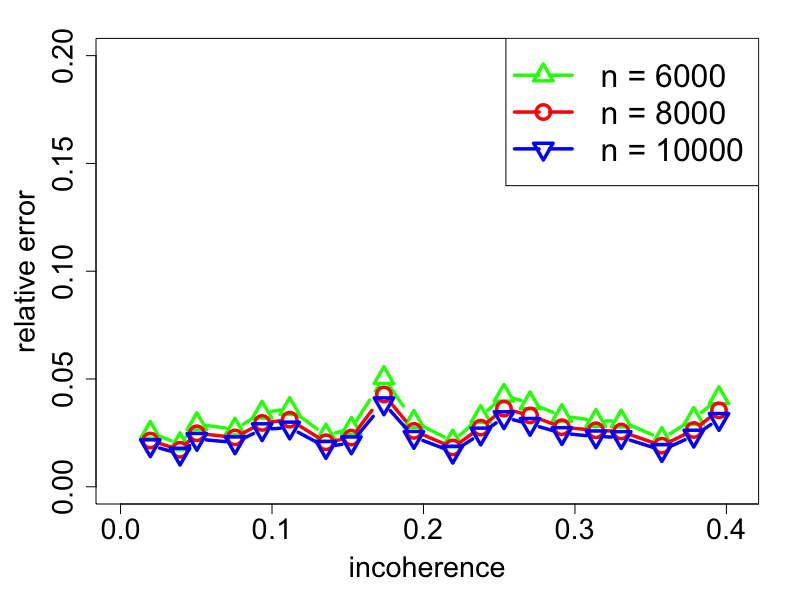}
			\caption{Left panel: incoherence parameter $\Gamma$ with varying sparsity. Here, the red line corresponds to the rate $\sqrt{s}$ required in the theoretical analysis. Right panel: average relative estimation error for tensors with varying incoherence.
			}
	\label{fig_incoherence}
	\end{figure}

Next, we demonstrate that the low-rank tensor parameter $\mT^*$ with randomly generated factors $\bbeta_k^*$ satisfies the incoherence condition \ref{con:incoherence_sym} with high probability. Set the CP-rank $K=3$ and the sparsity level $s/p = 0.3$ with the dimension $p$ ranging from 10 to 2000. We compute the incoherence parameter $\Gamma$ defined in Condition \ref{con:incoherence_sym}. The left panel of Figure \ref{fig_incoherence} shows that the incoherence parameter $\Gamma$ decays in a polynomial rate as $s$ grows, which matches the bound in Condition \ref{con:incoherence_sym}. Recall a theoretical justification on this point is also provided in Lemma \ref{lemma:proof_incoherence}. 

We further examine the performance of the proposed algorithm when the incoherence condition required in the theoretical analysis fails to hold. Specifically, we set the CP-rank $K=3$, $p=30$, and the sparsity level $s/p = 0.3$. We construct enormous copies of tensor parameter $\mT^*_j$ with i.i.d. standard normal factor vectors $\bbeta_k^*$. For each $\mT^*_j$, we calculate the incoherence $\Gamma_j$ defined in Condition \ref{con:incoherence_sym}, then manually pick 40 $\mT^*_{j'}$ such that 
$$0.01\cdot (j'-1)\leq \Gamma_{j'}\leq 0.01 \cdot j' \quad \text{for}\quad j' = \{1, 2, \ldots, 40\}.$$
In this way, we obtain a set of tensor parameters $\{\mT^*_{j'}\}$ with incoherence uniformly varying from 0 to 0.4. The right panel of Figure \ref{fig_incoherence} plots the relative error for estimating $\mT^*$ based on observations from cubic sketchings of $\mT^*_{j'}$ based on 1000 repetitions. We can see that the proposed algorithm  achieves small relative errors even when the true factors are highly coherent.

\begin{figure}
		\centering
		\includegraphics[width =0.48\linewidth,height=2.0in]{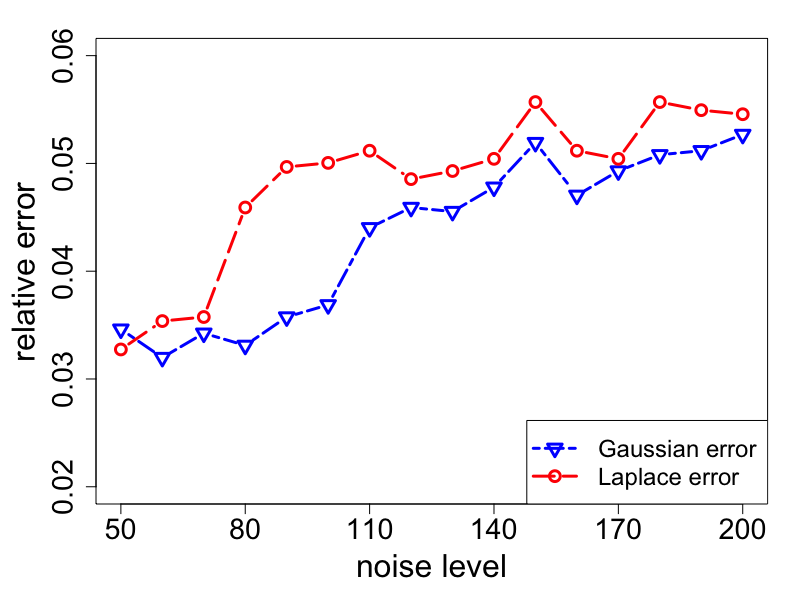}
		\caption{Comparison of estimation errors between Laplace error and Gaussian error}
		\label{fig_laplace_error}
	\end{figure}

Moreover, we consider a setting with Laplacian noise. Suppose $\{\epsilon_i\}_{i=1}^n \overset{iid}{\sim}Lap(\sigma)$ with density $f(x) = \frac{1}{\sigma}\exp(-2|x|/\sigma)$. With $n=3000$, $p=30$, and varying values of $\sigma$, the average estimation error and its comparison with Gaussian noise setting are provided in Figure \ref{fig_laplace_error}. We note that the estimation errors under Laplace noise are slightly higher than those under Gaussian noise.

\begin{figure}
		\centering
		\includegraphics[width =0.48\linewidth,height=2.0in]{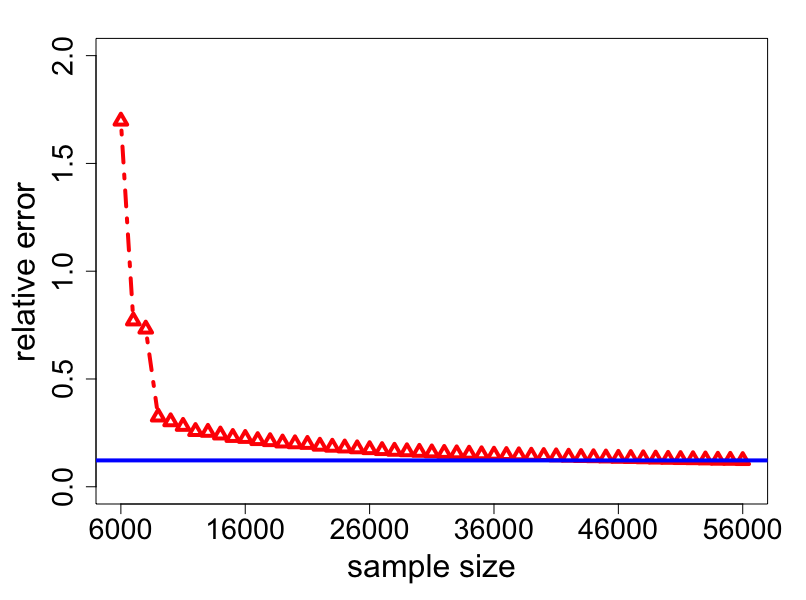}
		\includegraphics[width =0.48\linewidth,height=2.0in]{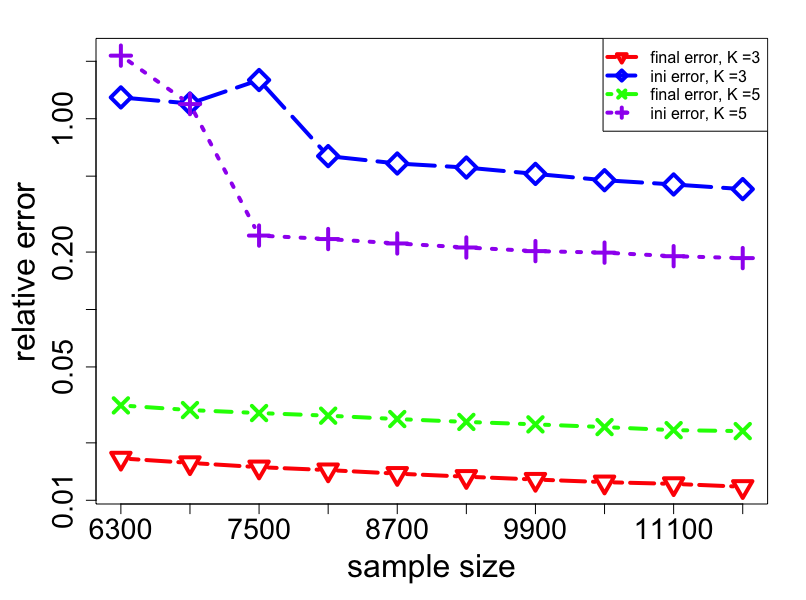}
		\caption{Log relative estimation error of initial estimation error (left panel) and  initialization/final estimation error (right panel)}
		\label{fig_initial_effect}
	\end{figure}

	We also compare the estimation errors of initial and final estimators for different ranks and sample sizes. Set $K=3, p =30, s/p = 0.3$ and consider the noiseless setting. It is clear from Figure \ref{fig_initial_effect} that the initialization error decays sufficiently, but does not converge to zero as sample size $n$ grows. This result matches our theoretical findings in Theorem \ref{thm:initial_symmetric}: as discussed in Remark \ref{rm:initialization error}, the initial stage may yield an inconsistent estimator due to the incoherence among $\bbeta_k$'s. 
	We also evaluate and compare the estimation errors for both initial and final estimators. From the right panel of Figure \ref{fig_initial_effect}, we can see that the final estimator is more stable and accurate compared to the initial one, which illustrates the merit of thresholded gradient descent step of the proposed procedure.

Finally, we compare the performance of the proposed method with the alternating least square (ALS)-based tensor regression method \cite{zhou2013}. We specifically consider two schemes for the initialization of ALS: (a) $\{\bbeta_k^{(0)}\}$ are i.i.d. standard Gaussian (cold start), and (b) $\{\bbeta_k^{(0)}\}$ are generated from the proposed Algorithm \ref{alg:overall_symmetric} (warm start). Setting $K=2$, $s/p=0.2$, $p=30$, $\{\epsilon_i\}_{i=1}^n \overset{iid}{\sim} N(0, 200^2)$, we apply both the proposed procedure and the ALS-based algorithm and record the average estimation errors with standard deviations for both initial and final estimators. From the result in Table \ref{tab:sim_result_tensorreg}, one can see the proposed algorithm significantly outperforms the ALS under both cold and warm start schemes. The main reason is pointed out in Remark \ref{rm:comparison-matrix-setting}: the cubic sketching setting possesses distinct aspects compared with the i.i.d. random Gaussian sketching setting, so that the method proposed by \cite{zhou2013} does not exactly fit here.
\begin{table}[htb]
	\centering
	\label{tab:sim_result_tensorreg}
	\vskip 1 em
	\begin{tabular}{c|c|c|c|c}
		\hline
		Sample size &ours & warm start & cold start & initial \\
		\hline
	   $n=4000$ &$4.023_{0.135}$   & $32.828_{1.798}$ & $37.785_{1.233}$& $38.032_{1.748}$  \\
		\hline
		$n=5000$&$1.945_{0.097}$  &$32.346_{2.343}$ & $36.962_{2.106}$ &$33.716_{1.786}$  \\
		\hline   
		$n=6000$ &$1.773_{0.092}$ &$22.220_{1.215}$  &$59.972_{3.407}$ & $25.579_{1.483}$ \\
		\hline
	\end{tabular}
	\caption{Estimation error and standard deviation (in subscript) of the proposed method and ALS-based method}
\end{table}

\section{Discussions}\label{sec:discussion}

This paper focuses on the third order tensor estimation via cubic sketchings. Moreover, all results can be extended to the higher-order case via high-order sketchings. To be specific, suppose
$$y_i = \langle \mT^\ast, \bx_i^{\otimes m} \rangle + \epsilon_i, \quad i=1,\ldots, n,$$
where $\mT^\ast \in (\mathbb{R}^{p})^{\otimes m}$ is an order-$m$, sparse, and low-rank tensor. In order to estimate $\mT^\ast$ based on $\{y_i, \bx_i\}_{i=1}^n$, one can first construct the order-$m$ moment-based estimator using a generalized version of Theorem \ref{thm:steins_higher} and the fact that the score functions $\cS_m(\bx) = (-1)^m \nabla^{m} p(\bx)/p(\bx)$ for the density function $p(\bx)$ satisfy a nice recursive equation: 
$$\cS_m(\bx):=-\cS_{m-1}(\bx)\circ \nabla\log p(\bx)-\nabla \cS_{m-1}(\bx).$$
Then, one can similarly perform high-order sparse tensor decomposition and thresholded gradient descent to estimate $\mT^*$. On the theoretical side,  
we can show if mild conditions hold and $n\geq C(\log n)^m(s\log p)^{m/2}$, the proposed procedure achieves
\begin{equation*}
\left\|\tilde{\mT} - \mT^*\right\|_F^2\lesssim \sigma^2\frac{Kms\log (p/s)}{n}
\end{equation*}
with high probability. 
The minimax optimality can be shown similarly.

\bibliographystyle{ieeetr}
\bibliography{cubic_sketching}

\newpage
\clearpage


\begin{center}
	\title{\Large Supplement to ``Sparse and Low-Rank Tensor Estimation via Cubic Sketchings"}
\end{center}

\medskip
\begin{center}
	Botao Hao\footnote{Ph.D student, Department of Statistics, Purdue University, West Lafayette, IN 47906. E-mail: hao22@purdue.edu.},~
	Anru Zhang\footnote{Assistant Professor, Department of Statistics, University of Wisconsin-Madison, Madison, WI 53706. E-mail: anruzhang@stat.wisc.edu.},~
	Guang Cheng\footnote{Professor, Department of Statistics, Purdue University, West Lafayette, IN 47906. E-mail: chengg@purdue.edu. Research Sponsored by NSF CAREER Award DMS-1151692, DMS-1418042, DMS-1712907 and Office of Naval Research ONR N00014-15-1-2331.}
\end{center}

\appendix
\bigskip

This supplementary contains five parts: (1) Section \ref{sec:proof} provides detailed proofs for empirical moment estimator and concentration results; (2) Section \ref{supp_sec:additional-proofs-main} provides additional proofs for the main theoretical results of this paper; (3) Section \ref{supp_sec:proof} contains detailed proofs for the theoretical developments in the main theorems; (4) Section \ref{sec:non_symmetric} covers the pseudo-code, conditions and main proofs of non-symmetric tensor estimation; (5) Section \ref{supp_sec:sgd} discusses the matrix form of gradient function and stochastic gradient descent; (6) Section \ref{sec:tech_lemma} provides several technical lemmas and their proofs.

\section{Proofs}\label{sec:proof}

\subsection{Moment Calculation} \label{sec:moment_calculation}

We first introduce three lemmas to show that the empirical moment based tensors \eqref{eqn:moment_tensor_sym}, \eqref{eqn:moment_inter}, and \eqref{eq:moment_asy} are all unbiased estimators for the target low-rank tensor in the corresponding scenarios. Detail proofs of three lemmas are postponed to Sections \ref{proof_moment_asy}, \ref{proof_moment_sym} and \ref{proof_moment_inter} in the supplementary materials.
	\begin{Lemma}[(Unbiasedness of moment estimator under non-symmetric sketchings)]
		\label{lemma:moment_asy}
		For non-symmetric tensor estimation model \eqref{eq:model_asymmetric} \& \eqref{eq:low_rank_nonsym}, define the empirical moment-based tensor $\cT$ by
		$$\cT:=\frac{1}{n}\sum_{i=1}^n y_i \bu_i\circ \bv_i\circ \bw_i.$$ Then $\cT$ is an unbiased estimator for $\mT^*$, i.e.,
		\begin{eqnarray*}
			\mathbb E(\cT) = \sum_{k=1}^K\eta_k^*\bbeta_{1k}^*\circ\bbeta_{2k}^*\circ\bbeta_{3k}^*.
		\end{eqnarray*}
	\end{Lemma}
	
	The extension to the symmetric case is non-trivial due to the dependency among three identical sketching vectors. We borrow the idea of high-order Stein's identity, which was originally proposed in \cite{janzamin2014score}. To fix the idea, we present only third order result for simplicity. The extension to higher-order is straightforward.
	\begin{Theorem}[(Third-order Stein's Identity, \cite{janzamin2014score})] \label{thm:steins_higher}
		Let $\bx \in \mathbb{R}^{p}$ be a random vector with joint density function $p(\bx)$. Define the third order score function $\cS_3(\bx):\mathbb R^p\rightarrow \mathbb R^{p\times p \times p}$ as $\cS_3(\bx) = -\nabla^{3}p(\bx)/p(\bx)$. Then for continuously differentiable function $G(\bx):\mathbb{R}^{p} \rightarrow \mathbb{R}$, 
		we have
		\begin{equation}\label{eqn:stein_identity}
		\mathbb E \left[ G(\bx) \cdot \cS_3(\bx) \right] = \mathbb E \left[ \nabla^{3} G(\bx) \right].
		\end{equation}
	\end{Theorem}
	In general, the order-$m$ high-order score function is defined as
	\begin{equation*} 
	\cS_m(\bx)=(-1)^m \frac{\nabla^{m} p(\bx)}{p(\bx)}.
	\end{equation*}
	Interestingly, the high-order score function has a recursive differential representation   
	\begin{equation}\label{eqn:recursive_form}
	\cS_m(\bx):=-\cS_{m-1}(\bx)\circ \nabla\log p(\bx)-\nabla \cS_{m-1}(\bx),
	\end{equation}
	with $\cS_0(\bx) = 1$. This recursive form is helpful for constructing unbiased tensor estimator under symmetric cubic sketchings. Note that the first order score function $\cS_1(\bx)=-\nabla \log p(\bx)$ is the same as score function in Lemma \ref{lemma:stein_lemma} (Stein's lemma \cite{stein2004use}).  The proof of Theorem \ref{thm:steins_higher} relies on iteratively applying the recursion representation of score function \eqref{eqn:recursive_form} and the first-order Stein's lemma (Lemma  \ref{lemma:stein_lemma}). We provide the detailed proof in Section \ref{proof:stein_higher} for the sake of completeness.
	
	In particular, if $\bx$ follows a standard Gaussian vector, each order score function can be calculated based on \eqref{eqn:recursive_form} as follows,
	\begin{equation}\label{eqn:S_x}
	\begin{split}
	&\cS_1(\bx) = \bx, \cS_2(\bx) = \bx\circ \bx -I_{d\times d},\\
	&\cS_3(\bx) = \bx\circ \bx\circ \bx - \sum_{j=1}^p\Big(\bx\circ \be_j\circ \be_j+ \be_j\circ \bx \circ \be_j+\be_j\circ \be_j\circ \bx\Big).
	\end{split}
	\end{equation}
	Interestingly, if we let $G(\bx) = \sum_{k=1}^K\eta_k^*(\bx^{\top}\bbeta_k^*)^3$, then 
	\begin{equation}\label{eqn:G_x}
	\frac{1}{6}\nabla^{3} G(\bx) = \sum_{k=1}^K\eta_k^*\bbeta_k^*\circ \bbeta_k^* \circ \bbeta_k^*,
	\end{equation}
	 which is exactly $\mT^*$. Connecting this fact with \eqref{eqn:stein_identity}, we are able to construct the unbiased estimator in the following lemma through high-order Stein's identity.

	\begin{Lemma}[(Unbiasedness of moment estimator under symmetric sketchings)]
		\label{lemma:moment_symmetric}
		Consider the symmetric tensor estimation model \eqref{eq:model} \& \eqref{eq:mathcal_F_symmetric}. Define the empirical first-order moment $\boldsymbol{m}_1 :=\frac{1}{n}\sum_{i=1}^n y_i\bx_i$.
		If we further define an empirical third-order-moment-based tensor $\cT_s$ by
		$$\cT_s:= \frac{1}{6}\Big[\frac{1}{n}\sum_{i=1}^n y_i \bx_i\circ \bx_i\circ \bx_i - \sum_{j=1}^p\Big(\boldsymbol{m}_1\circ \be_j\circ \be_j+ \be_j\circ \boldsymbol{m}_1 \circ \be_j+\be_j\circ \be_j\circ \boldsymbol{m}_1\Big)\Big],
		$$
		then
		\begin{eqnarray*}
			\mathbb E(\cT_s) = \sum_{k=1}^K\eta_k^*\bbeta_{k}^*\circ\bbeta_{k}^*\circ\bbeta_{k}^*.
		\end{eqnarray*}
	\end{Lemma}
	\emph{Proof.} Note that $y_i=G(\bx_i) +\epsilon_i$. Then we have
	\begin{equation*}
	\mathbb E\Big(\frac{1}{n}\sum_{i=1}^{n}y_i\cS_3(\bx)\Big) = \mathbb E\Big(\frac{1}{n}\sum_{i=1}^{n}(G(\bx_i)+\epsilon_i)\cS_3(\bx_i)\Big),
	\end{equation*}
	where $\cS_3(\bx)$ is defined in \eqref{eqn:S_x}. By using the conclusion in Theorem \ref{thm:steins_higher} and the fact \eqref{eqn:G_x}, we obtain
	\begin{equation*}
	\mathbb E(\cT_s) = \mathbb E\Big(\frac{1}{6n}\sum_{i=1}^{n}y_i\cS_3(\bx)\Big) = \sum_{k=1}^K\eta_k^*\bbeta_{k}^*\circ\bbeta_{k}^*\circ\bbeta_{k}^*,
	\end{equation*}
	since $\epsilon_i$ is independent of $\bx_i$. This ends the proof. \hfill $\blacksquare$\\

	Although the interaction effect model \eqref{eq:interaction_model} is still based on symmetric sketchings, we need much more careful construction for the moment-based estimator, since the first coordinate of the sketching vector is always constant 1. We give such an estimator in the following lemma.%
	\begin{Lemma}[(Unbiasedness of moment estimator in interaction model)]
		\label{lemma:moment_interaction}
		For interaction effect model \eqref{eq:interaction_model}, construct the empirical moment based tensor $\cT_{s'}$ as following 
	\begin{itemize}
		\item For $i,j,k\neq 0$, $\cT_{s'[i,j,k]} = \cT_{s[i,j,k]}$. And $\cT_{s'[i,j,0]} = \cT_{s[i,j,0]}, \cT_{s'[0,j,k]} = \cT_{s[0,j,k]}, \cT_{s'[i,0,k]} = \cT_{s[i,0,k]}$.
		\item For $i\neq 0$, $\cT_{s'[0,0,i]} = \cT_{s'[ 0, i, 0]} = \cT_{s'[i,0,0]} = \frac{1}{3}\cT_{s[0,0,i]} - \frac{1}{6}(\sum_{k=1}^p\cT_{s[k,k,i]}-(p+2)a_i)$.
		\item $\cT_{s'[0,0,0]} = \frac{1}{2p-2}(\sum_{k=1}^p\cT_{s[0,k,k]} - (p+2)\cT_{s[0,0,0]})$. 
	\end{itemize}
		The $\cT_{s'}$ is an unbiased estimator for $\mathcal{B}$, i.e.,
		\begin{eqnarray*}
			\mathbb E(\cT_{s'}) = \sum_{k=1}^K\eta_k\bbeta_{k}\circ\bbeta_{k}\circ\bbeta_{k}.
		\end{eqnarray*}
	\end{Lemma}

\subsection{Proofs of Lemmas \ref{cor:concentr_high_order} and \ref{lemma:tensor_spectral}: Concentration Inequalities}\label{subsec:concentration}

We aim to prove Lemmas \ref{cor:concentr_high_order} and \ref{lemma:tensor_spectral} in this subsection. These two lemmas provide key concentration inequalities of the theoretical analysis for the main result. Before going into technical details, we introduce a quasi-norm called $\psi_{\alpha}$-norm.

\begin{Definition}[($\psi_{\alpha}$-norm \cite{ALPT2011})]\label{def:phi_norm}
	The $\psi_{\alpha}$-norm of any random variable $X$ and $\alpha>0$ is defined as
	\begin{eqnarray*}
		\|X\|_{\psi_{\alpha}}:=\inf\Big\{C\in(0, \infty): ~ \mathbb E[\exp(|X|/C)^{\alpha}]\leq 2\Big\}.
	\end{eqnarray*}
\end{Definition}
Particularly, a random variable who has a bounded $\psi_2$-norm or bounded $\psi_1$-norm is called sub-Gaussian or sub-exponential random variable, respectively. 
Next lemma provides an upper bound for the $p$-th moment of sum of random variables with bounded $\psi_{\alpha}$-norm.
\begin{Lemma}\label{lemma:psi_moment_bound}
	Suppose $X_1,\ldots, X_n$ are $n$ independent random variables satisfying $\|X_i\|_{\psi_{\alpha}}\leq b$ with $\alpha> 0$, then for all $\ba=(a_1,\ldots, a_n)\in\mathbb R^n$ and $p\geq 2$,
	\begin{eqnarray}\label{eqn:psi_moment_bound}
	&&\Big(\mathbb E\Big|\sum_{i=1}^na_i X_i-\mathbb E(\sum_{i=1}^na_iX_i)\Big|^p\Big)^{\tfrac{1}{p}}\nonumber\\
	&\leq&
	\left\{\begin{array}{ll}
	C_1(\alpha)b\big(\sqrt{p}\|\ba\|_2+p^{1/\alpha}\|\ba\|_{\infty}\big), & \text{ if } 0<\alpha<1;\\
	C_2(\alpha)b\big(\sqrt{p}\|\ba\|_2+p^{1/\alpha}\|\ba\|_{\alpha^*}\big), & \text{ if } \alpha\geq 1.\\	
	\end{array}\right.
	\end{eqnarray}
	where $1/\alpha^*+1/\alpha = 1$, $C_1(\alpha),C_2(\alpha)$ are some absolute constants only depending on $\alpha$.
\end{Lemma}

If $0<\alpha<1$, \eqref{eqn:psi_moment_bound} is a combination of Theorem 6.2 in \cite{HMO1997} and the fact that the $p$-th moment of a Weibull variable with parameter $\alpha$ is of order $p^{1/\alpha}$. If $\alpha\geq 1$, \eqref{eqn:psi_moment_bound} follows from a combination of Corollaries 2.9 and 2.10 in \cite{talagrand1994supremum}. Continuing with standard symmetrization arguments, we reach the conclusion for general random variables. When $\alpha=1$ or 2, \eqref{eqn:psi_moment_bound} coincides with standard moment bounds for a sum of sub-Gaussian and sub-exponential random variables in \cite{RV12}. The detailed proof of Lemma \ref{lemma:psi_moment_bound} is postponed to Section \ref{proof_psi_moment_bound}.

When $0<\alpha<1$, by Chebyshev's inequality, one can obtain the following exponential tail bound for the sum of random variables with bounded $\psi_{\alpha}$-norm. This lemma generalizes the Hoeffding-type concentration inequality for sub-Gaussian random variables (see, e.g. Proposition 5.10 in \cite{RV12}), and Bernstein-type concentration inequality for sub-exponential random variables (see, e.g. Proposition 5.16 in \cite{RV12}).
\begin{Lemma}\label{lemma:orlicz_concentra}
	Suppose $0<\alpha<1$, $X_1,\ldots, X_n$ are independent random variables satisfying $\|X_i\|_{\psi_{\alpha}}\leq b$. Then there exists absolute constant $C(\alpha)$ only depending on $\alpha$ such that for any $\ba=(a_1,\ldots,a_n)\in \mathbb R^n$ and $0<\delta<1/e^2$,
	\begin{eqnarray*}
		\Big|\sum_{i=1}^na_iX_i-\mathbb E(\sum_{i=1}^na_iX_i)\Big|\leq C(\alpha)b\|\ba\|_2(\log\delta^{-1})^{1/2}+C(\alpha)b\|\ba\|_{\infty}(\log \delta^{-1})^{1/\alpha}
	\end{eqnarray*}
	with probability at least $1-\delta$.
\end{Lemma}

\emph{Proof.} For any $t>0$, by Markov's inequality,
\begin{equation*}
\begin{split}
&\mathbb P\Big(\Big|\sum_{i=1}^n a_iX_i-\mathbb E\Big(\sum_{i=1}^na_iX_i\Big)\Big|\geq t \Big)=\mathbb P\Big(\Big|\sum_{i=1}^n a_iX_i-\mathbb E\Big(\sum_{i=1}^na_iX_i\Big)\Big|^p\geq t^p \Big)\\
\leq& \frac{\mathbb E\Big|\sum_{i=1}^n a_iX_i-\mathbb E\Big(\sum_{i=1}^na_iX_i\Big)\Big|^p}{t^p}\leq \frac{C(\alpha)^pb^p\Big(\sqrt{p}\|\ba\|_2+p^{1/\alpha}\|\ba\|_{\infty}\Big)^p}{t^p},
\end{split}
\end{equation*}
where the last inequality is from Lemma \ref{lemma:psi_moment_bound}. We set $t$ such that $\exp(-p) =C(\alpha)^pb^p(\sqrt{p}\|\ba\|_2+p^{1/\alpha}\|\ba\|_{\infty})^p/t^p$. Then for $p\geq 2$,
\begin{equation*}
    \Big|\sum_{i=1}^n a_iX_i-\mathbb E\Big(\sum_{i=1}^na_iX_i\Big)\Big|\leq eC(\alpha)b\Big(\sqrt{p}\|\ba\|_2+p^{1/\alpha}\|\ba\|_{\infty}\Big)
\end{equation*}
holds with probability at least $1-\exp(-p)$. Letting $\delta= \exp(-p)$, we have that for any $0<\delta<1/e^2$,
\begin{equation*}
    \Big|\sum_{i=1}^n a_iX_i-\mathbb E\Big(\sum_{i=1}^na_iX_i\Big)\Big|\leq C(\alpha)b\Big(\|\ba\|_2(\log\delta^{-1})^{1/2}+\|\ba\|_{\infty}(\log \delta^{-1})^{1/\alpha}\Big),
\end{equation*}
holds with probability at least $1-\delta$. This ends the proof. \hfill $\blacksquare$

The next lemma provides an upper bound for the product of random variables in $\psi_{\alpha}$-norm.
\begin{Lemma}[($\psi_{\alpha}$ for product of random variables)]\label{lemma:sum_psi_norm}
	Suppose $X_1,\ldots, X_m$ are $m$ random variables (not necessarily independent) with $\psi_{\alpha}$-norm bounded by $\|X_j\|_{\psi_{\alpha}}\leq K_j$. Then the $\psi_{\alpha/m}$-norm of $\prod_{j=1}^mX_j$ is bounded as $$\left\|\prod_{j=1}^m X_j\right\|_{\psi_{\alpha/m}}\leq\prod_{j=1}^m K_j.$$ 
\end{Lemma}
\emph{Proof.} For any $\{x_j\}_{j=1}^m$ and  $\alpha>0$, by using the inequality of arithmetic and geometric means we have
\begin{equation*}
\Big(|\prod_{j=1}^m\frac{x_j}{K_j}|\Big)^{\alpha/m} = \Big(\prod_{j=1}^m|\frac{x_j}{K_j}|^{\alpha}\Big)^{1/m}\leq \frac{1}{m}\sum_{j=1}^m|\frac{x_j}{K_j}|^{\alpha}.
\end{equation*} 
Since exponential function is a monotone increasing function, it shows that
\begin{equation}\label{eqn:psi_1}
\begin{split}
&\exp\Big(|\prod_{j=1}^m\frac{x_j}{K_j}|\Big)^{\alpha/m}\leq  \exp\Big(\frac{1}{m}\sum_{j=1}^m|\frac{x_j}{K_j}|^{\alpha}\Big)\\
& = \Big(\prod_{j=1}^m \exp (|\frac{x_j}{K_j}|^{\alpha})\Big)^{1/m}\leq \frac{1}{m}\sum_{j=1}^m\exp\Big(|\frac{x_j}{K_j}|^{\alpha}\Big).
\end{split}
\end{equation}
From the definition of $\psi_{\alpha}$-norm, for $j=1,2,\ldots,m$, each individual $X_j$ has 
\begin{equation}\label{eqn:psi_2}
\mathbb E\Big(\exp(\frac{|X_j|}{K_j})^{\alpha}\Big)\leq 2.
\end{equation}
Putting \eqref{eqn:psi_1} and \eqref{eqn:psi_2} together, we obtain
\begin{equation*}
\begin{split}
	&\mathbb E\Big[\exp\Big(|\frac{\prod_{j=1}^mX_j}{\prod_{j=1}^mK_j}|\Big)^{\alpha/m}\Big] = \mathbb E\Big[\exp\Big(|\prod_{j=1}^m\frac{X_j}{K_j}|\Big)^{\alpha/m}\Big]\\
	&\leq \frac{1}{m}\sum_{j=1}^m\mathbb E\Big[\exp\Big(|\frac{X_j}{K_j}|\Big)^{\alpha}\Big]\leq 2.
	\end{split}
\end{equation*}
Therefore, we conclude that the $\psi_{\alpha/m}$-norm of $\prod_{j=1}^mX_j$ is bounded by $\prod_{j=1}^mK_j$.  \hfill $\blacksquare$\\

\emph{Proof of Lemma \ref{cor:concentr_high_order}.} 
Note that for any $j=1,2,\ldots, m$, the $\psi_2$-norm of $\bX_j^{\top}\bbeta_j$ is bounded by $\|\bbeta_j\|_2$ \cite{RV12}. According to Lemma \ref{lemma:sum_psi_norm}, the $\psi_{2/m}$-norm of $\prod_{j=1}^m(\bX_j^{\top}\bbeta_j)$ is bounded by $\prod_{j=1}^m\|\bbeta_j\|_2$. Directly applying Lemma \ref{lemma:orlicz_concentra}, we reach the conclusion.
\hfill $\blacksquare$\\

\emph{Proof of Lemma \ref{lemma:tensor_spectral}.}  We first focus on the non-symmetric version and the proof follows three steps:
\begin{enumerate}
	\item Truncate the first coordinate of $\bx_{1i}, \bx_{2i}, \bx_{3i}$ by a carefully chosen truncation level;
	\item Utilize the high-order concentration inequality in Lemma \ref{lemma:tensor_concentration} at order three;
	\item Show that the bias caused by truncation is negligible. 
\end{enumerate}

With slightly abuse of notations, we denote $a,x,y$ etc. as their \emph{first coordinate} of $\ba,\bx,\by$ etc.  Without loss of generality, we assume $p:=\max\{p_1,p_2,p_3\}$. By unitary invariance, we assume $\bbeta_1=\bbeta_2=\bbeta_3=\be_1$, where $\be_1 = (1,0,\ldots,0)^{\top}$. Then, it is equivalent to prove
\begin{eqnarray*}
	&&\Big\|M_{\text{nsy}} - \mathbb E(M_{\text{nsy}})\Big\|_s =\Big\|\frac{1}{n}\sum_{i=1}^nx_{1i}x_{2i}x_{3i}\bx_{1i}\circ\bx_{2i}\circ\bx_{3i}-\be_1\circ\be_1\circ\be_1\Big\|_{s}\\
	&\leq&  C(\log n)^3\Big(\sqrt{\frac{s^3\log^3 (p/s)}{n^2}} + \sqrt{\frac{s\log (p/s)}{n}}\Big) .
\end{eqnarray*}

Suppose $\bx_1\sim \cN(0, \bI_{p_1}), \bx_2\sim \cN(0, \bI_{p_2}), \bx_3\sim \cN(0, \bI_{p_3})$ and $\{\bx_{1i}, \bx_{2i}, \bx_{3i}\}_{i=1}^n$ are $n$ independent samples of $\{\bx_1, \bx_2, \bx_3\}$. And define a bounded event $\cG_n$ for the first coordinate and its corresponding population version,
\begin{eqnarray*}
	\cG_n=\{\max_{i}\{|x_{1i}|,|x_{2i}|, |x_{3i}|\}\leq M \},
	\cG=\{\max\{|x_1|,|x_2|, |x_3|\}\leq M\},
\end{eqnarray*}
where $M$ is a large constant to be specified later. Decomposing $\|M_{\text{nsy}} - \mathbb E(M_{\text{nsy}})\|_s$ as 
\begin{equation*}
\begin{split}
	&\Big\|M_{\text{nsy}} - \mathbb E(M_{\text{nsy}})\Big\|_s\\
	&\leq\underbrace{\Big\|\frac{1}{n}\sum_{i=1}^nx_{1i}x_{2i}x_{3i}\bx_{1i}\circ\bx_{2i}\circ\bx_{3i}-\mathbb E\Big(x_1x_2x_3\bx_1\circ\bx_2\circ\bx_3\big|\cG\Big)\Big\|_{s}}_{M_1:\text{main term}}\\
	&+\underbrace{\Big\|\mathbb E\Big(x_1x_2x_3\bx_1\circ\bx_2\circ\bx_3\big|\cG\Big)-\be_1\circ\be_1\circ\be_1\Big\|_{s}}_{M_2:\text{bias term}},\nonumber
\end{split}
\end{equation*}
 we will prove that $M_2$ is negligible in terms of convergence rate of $M_1$. 

\textbf{Bounding $M_1$.} For simplicity, we define
$\bx_{1}' = \bx_{1}| \cG, \ \bx_{2}' = \bx_{2}|\cG , \ \bx_{3}' = \bx_{3}| \cG
$, and $\{\bx_{1i}',\bx_{2i}', \bx_{3i}'\}_{i=1}^n$ are $n$ independent samples of $\{\bx_1', \bx_2', \bx_3'\}$. According to the law of total probability, we have
\begin{equation*}
\begin{split}
	&\mathbb P\Big(M_1\geq t\Big)\leq \mathbb P\Big(\cG_n^c\Big)\\
	&+ \mathbb P\Big(\underbrace{\Big\|\frac{1}{n}\sum_{i=1}^nx_{1i}'\bx_{1i}'\circ x_{3i}'\bx_{2i}'\circ x_{i1}'\bx_{3i}'-\mathbb E\Big(x_1'\bx_1'\circ x_2'\bx_2'\circ x_3'\bx_3'\Big)\Big\|_{s}}_{M_{11}}\geq t\Big).
	\end{split}
\end{equation*}
According to Lemma \ref{lemma:bounded_gaussian}, the entry of $x_{1i}'\bx_{1i}', x_{2i}'\bx_{2i}', x_{3i}'\bx_{3i}'$ are sub-Gaussian random variable with $\psi_2$-norm $M^2$. Applying Lemma \ref{lemma:tensor_concentration}, we obtain
\begin{eqnarray*}
	\mathbb P\Big(M_{11}\geq C_1M^{6}\delta_{n,s}\Big)\leq \frac{1}{p},
\end{eqnarray*}
where $\delta_{n,s} = ((s\log(p/s))^{3}/n^2)^{1/2}+(s\log(p/s)/n)^{1/2}$.

On the other hand, 
\begin{eqnarray*}
	\mathbb P(\cG_n^c)\leq 3\sum_{i=1}^n \mathbb P(|x_{1i}|\geq M)\leq 3ne^{1-C_2M^2}
\end{eqnarray*}
Putting the above bounds together, we obtain
\begin{eqnarray*}
	\mathbb P\Big(M_1\geq C_1 M^6 \delta_{n,s}\Big)\leq \frac{1}{p}+3n e^{1-C_2M^2}.
\end{eqnarray*}
By setting $M=2\sqrt{\log n/C_2}$, the bound of $M_1$ reduces to
\begin{eqnarray}
\label{bound_M1}
\mathbb P\Big(M_1\geq \frac{64C_1}{C_2^3}\delta_{n,s}(\log n)^3\Big)\leq \frac{1}{p}+\frac{3e}{n^3}.
\end{eqnarray}

\textbf{Bounding $M_2$.}  From the definitions of $M_2$ and sparse spectral norm, 
\begin{eqnarray*}
&&M_2 = \Big\|\mathbb E\Big(x_1x_2x_3\bx_1\circ\bx_2\circ\bx_3\big|\cG\Big)-\be_1\circ\be_1\circ\be_1\Big\|_{s}\\
&=& \sup_{\substack{\|\ba\|=\|\bb\|=\|\bc\|=1 \\ \max\{\|\ba\|_0,\|\bb\|_0, \|\bc\|_0\}\leq s}}\Big| \mathbb E\Big(x_1x_2x_3(\bx_1^{\top}\ba)(\bx_2^{\top}\bb)(\bx_3^{\top}\bc)\big|\cG\Big)-a_1b_1c_1\Big|.
\end{eqnarray*}
Since $x_{1j}$ is independent of $x_{1k}$ for any $j\neq k$, $\mathbb E(x_1(\bx_1^{\top}\varrho\ba)|\cG) = \mathbb E(x_1^2a_1|\cG)$. Similar results hold for $\bx_2, \bx_3$. Then we have
\begin{eqnarray*}
&&M_2 = \sup_{\substack{\|\ba\|=\|\bb\|=\|\bc\|=1 \\ \max\{\|\ba\|_0,\|\bb\|_0, \|\bc\|_0\}\leq s}} |a_1b_1c_1|\Big| \mathbb E\Big(x_1^2x_2^2x_3^2\big|\cG\Big)-1\Big|\leq \Big| \mathbb E\Big(x_1^2x_2^2x_3^2\big|\cG\Big)-1\Big|,\\
&=& \Big|\mathbb E\Big(x_1^2\Big||x_1|\leq M\Big) \mathbb E\Big(x_2^2\Big||x_2|\leq M\Big) \mathbb E\Big(x_3^2\Big||x_3|\leq M\Big) - 1\Big|.
\end{eqnarray*}

By the basic property of Gaussian random variable, we can show
\begin{eqnarray*}
	1\geq \mathbb E\Big(x_i^2\big||x_i|\leq M\Big)\geq 1-2Me^{-M^2/2} ,\quad  i=1,2,3.
\end{eqnarray*}
Plugging them into $M_2$, we have
\begin{eqnarray*}
	M_2 &\leq& \Big|\Big(1-2Me^{-M^2/2}\Big)^3-1\Big|\\
	&\leq& \Big|12M^2e^{-M^2}-6Me^{-M^2/2}-8M^3e^{-3M^2/2}\Big|\\
	&\leq& \Big|26M^3e^{-M^2/2}\Big|,
\end{eqnarray*}
where the last inequality holds for a large $M>0$.
By the choice of $M=2\sqrt{\log n/C_2}$, we have $M_2\leq 208/C_2^{3/2}(\log n)^{\frac{3}{2}}/n^2$ for some constant $C_2$. When $n$ is large, this rate is negligible comparing with \eqref{bound_M1}

\textbf{Bounding $M$:} We put the upper bounds of $M_1$ and $M_2$ together. After some adjustments for absolute constant, it suffices to obtain
\begin{eqnarray*}
	\mathbb P\Big(M_1+M_2\leq C(\log n)^3\Big(\sqrt{\frac{s^3\log^3 (p/s)}{n^2}} + \sqrt{\frac{s\log (p/s)}{n}}\Big)\Big)\geq 1- \frac{10}{n^3}-\frac{1}{p}.
\end{eqnarray*}
This concludes the proof of non-symmetric part. The proof of symmetric part remains similar and thus is omitted here.
\hfill $\blacksquare$\\

\section{Additional Proofs for main results}\label{supp_sec:additional-proofs-main}

\subsection{Proof of Theorem \ref{thm:initial_symmetric}: Initialization Effect}\label{subsec:proof_ini}

Theorem \ref{thm:initial_symmetric} gives an approximation error upper bound for the sparse-tensor-decomposition-based initial estimator. 
In Step I of Section \ref{sec:procedure-symmetric}, the original problem can be reformatted to a version of tensor denoising:
\begin{eqnarray}
\label{eqn:denoising_model}
\cT_s = \mT^*+\cE, \quad\text{where} \quad \cE = \cT_s - \mathbb E(\cT_s).
\end{eqnarray}
The key difference between our model \eqref{eqn:denoising_model} and recent works \cite{AGJ14, Sun2016} is that $\cE$ arises from empirical moment approximation, rather than the random observation noise considered in \cite{AGJ14} and \cite{Sun2016}.  
Next lemma gives an upper bound for the approximation error. The proof of Lemma \ref{lemma:perturbation_analysis} is deferred to Section \ref{proof_perturbation}.
\begin{Lemma}[(Approximation error of $\mathcal{T}_s$)]\label{lemma:perturbation_analysis}
	Recall that $\cE = \cT_s - \mathbb E(\cT_s)$, where $\cT_s$ is defined in (\ref{eqn:moment_tensor_sym}). Suppose Condition \ref{con:noise} is satisfied and $s\leq d \leq Cs$. Then
	\begin{eqnarray}\label{ineq:perturbation_error}
	\|\cE\|_{s+d} &\leq& 2C_1 \sum_{k=1}^K\eta_k^*\Big(\sqrt{\frac{s^3\log^3(p/s)}{n^2}}+\sqrt{\frac{s\log(p/s)}{n}}\Big)(\log n)^4
	\end{eqnarray}
	with probability at least $1-5/n$ for some uniform constant $C_1$. 
\end{Lemma}

Next we denote the following quantity for simplicity,
\begin{eqnarray}\label{def:gamma}
\gamma = C_2\min\Big\{\frac{R^{-1}}{6}-\frac{\sqrt{K}}{s}, \frac{R^{-1}}{4\sqrt{5}}-\frac{2}{\sqrt{s}}\Big(1+\sqrt{\frac{K}{s}}\Big)^2\}, 
\end{eqnarray}
where $R$ is the singular value ratio, $K$ is the CP-rank, $s$ is the sparsity parameter, $\Gamma$ is the incoherence parameter and $C_2$ is uniform constant. 

Next lemma provides theoretical guarantees for sparse tensor decomposition method.
\begin{Lemma}\label{lemma:ini_sym}
	Suppose that the symmetric tensor denoising model (\ref{eqn:denoising_model}) satisfies Conditions \ref{con:ident}, \ref{con:parameter} and \ref{con:incoherence_sym} (i.e., the identifiability, parameter space and incoherence). Assume the number of initializations $L \geq K^{C_3\gamma^{-4}}$ and the number of iterations $N\geq C_4\log \left(\gamma/\left(\frac{1}{\eta_{\min}^*}\|\cE\|_{s+d} + \sqrt{K}\Gamma^2\right)\right)$ for constants $C_3,C_4$, the truncation parameter $s\leq d \leq Cs$. Then the sparse-tensor-decomposition-based initialization satisfies 
	\begin{eqnarray}
	\label{ineq:initial_error}
	\max\left\{\|\bbeta_k^{(0)} - \bbeta_k^*\|_2, |\eta_k^{(0)} - \eta_k^*|\right\}\leq  \frac{C_4}{\eta_{\min}^*}\|\cE\|_{s+d} + \sqrt{K}\Gamma^2,
	\end{eqnarray}
	for any $k\in[K]$.
\end{Lemma}

The proof of Lemma \ref{lemma:ini_sym} essentially follows Theorem 3.9  in \cite{Sun2016}, we thus omit the detailed proof here. The upper bound in \eqref{ineq:initial_error} contains two terms: $\frac{C_4}{\eta^\ast_{\min}}\|\cE\|_{s+d}$ and $\sqrt{K}\Gamma^2$, which are due to the empirical moment approximation and the incoherence among different $\bbeta_k$, respectively. 

Although the sparse tensor decomposition is not optimal in statistical rate, it does offer a reasonable initial estimation provided enough samples. Equipped with \eqref{ineq:perturbation_error} and Condition \ref{con:parameter}, the right side of \eqref{ineq:initial_error} reduces to 
\begin{eqnarray*}
	&&\frac{C_4}{\eta_{\min}^*}\|\cE\|_{s+d} + \sqrt{K}\Gamma^2 \\
	&\leq& 2 C_1C_4KR\Big(\sqrt{\frac{s^3\log^3(p/s)}{n^2}}+\sqrt{\frac{s\log(p/s)}{n}}\Big)(\log n)^4+\sqrt{K}\Gamma^2,
\end{eqnarray*} 
with probability at least $1-5/n$. Denote $C_0=4\cdot 2160\cdot C_1C_4$. Using Conditions \ref{con:incoherence_sym} and \ref{con:sample}, we reach the conclusion that
$$\max\left\{\|\bbeta_k^{(0)} - \bbeta_k^*\|_2, |\eta_k^{(0)} - \eta_k^*|\right\}\leq  K^{-1}R^{-2}/2160,$$
with probability at least $1-5/n$. \hfill $\blacksquare$\\	

\subsection{Proof of Theorem \ref{thm:symmetric_main}: Gradient Update}
\label{subsec:gradient-update-proof}

We first introduce the following lemma to illustrate the improvement of one step thresholded gradient update under suitable conditions. 
The error bound includes two parts: the optimization error that describes one step effect for gradient update, and the statistical error that reflects the random noise effect. 
The proof of Lemma \ref{lemma:onestep_update_sym} is given in Section \ref{proof:onestep_update_sym}. For notation simplicity, we drop the superscript of $\eta_k^{(0)}$ in the following proof.

\begin{Lemma}
	\label{lemma:onestep_update_sym}
	Let $t\geq 0$ be an integer. Suppose Conditions \ref{con:ident}-\ref{con:sample} hold and $\{\bbeta_k^{(t)}, \eta_k\}$ satisfies the following upper bound
	\begin{equation}\label{con:initial_input}
	\sum_{k=1}^K\Big\|\sqrt[3]{\eta_k}\bbeta_k^{(t)}-\sqrt[3]{\eta_k^*}\bbeta_k^*\Big\|_2^2\leq 4K\eta_{\max}^{*\tfrac{2}{3}}\varepsilon_0^2, \ \max_{k\in[K]}\Big|\eta_k-\eta_k^*\Big|\leq \varepsilon_0,
	\end{equation}
	with probability at least $1-\cO(K/n)$, where $\varepsilon_0= K^{-1}R^{-\tfrac{4}{3}}/2160$. As long as the step size $\mu$ satisfies
	\begin{equation}\label{eqn:mu_0}
	0<\mu\leq \mu_0= \frac{32 R^{-20/3}}{3K[220+270K]^2},
	\end{equation}
	then $\{\bbeta_k^{(t+1)}\}$ can be upper bounded as
	\begin{equation*}
	\begin{split}
	&\sum_{k=1}^K\Big\|\sqrt[3]{\eta_k}\bbeta_k^{(t+1)}-\sqrt[3]{\eta_k^*}\bbeta_k^*\Big\|_2^2\\
	&\leq \underbrace{\Big(1-32\mu K^{-2}R^{-\tfrac{8}{3}}\Big)\sum_{k=1}^K\Big\|\sqrt[3]{\eta_k}\bbeta_k^{(t)}-\sqrt[3]{\eta_k^*}\bbeta_k^*\Big\|_2^2}_{\text{optimization error}}+\underbrace{2C_0\mu^2K^{-2}R^{-\tfrac{8}{3}}\eta_{\min}^{*-\tfrac{4}{3}}\frac{\sigma^2s\log p}{n}}_{\text{statistical error}},
	\end{split}
	\end{equation*}
	with probability at least $1-\cO(Ks/n)$. 
\end{Lemma}

In order to apply Lemma \ref{lemma:onestep_update_sym}, we prove that the required condition \eqref{con:initial_input} holds at every iteration step $t$ by induction. When $t=0$, by \eqref{ineq:initial} and Condition \ref{con:parameter},
\begin{eqnarray*}
	\Big\|\bbeta_k^{(0)}-\bbeta_k^*\Big\|_2\leq \varepsilon_0, \ \Big|\eta_k-\eta_k^*\Big|\leq \varepsilon_0, \ \text{for} \ k\in[K],
\end{eqnarray*}
holds with probability at least $1-\cO(1/n)$. Since the initial estimator output by first stage is normalized, i.e., $\|\bbeta_k^{(0)}\|_2 = \|\bbeta_k^\ast\|_2 = 1$, by triangle inequality we have 
\begin{eqnarray*}
	\Big\|\sqrt[3]{\eta_k}\bbeta_k^{(0)} -\sqrt[3]{\eta_k^*}\bbeta_k^{*} \Big\|_2&\leq& \Big\|\sqrt[3]{\eta_k}\bbeta_k^{(0)} - \sqrt[3]{\eta_k^*}\bbeta_k^{(0)} + \sqrt[3]{\eta_k^*}\bbeta_k^{(0)} -\sqrt[3]{\eta_k^*}\bbeta_k^{*} \Big\|_2\\
	&\leq& |\sqrt[3]{\eta_k}-\sqrt[3]{\eta_k^*}|+\sqrt[3]{\eta_k^*}\Big\|\bbeta_k^{(0)} -\bbeta_k^{*} \Big\|_2.
\end{eqnarray*}
Note that $$\Big|\sqrt[3]{\eta_k}-\sqrt[3]{\eta_k^*}\Big|\leq \frac{\varepsilon_0}{(\sqrt[3]{\eta_k})^2+\sqrt[3]{\eta_k\eta_k^*}+(\sqrt[3]{\eta_k^*})^2}\leq \varepsilon_0\sqrt[3]{\eta_k^*}.$$
This implies
\begin{eqnarray*}
	\Big\|\sqrt[3]{\eta_k}\bbeta_k^{(0)} - \sqrt[3]{\eta_k^*}\bbeta_k^*\Big\|_2\leq2\sqrt[3]{\eta_k^*}\varepsilon_0,
\end{eqnarray*}
with probability at least $1-\cO(1/n)$.
Taking the summation over $k\in[K]$, we have
\begin{equation*}
\sum_{k=1}^K\Big\|\sqrt[3]{\eta_k}\bbeta_k^{(0)} - \sqrt[3]{\eta_k^*}\bbeta_k^*\Big\|_2^2 \leq \sum_{k=1}^K4\eta_k^{*\tfrac{2}{3}}\varepsilon_0^2\leq 4 K\eta_{\max}^{*\tfrac{2}{3}}\varepsilon_0^2,
\end{equation*}
with probability at least $1-\cO(K/n)$, which means \eqref{con:initial_input} holds for $t= 0$.

Suppose \eqref{con:initial_input} holds at the iteration step $t-1$, which implies 
\begin{equation*}
\begin{split}
&\sum_{k=1}^{K}\Big\|\sqrt[3]{\eta_k}\bbeta_k^{(t)} -\sqrt[3]{\eta_k^*}\bbeta_k^{*} \Big\|_2^2\\
\leq& \Big(1-32\mu K^{-2}R^{-\tfrac{8}{3}}\Big)\sum_{k=1}^K\Big\|\sqrt[3]{\eta_k}\bbeta_k^{(t-1)} -\sqrt[3]{\eta_k^*}\bbeta_k^{*} \Big\|_2^2+\mu2C_0K^{-2}R^{-\tfrac{8}{3}}\eta_{\min}^{*\tfrac{4}{3}}\frac{\sigma^2s\log p}{n}\\
\leq& 4K\eta_{\max}^{*\tfrac{2}{3}}\varepsilon_0^2-\mu\Big(128KR^{-\tfrac{8}{3}}\eta_{\max}^{*\tfrac{2}{3}}\varepsilon_0^2-2C_0K^{-2}R^{-\tfrac{8}{3}}\eta_{\min}^{*\tfrac{4}{3}}\frac{\sigma^2s\log p}{n}\Big).
\end{split}
\end{equation*}
Since Condition \ref{con:sample} automatically implies
\begin{equation*}
\frac{n}{s\log p}\geq \frac{C_0\sigma^2R^{-\tfrac{2}{3}}\eta_{\min}^{*\tfrac{2}{3}}K}{64\varepsilon_0^2},
\end{equation*}
for a sufficiently large $C_0$, we can obtain
\begin{eqnarray*}
	\sum_{k=1}^{K}\Big\|\sqrt[3]{\eta_k}\bbeta_k^{(t)} -\sqrt[3]{\eta_k^*}\bbeta_k^{*} \Big\|_2^2\leq 4K\eta_{\max}^{*\tfrac{2}{3}}\varepsilon_0^2.
\end{eqnarray*}
By induction, \eqref{con:initial_input} holds at each iteration step.

Now we are able to use Lemma \ref{lemma:onestep_update_sym} recursively to complete the proof. Repeatedly using Lemma \ref{lemma:onestep_update_sym}, we have for $t=1,2,\ldots,$
\begin{equation*}
\begin{split}
&\sum_{k=1}^{K}\Big\|\sqrt[3]{\eta_k}\bbeta_k^{(t+1)} -\sqrt[3]{\eta_k^*}\bbeta_k^{*} \Big\|_2^2\\
&\leq \Big(1-32\mu K^{-2}R^{-\tfrac{8}{3}}\Big)^t\sum_{k=1}^K\Big\|\sqrt[3]{\eta_k}\bbeta_k^{(0)} -\sqrt[3]{\eta_k^*}\bbeta_k^{*} \Big\|_2^2+\frac{C_0\eta_{\min}^{*-\tfrac{4}{3}}}{16}\frac{\sigma^2s\log p}{n},
\end{split}
\end{equation*}
with probability at least $1-\cO(tKs/n)$. This concludes the first part of Theorem \ref{thm:symmetric_main}. 


When the total number of iterations is no smaller than 
\begin{eqnarray*}
	T^* = \frac{\log(C_3\eta_{\min}^{*-4/3}\sigma^2s\log p)-\log(64\eta_{\max}^{*2/3}K\varepsilon_0 n)}{\log(1-32\mu K^{-2}R^{-8/3})},
\end{eqnarray*}
the statistical error will dominate the whole error bound in the sense that
\begin{eqnarray}
\label{eqn:factor_error_sym}
\sum_{k=1}^{K}\Big\|\sqrt[3]{\eta_k}\bbeta_k^{(T^*)} -\sqrt[3]{\eta_k^*}\bbeta_k^{*} \Big\|_2^2\leq \frac{C_3\eta_{\min}^{*-\tfrac{4}{3}}}{8}\frac{\sigma^2s\log p}{n},
\end{eqnarray}
with probability at least $1-\cO(T^*Ks/n)$.

The next lemma shows that the Frobenius norm distance between two tensors can be bounded by the distances between each factors in their CP decomposition. The proof of this lemma is provided in Section \ref{proof:error_transfer}.
\begin{Lemma}\label{lemma:error_transfer}
	Suppose $\mT$ and $\mT^*$ have CP-decomposition $\mT = \sum_{k=1}^K\eta_k\bbeta_k\circ \bbeta_k\circ \bbeta_k$ and $\mT^* = \sum_{k=1}^K\eta_k^*\bbeta_k^*\circ \bbeta_k^*\circ \bbeta_k^*$. If $|\eta_k-\eta_k^*|\leq c$, then 
	\begin{equation*}
	\Big\|\mT-\mT^*\Big\|_F^2\leq 9(1+c)\Big(\sum_{k=1}^K\Big\|\sqrt[3]{\eta_k}\bbeta_k - \sqrt[3]{\eta_k^*}\bbeta_k^*\Big\|_2^2\Big)\Big(\sum_{k=1}^K(\sqrt[3]{\eta_k^*})^4\Big)
	\end{equation*}
\end{Lemma}

Denote $\hat{\mT} = \sum_{k=1}^K \eta_k\bbeta_k^{(T^*)}\circ \bbeta_k^{(T^*)}\circ \bbeta_k^{(T^*)} $. Combing  \eqref{eqn:factor_error_sym} and Lemma \ref{lemma:error_transfer},
we have
\begin{eqnarray*}
	\Big\|\hat{\mT}-\mT^*\Big\|_F^2&\leq& 9(1+\varepsilon_0)\frac{C_3\eta_{\min}^{*-\tfrac{4}{3}}}{8}\frac{\sigma^2s\log p}{n}K\eta_{\max}^{*\tfrac{4}{3}},\\
	&=&\frac{9C_3R}{4}\frac{\sigma^2Ks\log p}{n},
\end{eqnarray*}
with probability at least $1-\cO(TKs/n)$.  By setting $C_1 = 9C_2/4$, we complete the proof of Theorem \ref{thm:symmetric_main}. \hfill $\blacksquare$\\

\subsection{Proofs of Theorems \ref{th:lower_bound} and \ref{th:lower_bound_asymmetric}: Minimax Lower Bounds} \label{proof:thm_minimax}
We first consider the proof for Theorem \ref{th:lower_bound_asymmetric} on non-symmetric tensor estimation. Without loss of generality we assume $p = \max\{p_1, p_2, p_3\} $. We uniformly randomly generate $\{\Omega^{(k, m)}\}_{\substack{m=1,\ldots, M\\k=1,\ldots, K}}$ as $MK$ subsets of $\{1,\ldots, p\}$ with cardinality of $s$. Here $M>0$ is a large integer to be specified later. Then we construct $\{\bbeta^{(k, m)}\}_{\substack{m=1,\ldots, M\\k=1,\ldots, K}}\subseteq \mathbb{R}^{p}$ as
\begin{equation*}
\bbeta^{(k, m)}_j = \left\{\begin{array}{ll}
\sqrt{\lambda}, & \text{ if } j\in \Omega^{(k, m)};\\
0, & \text{ if } j\notin \Omega^{(k, m)}.\\	
\end{array}\right.
\end{equation*}
$\lambda>0$ will also be specified a little while later. Clearly, $\|\bbeta^{(k, m_1)} - \bbeta^{(k, m_2)}\|_2^2 \leq 2s\lambda$ for any $1\leq k \leq K$, $1\leq m_1, m_2\leq M$. Additionally, $|\Omega^{(k, m_1)}\cap \Omega^{(k, m_2)}|$ satisfies the hyper-geometric distribution: 
$\mathbb P\left(\left|\Omega^{(k, m_1)}\cap \Omega^{(k, m_2)}\right| = t\right) = \frac{\binom{s}{t}\binom{p-s}{s-t}}{\binom{p}{s}}$. 

Let 
\begin{equation}\label{eq:def-w}
w^{(k, m_1, m_2)} = \left|\Omega^{(k, m_1)}\cap \Omega^{(k, m_2)}\right|,
\end{equation}
then for any $s/2\leq t \leq s$, 
\begin{equation*}
\begin{split}
& \mathbb{P}\left(w^{(k, m_1, m_2)} = t\right) = \frac{\frac{s\cdots (s-t+1)}{t!} \cdot \frac{(p-s)\cdots (p-2s+t+1)}{(s-t)!}}{\frac{p\cdots (p-s+1)}{s!}} \leq \binom{s}{t}\cdot \left(\frac{s}{p-s+1}\right)^t\\
\leq & 2^s\left(\frac{s}{p-s+1}\right)^t \leq \left(\frac{4s}{p-s+1}\right)^t.
\end{split}
\end{equation*}
Thus, if $\eta >0$, the moment generating function of $w^{(k, m_1, m_2)} - \frac{s}{2}$ satisfies
\begin{equation*}
\begin{split}
& \mathbb{E}\exp\left(\eta \left(w^{(k, m_1, m_2)} - \frac{s}{2}\right)\right) \\
\leq & \exp(0)\cdot \mathbb P\left(w^{(k, m_1, m_2)}\leq \frac{s}{2}\right) + \sum_{t=\lfloor s/2\rfloor +1}^s \exp\left(\eta \left(t-\frac{s}{2}\right)\right)\cdot \mathbb P\left(w^{(k, m_1, m_2)} = t\right)\\
\leq & 1 + \sum_{t=\lfloor s/2\rfloor +1}^{s} \left(4s/(p-s+1)\right)^t\exp\left(\eta (t-s/2)\right)\\
= & 1 + \left(\frac{4s}{p-s+1}\right)^{\lfloor s/2\rfloor+1}\sum_{t=0}^{s-\lfloor s/2\rfloor-1} \left(\frac{4s}{p-s+1}\right)^t\exp\left(\eta (t+\lfloor s/2\rfloor+1-s/2)\right)\\
\overset{(*)}{\leq} & 1 + \left(\frac{4s}{p-s+1}\right)^{s/2}\sum_{t=0}^{s-\lfloor s/2\rfloor-1} \left(\frac{4s e^\eta}{p-s+1}\right)^t =  1 + \left(\frac{4s}{p-s+1}\right)^{s/2}\frac{1 - \left(4s e^{\eta}/(p-s+1) \right)^{s-\lfloor s/2\rfloor}}{1 - 4s e^{\eta}/(p-s+1)}\\
< & 1 + \left(4s/(p-s+1) \right)^{s/2} \frac{1}{1 - 4s/(p-s+1)\cdot e^\eta}.
\end{split}
\end{equation*}
Here, (*) is due to $\eta>0$ and $\lfloor s/2\rfloor + 1 \geq s/2$. By setting $\eta = \log((p-s+1)/(8s))$, we have
\begin{equation}\label{ineq:intermediate-prob}
\begin{split}
& \mathbb P\left(\sum_{k=1}^K w^{(k, m_1, m_2)} \geq \frac{3sK}{4}\right) = \mathbb P\left(\sum_{k=1}^K w^{(k, m_1, m_2)} - \frac{sK}{2} \geq \frac{sK}{4}\right)\\
\leq & \frac{ \mathbb{E}\exp\left(\eta(\sum_{k=1}^K w^{(k, m_1, m_2)}-\frac{sK}{2})\right)}{\exp\left(\eta\cdot \frac{sK}{4}\right)} = \frac{\prod_{k=1}^K \mathbb{E}\exp\left(\eta(w^{(k, m_1, m_2)}-\frac{s}{2})\right)}{\exp(\eta\cdot \frac{sK}{4})}\\
\leq & \left(1 + (4s/(p-s+1))^{s/2}\cdot 2 \right)^K\exp\left(-\frac{sK}{4}\log\left(\frac{p-s+1}{8s}\right)\right).
\end{split}
\end{equation}
Since $p\geq 20s$ and $s\geq 4$, we have
\begin{equation*}
    \begin{split}
        & \left(1 + 2(4s/(p-s+1))^{s/2}\right)^K \leq \exp\left(K\log\left(1+2\left(\frac{4}{p/s-1}\right)^{s/2}\right)\right)\\
        \leq & \exp\left(K\log\left(1+2\left(\frac{4}{19}\right)^2\right)\right) \leq \exp\left(K\cdot 0.085\right) \leq \exp\left(sK\log(p/s) \cdot 0.0144\right),
    \end{split}
\end{equation*}
\begin{equation*}
    \begin{split}
        & \exp\left(-\frac{sK}{4}\log\left(\frac{p-s+1}{8s}\right)\right) = \exp\left(-\frac{sK \log(p/s)}{4} + \frac{sK}{4}\log(8p/(p-s+1))\right)\\
        \leq & \exp\left(-\frac{sK\log(p/s)}{4}+\frac{sK}{4}\log(8\cdot 19/20)\right)\leq \exp\left(-sK\log(p/s)\cdot 0.08 \right).
    \end{split}
\end{equation*}
Combining the two inequalities above, we have
\begin{equation*}
    \begin{split}
    & \left(1 + (4s/(p-s+1))^{s/2}\cdot 2 \right)^K\exp\left(-\frac{sK}{4}\log\left(\frac{p-s+1}{8s}\right)\right)\\
    \leq & \exp\left(-c_0sK\log(p/s)\right)
    \end{split}
\end{equation*}
for $c_0 = 1/20$.

Next we choose $M = \lfloor\exp(c_0/2\cdot sK\log(p/s))\rfloor$. Note that 
\begin{equation*}
\begin{split}
& \|\bbeta^{(k, m_1)} - \bbeta^{(k, m_2)}\|_2^2 = \lambda\cdot \left(\left|\Omega^{(k, m_1)}\setminus \Omega^{(k, m_2)}\right| + \left|\Omega^{(k, m_2)}\setminus \Omega^{(k, m_1)}\right|\right)\\
= & \lambda \left(\left|\Omega^{(k, m_1)}\right| + \left|\Omega^{(k, m_2)}\right| - 2\left|\Omega^{(k, m_1)}\cap \Omega^{(k, m_2)}\right|\right) \\
=& 2\lambda\left(s - \left|\Omega^{(k, m_1)}\cap \Omega^{(k, m_2)}\right|\right) \overset{\eqref{eq:def-w}}{=} 2\lambda\left(s - w^{(k, m_1, m_2)}\right),
\end{split}
\end{equation*}
then we further have
\begin{equation*}
\begin{split}
& \mP\left(\sum_{k=1}^K\|\bbeta^{(k, m_1)} - \bbeta^{(k, m_2)}\|_2^2 \geq \frac{sK\lambda}{2}, \forall 1\leq m_1 < m_2 \leq M\right)\\
= & \mP\left(\sum_{k=1}^K2\lambda\left(s-w^{(k,m_1,m_2)}\right) \geq \frac{sK\lambda}{2}, \forall 1\leq m_1 < m_2 \leq M\right)\\
= & \mP\left(\sum_{k=1}^Kw^{(k, m_1, m_2)} \leq \frac{3K}{4}, \forall 1\leq m_1 < m_2 \leq M\right)\\
\overset{\eqref{ineq:intermediate-prob}}{\geq} &  1 - \frac{M(M-1)}{2}\exp\left(-c_0sK\log(p/s)\right)\\
> & 1 - M^2 \exp\left(-c_0sK\log(p/s)\right) \geq 0,
\end{split}
\end{equation*}
which means there are positive probability that $\left\{\bbeta^{(k, m)}\right\}_{\substack{k=1,\ldots, K\\m=1,\ldots, M}}$ satisfy
\begin{equation}\label{ineq:beta^m_1-beta^m_2}
\begin{split}
&\frac{sK\lambda}{2} \leq \min_{1\leq m_1< m_2\leq M} \sum_{k=1}^K\left\|\bbeta^{(k, m_1)}-\bbeta^{(k, m_2)}\right\|_2^2 \\
&\leq \max_{1\leq m_1< m_2\leq M} \sum_{k=1}^K\left\|\bbeta^{(k, m_1)}-\bbeta^{(k, m_2)}\right\|_2^2 \leq 2sK\lambda.
\end{split}
\end{equation}
For the rest of the proof, we fix $\left\{\bbeta^{(k, m)}\right\}_{\substack{k=1,\ldots, K\\m=1,\ldots, M}}$ to be the set of vectors satisfying \eqref{ineq:beta^m_1-beta^m_2}.

Next, recall the canonical basis $\be_{k} = (0,\ldots, \overbrace{1}^{k\text{-th}}, 0, \cdots, 0)\in \mathbb{R}^{p}$. Define
$$\mT^{(m)} = \sum_{k=1}^{K} \bbeta^{(k, m)}\circ \be_{k} \circ \be_{k},\quad 1\leq m \leq M.$$ 
For each tensor $\mT^{(m)}$ and $n$ i.i.d. Gaussian sketches $\bu_i, \bv_i, \bw_i\in \mathbb{R}^{p}$, we denote the response
$$\by^{(m)} = \left\{y_i^{(m)}\right\}_{i=1}^n, \quad y^{(m)}_i = \langle \bu_i\circ \bv_i\circ \bw_i, \mT^{(m)}\rangle + \epsilon_i, $$
where $\epsilon_i \overset{iid}{\sim} N(0, \sigma^2), \quad i=1,\ldots, n.$ Clearly, $\left(\by^{(m)}, \bu, \bv, \bw\right)$ follows a joint distribution, which may vary based on different values of $m$. 

In this step, we analyze the Kullback-Leibler divergence between different distribution pairs:
\begin{equation*}
\begin{split}
& D_{KL}\left((\by^{(m_1)}, \bu, \bv, \bw), (\by^{(m_2)}, \bu, \bv, \bw)\right) := \mathbb{E}_{(\by^{(m_1)}, \bu, \bv, \bw)} \log\Big(\frac{p(\by^{(m_1)}, \bu, \bv, \bw)}{p(\by^{(m_2)}, \bu, \bv, \bw)}\Big).\\
\end{split}
\end{equation*}
Note that conditioning on fixed values of $\bu, \bv, \bw$, 
$$y^{(m)}_i\sim N\left(\sum_{k=1}^K(\bbeta^{(k, m)\top}\bu_i )\cdot (\be^{(k)\top} \bv_i) \cdot (\be^{(k)\top}\bw_i), \sigma^2\right).$$ By the KL-divergence formula for Gaussian distribution,
\begin{equation*}
\begin{split}
& \mathbb{E}_{(\by^{(m_1)}, \bu, \bv, \bw)} \left(\frac{p(\by^{(m_1)}, \bu, \bv, \bw)}{p(\by^{(m_2)}, \bu, \bv, \bw)}\Big| \bu, \bv, \bw\right) \\
= & \frac{1}{2} \sum_{i=1}^n\left(\sum_{k=1}^K \left(\left(\bbeta^{(k,m_1)} - \bbeta^{(k,m_2)}\right)^\top \bu_i\right)\left(\be^{(k)\top} \bv_i\right)\left(\be^{(k)\top} \bw_i\right)\right)^2\sigma^{-2}. 
\end{split}
\end{equation*}
Therefore, for any $m_1 \neq m_2$,
\begin{equation*}
\begin{split}
& D_{KL}\left((\by^{(m_1)}, \bu, \bv, \bw), (\by^{(m_2)}, \bu, \bv, \bw)\right)\\
= & \mathbb{E}_{\bu, \bv, \bw} \frac{1}{2} \sum_{i=1}^n\left(\sum_{k=1}^K (\bbeta^{(k,m_1)} - \bbeta^{(k,m_2)})^\top \bu_i)(\be^{(k)\top} \bv_i)(\be^{(k)\top} \bw_i)\right)^2\sigma^{-2}\\
= & \frac{\sigma^{-2}}{2}\sum_{i=1}^n \sum_{k=1}^K \mathbb{E}_{\bu}((\bbeta^{(k,m_1)} - \bbeta^{(k,m_2)})^\top \bu_i)^2\mathbb{E}_{\bv}(\be^{(k)\top} \bv_i)^2\mathbb{E}_{\bw}(\be^{(k)\top} \bw_i)^2\\
= & \frac{n\sigma^{-2}}{2}\sum_{k=1}^K\|\bbeta^{(k, m_1)} - \bbeta^{(k, m_2)}\|_2^2 \leq \sigma^{-2}nKs\lambda. \\
\end{split}
\end{equation*}

Meanwhile, for any $1\leq m_1 < m_2 \leq M$,
\begin{equation*}
\begin{split}
\|\mT^{(m_1)} - \mT^{(m_2)}\|_F = & \left\|\sum_{k=1}^K(\bbeta^{(k, m_1)}-\bbeta^{(k, m_2)})\circ \be^{(k)} \circ \be^{(k)}\right\|_F\\
= & \sqrt{\sum_{k=1}^K\left\|\bbeta^{(k, m_1)}-\bbeta^{(k, m_2)}\right\|_2^2} \overset{\eqref{ineq:beta^m_1-beta^m_2}}{\geq} \sqrt{\frac{sK\lambda}{2}}.
\end{split}
\end{equation*}
By generalized Fano's Lemma (see, e.g., \cite{yu1997assouad}),
\begin{equation*}
\inf_{\hat{\mT}} \sup_{\mT \in \mathcal{F}} \mathbb{E}\|\hat{\mT} - \mT\|_F \geq \sqrt{\frac{sK\lambda}{2}}\left(1 - \frac{\sigma^{-2}nKs\lambda + \log 2}{\log M}\right).
\end{equation*}
Finally we set $\lambda = \frac{c\sigma^2}{n}\log(p/s)$ for some small constant $c>0$, then
\begin{equation*}
\inf_{\hat{\mT}} \sup_{\mT \in \mathcal{F}} \mathbb{E}\|\hat{\mT} - \mT\|_F^2 \geq \left(\inf_{\hat{\mT}} \sup_{\mT \in \mathcal{F}} \mathbb{E}\|\hat{\mT} - \mT\|_F\right)^2 \geq \frac{c\sigma^2 sK\log(p/s)}{n}.
\end{equation*}
which has finished the proof of Theorem \ref{th:lower_bound_asymmetric}.

For the proof for Theorem \ref{th:lower_bound}, without loss of generality we assume $K$ is a multiple of 3. We first partition $\{1,\ldots, p\}$ into two subintervals: $I_1  = \{1,\ldots, p-K/3\}, I_2 = \{p-K/3+1,\ldots, p\}$, randomly generate $\{\Omega^{(k, m)}\}_{\substack{m=1,\ldots, M\\k=1,\ldots, K/3}}$ as $(MK/3)$ subsets of $\{1,\ldots, p-K/3\}$, and construct $\{\bbeta^{(k,m)}\}_{\substack{m=1,\ldots, M\\k=1,\ldots, K}}\subseteq \mathbb{R}^{p-K/3}$ as
\begin{equation*}
\bbeta^{(k, m)} = \left\{\begin{array}{ll}
\sqrt{\lambda}, & \text{if } j \notin \Omega^{(k, m)};\\
0, & \text{if } j \notin \Omega^{(k, m)}.
\end{array}\right.
\end{equation*}
With $M = \exp(csK\log(p/s))$ and similar techniques as previous proof, one can show there exists positive possibility that
\begin{equation*}
\begin{split}
&\frac{sK\lambda}{6} \leq \min_{1\leq m_1 < m_2\leq M}\sum_{k=1}^{K/3}\|\bbeta^{(k, m_1)} - \bbeta^{(k, m_2)}\|_2^2 \\
&\leq \max_{1\leq m_1 < m_2\leq M}\sum_{k=1}^{K/3}\|\bbeta^{(k, m_1)} - \bbeta^{(k, m_2)}\|_2^2 \leq \frac{2sK}{3}\lambda.
\end{split}
\end{equation*}

We then construct the following candidate symmetric tensors by blockwise design,
\begin{equation*}
\cT^{(m)}\in \mathbb{R}^{p\times p\times p}, \quad \left\{\begin{array}{ll}
\cT^{(m)}_{[I_1, I_2, I_2]} = \sum_{k=1}^{K/3}\bbeta^{(k, m)}\circ \be^{(k)}\circ \be^{(k)},\\ 
\cT^{(m)}_{[I_2, I_1, I_2]} = \sum_{k=1}^{K/3}\be^{(k)}\circ \bbeta^{(k, m)}\circ \be^{(k)},\\
\cT^{(m)}_{[I_2, I_2, I_1]} = \sum_{k=1}^{K/3}\be^{(k)}\circ \be^{(k)} \circ \bbeta^{(k, m)},\\
\cT^{(m)}_{[I_1, I_1, I_1]}, \cT^{(m)}_{[I_1, I_1, I_2]}, \cT^{(m)}_{[I_1, I_2, I_1]}, \cT^{(m)}_{[I_2, I_1, I_1]}, \cT^{(m)}_{[I_2, I_2, I_2]} \text{ are all zeros}.
\end{array}\right.
\end{equation*}
Then we can see for any $\bu\in \mathbb{R}^p$,
\begin{equation*}
\langle \cT^{(m)}, \bu\circ\bu\circ\bu \rangle = 3\sum_{k=1}^{K/3}\left(\bbeta^{(k, m)\top}\bu_{I_1}\right)\cdot\left(\be^{(k)\top}\bu_{I_2}\right)^2.
\end{equation*}
The rest of the proof essentially follows from the proof of Theorem \ref{th:lower_bound_asymmetric}.
\hfill $\blacksquare$\\

\subsection{Proof of Theorem \ref{thm:steins_higher}: High-order Stein's Lemma}\label{proof:stein_higher}
The proof of this theorem follows from the one of Theorem 6 in \cite{janzamin2014score}. For the sake of completeness, we restate the detail here. Applying the recursion representation of score function \eqref{eqn:recursive_form}, we have
\begin{eqnarray*}
\mathbb E\Big[G(\bx)\cS_3(\bx)\Big] &=& \mathbb E\Big[G(\bx)\Big(-\cS_2(\bx)\circ \nabla_{\bx}\log p(\bx)-\nabla_{\bx} \cS_2(\bx)\Big)\Big]\\
&=& -\mathbb E\Big[G(\bx)\cS_2(\bx)\circ \nabla_{\bx}\log p(\bx)\Big]-\mathbb E\Big[G(\bx)\nabla_{\bx} \cS_2(\bx)\Big)\Big].
\end{eqnarray*}
Then, we apply the first-order Stein's lemma (see Lemma  \ref{lemma:stein_lemma}) on function $G(\bx)\cS_2(\bx)$ and obtain
\begin{eqnarray*}
    \mathbb E\Big[G(\bx)\cS_3(\bx)\Big] &=& \mathbb E\Big[\nabla_{\bx}\Big(G(\bx)\cS_2(\bx)\Big)\Big]-\mathbb E\Big[G(\bx)\nabla_{\bx} \cS_2(\bx)\Big)\Big]\\
    &=& \mathbb E\Big[\nabla_{\bx}G(\bx)\cS_2(\bx)+\nabla_{\bx}\cS_2(\bx)G(\bx)\Big]-\mathbb E\Big[G(\bx)\nabla_{\bx} \cS_2(\bx)\Big)\Big]\\
    &=&\mathbb E\Big[\nabla_{\bx}G(\bx)\cS_2(\bx)\Big].
\end{eqnarray*}
Repeating the above argument two more times, we reach the conclusion. \hfill $\blacksquare$\\

\section{Proofs of Several Lemmas}\label{sec:proof_lemmas} \label{supp_sec:proof}

\subsection{Proofs of Lemmas \ref{lemma:moment_asy}, \ref{lemma:moment_symmetric}, and \ref{lemma:moment_interaction}: Moment Calculation}
In this subsection, we present the detail proofs of moment calculation, including non-symmetric case, symmetric case, and interaction model. 

\subsubsection{Proof of Lemma \ref{lemma:moment_asy}}\label{proof_moment_asy}
 By the definition of $\{y_i\}$ in \eqref{eq:model_asymmetric} \& \eqref{eq:low_rank_nonsym}, we have
\begin{equation}
\label{equation:non-symmetric}
\begin{split}
&\mathbb E\Big(\frac{1}{n}\sum_{i=1}^n y_i\bu_i\circ\bv_i\circ\bw_i\Big) = \mathbb E\Big(\frac{1}{n}\sum_{i=1}^n\epsilon_i\bu_i\circ \bv_i\circ \bw_i\Big)\\
&+\mathbb E\Big(\frac{1}{n}\sum_{i=1}^n \sum_{k=1}^K\eta_k^*(\bbeta_{1k}^{*\top}\bu_i)(\bbeta_{2k}^{*\top}\bv_i)(\bbeta_{3k}^{*\top}\bw_i)\bu_i\circ\bv_i\circ\bw_i\Big). 
\end{split}
\end{equation}
First, we observe $\mathbb E(\epsilon_i\bu_i\circ \bv_i\circ \bw_i )=0$ due to the independence between $\epsilon_i$ and $\{\bu_i, \bv_i, \bw_i\}$. Then, we consider a single component from a single observation 
$$M=\mathbb E((\bbeta_{1k}^{*\top}\bu_i)(\bbeta_{2k}^{*\top}\bv_i)(\bbeta_{3k}^{*\top}\bw_i)\bu_i\circ\bv_i\circ\bw_i), \ i\in[n], k\in[K].
$$ For notation simplicity, we drop the subscript $i$ for $i$-th observation and $k$ for $k$-th component such that
\begin{eqnarray}
\label{equation:M_nonsymmetric}
M = \mathbb E \Big((\bbeta_1^{*\top}\bu)(\bbeta_2^{*\top}\bv)(\bbeta_3^{*\top}\bw)\bu\circ\bv\circ\bw\Big)\in\mathbb R^{p_1\times p_2\times p_3}.
\end{eqnarray}
Each entry of $M$ can be calculated as follows
\begin{eqnarray*}
	M_{ijk} &=& \mathbb E \Big((\bbeta_1^{*\top}\bu)(\bbeta_2^{*\top}\bv)(\bbeta_3^{*\top}\bw)u_iv_jw_k\Big)\\
	&=&\mathbb E\Big((\beta_{1i}^*u_i+\sum_{m\neq i}\beta_{1m}^*u_m)u_i\Big)\mathbb E\Big((\beta_{2j}^*u_i+\sum_{m\neq j}\beta_{2m}^*v_m)v_j\Big)\\
	&&\times\mathbb E\Big((\beta_{3k}^*w_k+\sum_{m\neq k}\beta_{3m}^*w_m)w_k\Big)\\
	&=& \beta_{1i}^*\beta_{2j}^*\beta_{3k}^*,
\end{eqnarray*}
which implies $M = \bbeta_1\circ\bbeta_2\circ\bbeta_3$. Combining with $n$ observations and $K$ components, we can obtain
\begin{eqnarray*}
	\mathbb E\big(\cT\big) = \frac{1}{n}\sum_{i=1}^n \sum_{k=1}^K\eta_k^* \bbeta_{1k}\circ\bbeta_{2k}\circ\bbeta_{3k}.
\end{eqnarray*}
This finished our proof. \hfill $\blacksquare$\\

\subsubsection{Proof of Lemma \ref{lemma:moment_symmetric}} \label{proof_moment_sym}

In this subsection, we provide an alternative and more direct proof for Lemma \ref{lemma:moment_symmetric}. We consider a similar single component of (\ref{equation:M_nonsymmetric}) but with a symmetric structure, namely, $M_s = \mathbb E \Big((\bbeta^{*\top}\bx)^3\bx\circ\bx\circ\bx\Big)$. Based on the symmetry of both underlying tensor and sketchings, we will verify the following three cases:
\begin{itemize}
	\item When $i=j=k$, then
	\begin{eqnarray*}
		M_{s_{iii}} &=& \mathbb E\Big(\beta_{i}^*x_i+\sum_{m\neq i}\beta_m^*x_m\Big)^3x_i^3\\
		&=&\mathbb E\Big(\beta_{i}^{*3}x_i^3+3\beta_{i}^{*2}x_i^2\big(\sum_{m\neq i}\beta_{m}^*x_m\big)\\
		&&+3\beta_{i}^*x_i\big(\sum_{m\neq i}\beta_m^*x_m\big)^2+\big(\sum_{m\neq i}\beta_m^* x_m\big)^3\Big)x_i^3\\
		&=& 15 \beta_i^{*3}+9\beta_i^*\sum_{m\neq i}\beta_m^{*2} = 9\beta_i^*+6\beta_i^{*3}.
	\end{eqnarray*}
	The last equation is due to $\|\bbeta^*\|_2=1$.
	\item  When $i\neq j\neq k$, then
	\begin{eqnarray*}
		M_{s_{ijk}} &=& \mathbb E\Big(\beta_{i}^*x_i+ \beta_{j}^*x_j + \beta_{k}^*x_k+\sum_{m\neq i,j,k}\beta_m^*x_m\Big)^3x_ix_jx_k\\
		&=& \mathbb E\big(\beta_{i}^*x_i+ \beta_{j}^*x_j + \beta_{k}^*x_k\big)^3x_ix_jx_k\\
		&=& 6\beta_i^*\beta_j^*\beta_k^*.
	\end{eqnarray*}
	\item When $i=j\neq k$, then
	\begin{eqnarray*}
		M_{s_{iik}} &=&  \mathbb E\Big(\beta_{i}^*x_i + \beta_{k}^*x_k+\sum_{m\neq i,k}\beta_m^*x_m\Big)^3x_i^2x_k\\
		&=& 9\beta_i^{*2}\beta_k^*+3\beta_k^{*3}+3\beta_k^*\big(\sum_{m\neq i,k}
		\beta_m^{*2}\big)\\
		&=&  9\beta_i^{*2}\beta_k^*+3\beta_k^*\big(\sum_{m\neq i}\beta_m^{*2}\big)\\
		&=& 3\beta_k^* + 6\beta_i^{*2}\beta_k^*.
	\end{eqnarray*}
\end{itemize}
Therefore, it is sufficient to calculate $M_s$ by
\begin{equation*}
\begin{split}
M_s =& 3\sum_{k=1}^{K}\eta_k^*\Big(\sum_{m=1}^p\bbeta_k^*\circ \be_m\circ \be_m+ \be_m\circ\bbeta_k^*\circ \be_m+  \be_m\circ \be_m\circ\bbeta_k^*\Big)\nonumber\\
&+6\sum_{k=1}^K \eta_k^*\bbeta_k^*\circ\bbeta_k^*\circ\bbeta_k^*.\nonumber
\end{split}
\end{equation*}
The first term is the bias term due to correlations among symmetric sketchings. Denote $M_1 = \frac{1}{n}\sum_{i=1}^ny_i\bx_i$ and note that  $\mathbb E\Big(\frac{1}{n}\sum_{i=1}^ny_i\bx_i\Big)=3\sum_{k=1}^K\eta_k^*\bbeta_k^*$. Therefore, the empirical first-order moment $M_1$ could be used to remove the bias term as follows
\begin{eqnarray*}
	&&\mathbb E \Big(M_s - \sum_{m=1}^{p}\Big(M_1\circ \be_m\circ \be_m+ \be_m\circ M_1\circ \be_m+  \be_m\circ \be_m\circ M_1\Big)\Big)\\
	 &=& 6\sum_{k=1}^K \eta_k^*\bbeta_k^*\circ\bbeta_k^*\circ\bbeta_k^*.
\end{eqnarray*}
This finishes our proof. \hfill $\blacksquare$\\

\subsubsection{Proof of Lemma \ref{lemma:moment_interaction}}\label{proof_moment_inter}

As before, consider a single component first. For notation simplicity, we drop the subscript $l$ for $l$-th observation and $k$ for $k$-th component. Since each component is normalized, the entry-wise expectation of $(\bbeta^{\top}\bx)^3\bx\circ \bx\circ \bx$ can be calculated as
\begin{eqnarray*}
	&&\Big[\mathbb E(\bbeta^{\top}\bx)^3\bx\circ \bx\circ \bx\Big]_{0,0,0} = 3\beta_0-2\beta_0^3\\
	&&\Big[\mathbb E(\bbeta^{\top}\bx)^3\bx\circ \bx\circ \bx\Big]_{0,0,i} = 3\beta_i\\
	&&\Big[\mathbb E(\bbeta^{\top}\bx)^3\bx\circ \bx\circ \bx\Big]_{0,i,i} = 6\beta_0\beta_i^2+3\beta_0\\
	&&\Big[\mathbb E(\bbeta^{\top}\bx)^3\bx\circ \bx\circ \bx\Big]_{0,i,j} = 6\beta_0\beta_i\beta_j\\
	&&\Big[\mathbb E(\bbeta^{\top}\bx)^3\bx\circ \bx\circ \bx\Big]_{i,i,i} = 6\beta_i^3+9\beta_i\\
	&&\Big[\mathbb E(\bbeta^{\top}\bx)^3\bx\circ \bx\circ \bx\Big]_{i,i,j} = 6\beta_i^2\beta_j+3\beta_j\\
	&&\Big[\mathbb E(\bbeta^{\top}\bx)^3\bx\circ \bx\circ \bx\Big]_{i,j,k} = 6\beta_i\beta_j\beta_k.
\end{eqnarray*}
Due to the symmetric structure and non-randomness of first coordinate, there are bias appearing for each entry. For $i,j,k\neq 0$, we could use $ \sum_{m=1}^p(\ba\circ \be_m\circ \be_m+ \be_m\circ \ba \circ \be_m+\be_m\circ \be_m\circ \ba)$ to remove the bias as shown in the previous proof of Lemma \ref{lemma:moment_symmetric}. For the subscript involving 0, the following two calculations work for removing the bias,
\begin{eqnarray*}
	&&\mathbb E\Big(\frac{1}{3}\cT_s - \frac{1}{6}(\sum_{k=1}^p\cT_{s,[k,k,i]} - (p+1)\ba_i)\Big) = \beta_0^2\beta_i.\\
	&&\mathbb E\Big(\frac{1}{2p-2}(\sum_{k=1}^p\cT_{s[0,k,k]}-(p+2)\cT_{s[0,0,0]})\Big) = \beta_0^3.
\end{eqnarray*} 
This ends the proof. \hfill $\blacksquare$\\

\subsection{Proof of Lemma \ref{lemma:psi_moment_bound}}\label{proof_psi_moment_bound}

Recall the $\|X\|_{\psi_{\alpha}}$ is defined in Definition \ref{def:phi_norm}.  Without loss of generality, we assume $\|X_i\|_{\psi_{\alpha}}=1$ and $\mathbb{E}X_i = 0$ throughout this proof. Let $\beta = (\log 2)^{1/\alpha}$ and $Z_i=(|X_i|-\beta)_{+}$, where $(x)_{+} = x$ if $x\geq 0$ and $(x)_{+}=0$ if else. For notation simplicity, we define $ \|X\|_p = (\mathbb E|X|^p)^{1/p}$ for a random variable $X$. The following step is to estimate the moment of linear combinations of variables $\{X_i\}_{i=1}^n$.

According to the symmetrization inequality (e.g., Proposition 6.3 of \cite{ledoux2013probability}), we have 
\begin{equation}\label{ineqn:contration}
    \Big\|\sum_{i=1}^n a_iX_i\Big\|_p\leq 2\Big\|\sum_{i=1}^n a_i\varepsilon_i X_i\Big\|_p = 2\Big\|\sum_{i=1}^na_i\varepsilon_i|X_i|\Big\|_p,
\end{equation}
where $\{\varepsilon_i\}_{i=1}^n$ are independent Rademacher random variables and we notice that $\varepsilon_i X_i$ and $\varepsilon_i|X_i|$ are identically distributed. Moreover, if $|X_i|\geq \beta$, the definition of $Z_i$ implies that $|X_i| = Z_i+\beta$. And if $|X_i|<\beta$, we have $Z_i=0$. Thus, we have $|X_i|\leq Z_i+\beta$ at any time and it leads to
\begin{eqnarray}\label{ineqn:contration_1}
2\Big\|\sum_{i=1}^na_i\varepsilon_i|X_i|\Big\|_p &\leq& 2\Big\|\sum_{i=1}^n a_i\varepsilon_i(\beta+Z_i)\Big\|_p.
\end{eqnarray}
By triangle inequality, 
\begin{eqnarray}\label{ineqn:contration_2}
2\Big\|\sum_{i=1}^n a_i\varepsilon_i(\beta+Z_i)\Big\|_p \leq 2\Big\|\sum_{i=1}^n a_i\varepsilon_iZ_i\Big\|_p + 2\Big\|\sum_{i=1}^n a_i\varepsilon_i\beta\Big\|_p.
\end{eqnarray}

Next, we will bound the second term of the RHS of \eqref{ineqn:contration_2}. In particular, we will utilize Khinchin-Kahane inequality, whose formal statement is included in Lemma \ref{lemma:Khi} for the sake of completeness. From Lemma \ref{lemma:Khi} we have
\begin{eqnarray}\label{ineqn:contration_3}
    \Big\|\sum_{i=1}^n a_i\varepsilon_i\beta\Big\|_p &\leq& \Big(\frac{p-1}{2-1}\Big)^{1/2} \Big\|\sum_{i=1}^n a_i\varepsilon_i\beta\Big\|_2\nonumber\\
    &\leq& \beta\sqrt{p}\Big\|\sum_{i=1}^n a_i\varepsilon_i\Big\|_2.
\end{eqnarray}
Since $\{\varepsilon_i\}_{i=1}^n$ are independent Rademacher random variables, some simple calculations implies
\begin{eqnarray}\label{ineqn:contration_4}
    \Big(\mathbb E\Big(\sum_{i=1}^n \varepsilon_i a_i\Big)^2\Big)^{1/2}&=& \Big(\mathbb E\Big(\sum_{i=1}^n \varepsilon_i^2 a_i^2+2\sum_{1\leq i<j\leq n}\varepsilon_i\varepsilon_j a_i a_j\Big)\Big)^{1/2}\nonumber\\
    &=& \Big(\sum_{i=1}^n a_i^2\mathbb E\varepsilon_i^2 + 2\sum_{1\leq i<j\leq n}a_ia_j\mathbb E\varepsilon_i\mathbb E\varepsilon_j\Big)^{1/2}\nonumber\\
    &=& \Big(\sum_{i=1}^n a_i^2\Big)^{1/2} = \|\ba\|_2.
\end{eqnarray}
Combining inequalities \eqref{ineqn:contration_1}-\eqref{ineqn:contration_4},
\begin{eqnarray}\label{eqn:proof_1}
2\Big\|\sum_{i=1}^na_i\varepsilon_i|X_i|\Big\|_p \leq 2\Big\|\sum_{i=1}^n a_i\varepsilon_iZ_i\Big\|_p + 2\beta \sqrt{p}\|\ba\|_2.
\end{eqnarray}
Let $\{Y_i\}_{i=1}^n$ are independent symmetric random variables satisfying $\mathbb P(|Y_i|\geq t)=\exp(-t^{\alpha})$ for all $t\geq 0$. Then we have 
\begin{equation*}
\begin{split}
\mathbb P(Z_i\geq t) \leq& \mathbb P(|X_i|\geq t+\beta) = \mathbb P\left(\exp(|X_i|^\alpha)\geq \exp((t+\beta)^\alpha)\right) \\
\leq& \mathbb{E} (\exp(|X_i|)^{\alpha})  \cdot \exp(-(t+\beta)^\alpha) \leq 2\exp(-(t+\beta)^{\alpha})\\
\leq& 2\exp(-t^{\alpha}-\beta^{\alpha}) = \mathbb P(|Y_i|\geq t),
\end{split}
\end{equation*}
which implies
\begin{equation}\label{eqn:proof_2}
    \Big\|\sum_{i=1}^n a_i\varepsilon_iZ_i\Big\|_p\leq \Big\|\sum_{i=1}^n a_i \varepsilon_i Y_i\Big\|_p = \Big\|\sum_{i=1}^n a_i  Y_i\Big\|_p,
\end{equation}
since $\varepsilon_i Y_i$ and $Y_i$ have the same distribution due to symmetry. 
Combining \eqref{eqn:proof_1} and \eqref{eqn:proof_2} together, we reach
\begin{equation}\label{eqn:proof_B1}
    \Big\|\sum_{i=1}^n a_iX_i\Big\|_p\leq 2\beta \sqrt{p}\|\ba\|_2 + 2\Big\|\sum_{i=1}^n a_iY_i\Big\|_p.
\end{equation}

For $0<\alpha <1$, it follows Lemma \ref{lemma:alpha2} that
\begin{equation}\label{eqn:proof_B2}
    \Big\|\sum_{i=1}^n a_i Y_i\Big\|_p \leq C_1(\alpha)(\sqrt{p}\|\ba\|_2 + p^{1/\alpha}\|\ba\|_{\infty}),
\end{equation}
where $C_1(\alpha)$ is some absolute constant only depending on $\alpha$.

For $\alpha\geq 1$, we will combine Lemma \ref{lemma:alpha_1} and the method of the integration by parts to pass from tail bound result to moment bound result. Recall that for every non-negative random variable $X$, integration by parts yields the identity 
\begin{equation*}
    \mathbb E X = \int_{0}^{\infty}\mathbb P(X\geq t)dt.
\end{equation*}
Applying this to $X = |\sum_{i=1}^n a_i Y_i|^p$ and changing the variable $t = t^p$, then we have 
\begin{eqnarray}\label{ineqn:integral}
    \mathbb E|\sum_{i=1}^n a_iY_i|^p &=& \int_{0}^{\infty} \mathbb P\Big(|\sum_{i=1}^n a_iY_i|\geq t\Big) p t^{p-1}dt\nonumber\\
    &\leq &\int_0^{\infty} 2\exp\Big(-c\min\Big(\frac{t^2}{\|\ba\|_2^2}, \frac{t^{\alpha}}{\|\ba\|_{\alpha^*}^{\alpha}}\Big)\Big)pt^{p-1}dt,
\end{eqnarray}
where the inequality is from Lemma \ref{lemma:alpha_1} for all $p\geq 2$ and $1/\alpha+1/\alpha^* = 1$.
In this following, we bound the integral in three steps:
\begin{enumerate}
    \item If $\frac{t^2}{\|\ba\|_2^2}\leq \frac{t^{\alpha}}{\|\ba\|_{\alpha^*}^{\alpha}}$, \eqref{ineqn:integral} reduces to
    \begin{equation*}
          \mathbb E|\sum_{i=1}^n a_iY_i|^p\leq 2p\int_{0}^{\infty}\exp\Big(-c\frac{t^2}{\|\ba\|_2^2}\Big)\Big)t^{p-1}dt.
    \end{equation*}
    Letting $t' = ct^2/\|\ba\|_2^2$, we have 
    \begin{eqnarray*}
        2p\int_{0}^{\infty}\exp\Big(-c\frac{t^2}{\|\ba\|_2^2}\Big)\Big)t^{p-1}dt& =& \frac{p\|\ba\|_2^p}{c^{p/2}}\int_0^{\infty}e^{-t'}t'^{p/2-1}dt'\\
        &=&\frac{p\|\ba\|_2^p}{c^{p/2}} \Gamma(\frac{p}{2})\leq \frac{p\|\ba\|_2^p}{c^{p/2}} (\frac{p}{2})^{p/2},
    \end{eqnarray*}
    where the second equation is from the density of Gamma random variable. Thus,
    \begin{equation}\label{ineqn:moment1}
        \Big(\mathbb E|\sum_{i=1}^n a_iY_i|^p\Big) ^{\tfrac{1}{p}} \leq \frac{p^{1/p}}{(2c)^{1/2}}\sqrt{p}\|\ba\|_2 \leq \frac{\sqrt{2}}{\sqrt{c}}\sqrt{p}\|\ba\|_2.
    \end{equation}

    \item If $\frac{t^2}{\|\ba\|_2^2}> \frac{t^{\alpha}}{\|\ba\|_{\alpha^*}^{\alpha}}$, \eqref{ineqn:integral} reduces to
    \begin{equation*}
         \mathbb E|\sum_{i=1}^n a_iY_i|^p\leq 2p\int_{0}^{\infty}\exp\Big(-c\frac{t^{\alpha}}{\|\ba\|_{\alpha^*}^{\alpha}}\Big)\Big)t^{p-1}dt.
    \end{equation*}
    Letting $t' = ct^{\alpha}/\|\ba\|_{\alpha^*}^{\alpha}$, we have
    \begin{eqnarray*}
        2p\int_{0}^{\infty}\exp\Big(-c\frac{t^{\alpha}}{\|\ba\|_{\alpha^*}^{\alpha}}\Big)\Big)t^{p-1}dt& =& \frac{2p\|\ba\|_{\alpha^*}^p}{\alpha c^{p/{\alpha}}}\int_0^{\infty}e^{-t'}t'^{p/{\alpha}-1}dt'\\
        &=&\frac{2}{\alpha}\frac{p\|\ba\|_{\alpha^*}^p}{c^{p/{\alpha}}} \Gamma(\frac{p}{{\alpha}})\leq \frac{2}{\alpha}\frac{p\|\ba\|_{\alpha^*}^p}{c^{p/{\alpha}}} (\frac{p}{{\alpha}})^{p/{\alpha}}.
    \end{eqnarray*}
    Thus, 
     \begin{equation}\label{ineqn:moment2}
        \Big(\mathbb E|\sum_{i=1}^n a_iY_i|^p\Big) ^{\tfrac{1}{p}} \leq \frac{{2p}^{1/p}}{(c\alpha)^{1/\alpha}}p^{1/\alpha}\|\ba\|_{\alpha^*} \leq \frac{4}{(c\alpha)^{1/\alpha}}p^{1/\alpha}\|\ba\|_{\alpha^*} .
    \end{equation}
    \item Overall, we have the following by combining \eqref{ineqn:moment1} and \eqref{ineqn:moment2},
    \begin{eqnarray*}
    \Big(\mathbb E|\sum_{i=1}^n a_iY_i|^p\Big) ^{\tfrac{1}{p}}\leq \max\Big(\sqrt{\frac{2}{c}},\frac{4}{(c\alpha)^{1/\alpha}}\Big)\Big(\sqrt{p}\|\ba\|_2 + p^{1/\alpha}\|\ba\|_{\alpha^*} \Big).
    \end{eqnarray*}
    After denoting $C_2(\alpha) = \max\Big(\sqrt{\frac{2}{c}},\frac{4}{(c\alpha)^{1/\alpha}}\Big)$, we reach
    \begin{equation}\label{eqn:proof_B3}
         \Big\|\sum_{i=1}^n a_i Y_i\Big\|_p \leq C_2(\alpha)\Big(\sqrt{p}\|\ba\|_2 + p^{1/\alpha}\|\ba\|_{\alpha^*} \Big).
    \end{equation}
\end{enumerate}
Since $0 < \beta < 1$, the conclusion can be reached by combining \eqref{eqn:proof_B1},\eqref{eqn:proof_B2} and \eqref{eqn:proof_B3}. \hfill $\blacksquare$\\

\subsection{Proof of Lemma \ref{lemma:perturbation_analysis}}\label{proof_perturbation}

Firstly, let us consider the non-symmetric perturbation error analysis. According to Lemma \ref{lemma:moment_asy}, the exact form of $\cE=\cT- \mathbb E(\cT)$ is given by
\begin{eqnarray*}
	\cE = \frac{1}{n}\sum_{i=1}^n y_i\bu_i\circ\bv_i\circ \bw_i - \sum_{k=1}^K\eta_k^*\bbeta_{1k}^*\circ\bbeta_{2k}^*\circ\bbeta_{3k}^*.
\end{eqnarray*}
We decompose it by a concentration term $(\cE_1)$ and a noise term $(\cE_2)$ as follows,
\begin{equation*}
\begin{split}
	\cE =& \underbrace{\frac{1}{n}\sum_{i=1}^n\langle\bu_i\circ\bv_i\circ\bw_i, \sum_{k=1}^K\eta_k^*\bbeta_{1k}^*\circ\bbeta_{2k}^*\circ\bbeta_{3k}^*\rangle\bu_i\circ\bv_i\circ\bw_i-\sum_{k=1}^K\eta_k^*\bbeta_{1k}^*\circ\bbeta_{2k}^*\circ\bbeta_{3k}^*}_{\cE_1}\\
	&+\underbrace{\frac{1}{n}\sum_{i=1}^n\epsilon_i \bu_i\circ \bv_i\circ \bw_i}_{\cE_2}.
	\end{split}
\end{equation*}

\noindent\textbf{Bounding $\cE_1$:} For $k$-th componet of $\cE_1$, we denote 
$$\cE_{1k} = \frac{1}{n}\sum_{i=1}^n\langle\bu_i\circ\bv_i\circ\bw_i, \bbeta_{1k}^*\circ\bbeta_{2k}^*\circ\bbeta_{3k}^*\rangle\bu_i\circ\bv_i\circ\bw_i - \bbeta_{1k}^*\circ\bbeta_{2k}^*\circ\bbeta_{3k}^*.
$$
By using Lemma \ref{lemma:tensor_spectral} and $s\leq d\leq Cs$, it suffices to have for some absolute constant $C_{11}$,
$$
\|\cE_{1k}\|_{s+d} \leq  C_{11}\delta_{n,p,s}, \ \text{where} \ \delta_{n,p,s} = (\log n)^3\Big(\sqrt{\frac{s^3\log^3(p/s)}{n^2}} + \sqrt{\frac{s\log (p/s)}{n}}\Big),
$$ 
with probability at least $1-10/n^3$, where $\|\cdot\|_{s+d}$ is the sparse tensor spectral norm defined in \eqref{def:sparse_spectral}. Equipped with the triangle inequality, the sparse tensor spectral norm for $\cE_1$ can be bounded by 
\begin{equation}
\label{eqn:epsilon_1}
\|\cE_1\|_{s+d} \leq  C_{11}\delta_{n,p,s}\sum_{k=1}^K\eta_k^*,
\end{equation}
with probability at least $1-10K/n^3$.

\noindent\textbf{Bounding $\cE_2$:} Note that the random noise $\{\epsilon_i\}_{i=1}^n$ is independent of sketching vector $\{\bu_i,\bv_i, \bw_i\}$. For fixed $\{\epsilon_i\}_{i=1}^n$, applying Lemma \ref{lemma:tensor_concentration}, we have for some absolute constant $C_{12}$
\begin{eqnarray*}
	\Big\|\frac{1}{n}\sum_{i=1}^n\epsilon_i \bu_i\circ \bv_i\circ \bw_i\Big\|_{s+d} \leq C_{12}\|\bepsilon\|_{\infty} C_{11}\delta_{n,p,s},
\end{eqnarray*}
with probability at least $1-1/p$. According to Lemma \ref{lemma:sub_exp}, we have
\begin{equation}
\label{eqn:epsilon_2}
\mathbb P\Big(\|\cE_2\|_{s+d}\geq C_{12}\sigma\log n\delta_{n,p,s}\Big)\leq  \frac{1}{p}+\frac{3}{n}\leq \frac{4}{n}.
\end{equation}

\noindent\textbf{Bounding $\cE$:}
Putting \eqref{eqn:epsilon_1} and \eqref{eqn:epsilon_2} together, we obtain
\begin{eqnarray*}
	\|\cE\|_{s+d}\leq\Big(C_{11}\sum_{k=1}^K\eta_k^* + C_{12}\sigma\log n\Big)\delta_{n,p,s},
\end{eqnarray*}
with probability at least $1-5/n$. Under Condition \ref{con:noise_non}, we have
\begin{eqnarray*}
	\|\cE\|_{s+d} \leq 2C_1 \sum_{k=1}^K\eta_k^* \delta_{n,p,s}\log n,
\end{eqnarray*}
with probability at least $1-5/n$. 

The perturbation error analysis for the symmetric tensor estimation model and the interaction effect model is similar since the empirical first-order moment converges much faster than the empirical third-order moment. So we omit the detailed proof here. \hfill $\blacksquare$\\

\subsection{Proof of Lemma \ref{lemma:onestep_update_sym}}\label{proof:onestep_update_sym}

Lemma \ref{lemma:onestep_update_sym} quantifies one step update for thresholded gradient update. The proof consists of two parts. 

First, we evaluate an oracle estimator $\{\tilde{\bbeta}_k^{(t+1)}\}_{k=1}^K$ with known support information, which is defined as
\begin{eqnarray}\label{def:oracle_est}
\tilde{\bbeta}_k^{(t+1)} = \varphi_{\frac{\mu}{\phi}h(\bbeta_k^{(t)})}\Big(\bbeta_k^{(t)}-\frac{\mu}{\phi}\nabla_k\cL(\bbeta_k^{(t)})_{F^{(t)}}\Big).
\end{eqnarray}
Here,
\begin{itemize}
	\item $h(\bbeta_k^{(t)})$ is the $k$-th component of $h(\bB^{(t)})$ defined in \eqref{eqn:threshold_level_symmetric}.
	\item $\nabla_{\bB}\cL(\bB) = (\nabla_1\cL(\bbeta_1), \cdots, \nabla_K\cL(\bbeta_K)).$
	\item $F^{(t)}= \cup_{k=1}^K F_k^{(t)}$, where $F_k^{(t)} = \supp(\bbeta_k^*)\cup \supp(\bbeta_k^{(t)})$. 
	\item For a vector $\bx\in\mathbb R^p$ and a subset $A\subset \{1,\ldots, p\}$, we denote $\bx_{A}\in \mathbb R^p$ by keeping the coordinates of $\bx$ with indices in $A$ unchanged, while changing all other components to zero.
\end{itemize}
We will show that $\tilde{\bbeta}_k^{(t+1)}$ converges as a geometric rate for optimization error and an optimal rate for statistical error. See Lemma \ref{lemma:oracle_update_sym} for details. 

Second, we aim to prove that $\tilde{\bbeta}_k^{(t+1)}$ and $\bbeta_k^{(t+1)}$ are almost equivalent with high probability. See Lemma \ref{lemma:oracle_equivalence_sym} for details. For simplicity, we drop the superscript of $\bbeta_k^{(t)}, F^{(t)}$ in the following proof, and denote $\tilde{\bbeta}_k^{(t+1)}$, $\bbeta_k^{(t+1)}$ and $F^{(t+1)}$ by $\tilde{\bbeta}_k^+$, $\tilde{\bbeta}_k^+$ and $F^+$, respectively.

\begin{Lemma}
	\label{lemma:oracle_update_sym}
	Suppose Conditions \ref{con:ident}-\ref{con:sample} hold. Assume \eqref{con:initial_input} is satisfied and $|F|\lesssim Ks$. As long as the step size $\mu\leq 32 R^{-20/3}/(3K[220+270K]^2)$, we obtain the upper bound for $\{\tilde{\bbeta}_k^+\}$,
	\begin{equation}
	\begin{split}
	\sum_{k=1}^K\Big\|\sqrt[3]{\eta_k}\tilde{\bbeta}_k^+-\sqrt[3]{\eta_k^*}\bbeta_k^*\Big\|_2^2\leq&\Big(1-32\mu \frac{R^{-\tfrac{8}{3}}}{K^2}\Big)\sum_{k=1}^K\Big\|\sqrt[3]{\eta_k}\bbeta_k-\sqrt[3]{\eta_k^*}\bbeta_k^*\Big\|_2^2\\
	&+2C_3\mu^2R^{-\tfrac{8}{3}}\eta_{\min}^{*-\tfrac{4}{3}}\frac{\sigma^2K^{-2}s\log p}{n},
	\end{split}
	\end{equation}
	with probability at least $1-(21K^2+11K+4Ks)/n$. 
	\end{Lemma}


The proof of Lemma \ref{lemma:oracle_update_sym} is postponed to the Section \ref{proof:one-step_sym}. Next lemma guarantees that with high probability, $\{\bbeta_k^+\}_{k=1}^K$ is equivalent to the oracle update $\{\tilde{\bbeta}_k^+\}_{k=1}^K$ with high probability. 
\begin{Lemma}
	\label{lemma:oracle_equivalence_sym}
	Recall that the truncation level $h(\bbeta_k)$ is defined as 
	\begin{equation}\label{def:truncation_level}
	h(\bbeta_k) = \frac{\sqrt{4\log np}}{n}\sqrt{\sum_{i=1}^n\Big(\sum_{k=1}^K\eta_k(\bx_i^{\top}\bbeta_k)^3-y_i\Big)^2\Big(\eta_k(\bx_i^{\top}\bbeta_k)^2\Big)^2}.
	\end{equation}	 
	If $|F|\lesssim Ks$, we have $\bbeta_k^{+}=\tilde{\bbeta}_k^+$ for any $k\in[K]$ with probability at least $1-(n^2p)^{-1}$ and $F^+\subset F$.
\end{Lemma}

The proof of Lemma \ref{lemma:oracle_equivalence_sym} is postponed to the Section \ref{proof:one-step_sym}. By using Lemma \ref{lemma:oracle_equivalence_sym} and induction, we have 
\begin{eqnarray*}
	F^{(t+1)}\subset \cdots F^{(1)}\subset F^{(0)} = \cup_{k=1}^K\supp(\bbeta_k^*)\cup\supp(\bbeta_k^{(0)}).
\end{eqnarray*} 
It implies for every $t$, we have $|F^{(t)}|\lesssim Ks$. Combining with Lemmas \ref{lemma:oracle_update_sym} and \ref{lemma:oracle_equivalence_sym} together, we obtain with probability at least $1-(21K^2+11K+4Ks)/n$,
\begin{equation}
\begin{split}
\sum_{k=1}^K\Big\|\sqrt[3]{\eta_k}\bbeta_k^+-\sqrt[3]{\eta_k^*}\bbeta_k^*\Big\|_2^2\leq &\Big(1-32\mu K^{-2}R^{-\tfrac{8}{3}}\Big)\sum_{k=1}^K\Big\|\sqrt[3]{\eta_k}\bbeta_k-\sqrt[3]{\eta_k^*}\bbeta_k^*\Big\|_2^2\\
&+2C_3\mu^2R^{-\tfrac{8}{3}}\eta_{\min}^{*-\tfrac{4}{3}}\frac{\sigma^2K^{-2}s\log p}{n},
\end{split}
\end{equation}
This ends the proof.
\hfill $\blacksquare$\\

\subsection{Proof of Lemma \ref{lemma:error_transfer}}\label{proof:error_transfer}

 Based on the CP low-rank structure of true tensor parameter $\mT^*$, we can explicitly write down the distance between $\mT$ and $\mT^*$ under tensor Frobenius norm as follows
\begin{eqnarray*}
	\Big\|\mT - \mT^*\Big\|_F^2 = \sum_{i_1, i_2, i_3}\Big(\sum_{k=1}^{K}\eta_k\beta_{ki_1}\beta_{ki_2}\beta_{ki_3}-\sum_{k=1}^K\eta_k^*\beta_{ki_1}^*\beta_{ki_2}^*\beta_{ki_3}^*\Big)^2.
\end{eqnarray*}
For notation simplicity, denote $\bar{\bbeta}_k=\sqrt[3]{\eta_k}\bbeta_k, \bar{\bbeta}_k^*=\sqrt[3]{\eta_k^*}\bbeta_k^*$. Then
\begin{equation*}
\begin{split}
	\Big\|\mT - \mT^*\Big\|_F^2=&\sum_{i_1, i_2, i_3}\Big(\sum_{k=1}^{K}\bar{\beta}_{ki_1}\bar{\beta}_{ki_2}\bar{\beta}_{ki_3}-\sum_{k=1}^K\bar{\beta}_{ki_1}^*\bar{\beta}_{ki_2}^*\bar{\beta}_{ki_3}^*\Big)^2\\
	=&\sum_{i_1, i_2, i_3}\Big(\sum_{k=1}^{K}(\bar{\beta}_{ki_1}-\bar{\beta}_{ki_1}^*)\bar{\beta}_{ki_2}^*\bar{\beta}_{ki_3}^*+\sum_{k=1}^{K}\bar{\beta}_{ki_1}(\bar{\beta}_{ki_2}-\bar{\beta}_{ki_2}^*)\bar{\beta}_{ki_3}^*\\
	&+\sum_{k=1}^{K}\bar{\beta}_{ki_1}\bar{\beta}_{ki_2}(\bar{\beta}_{ki_3}-\bar{\beta}_{ki_3}^*)\Big)^2 = \text{RHS}.
	\end{split}
\end{equation*}
Since $(a+b+c)^2\leq 3(a^2+b^2+c^2)$, we have 
\begin{equation*}
\begin{split}
	\text{RHS}\leq& 3\sum_{i_1, i_2, i_3}\Big[(\sum_{k=1}^{K}(\bar{\beta}_{ki_1}-\bar{\beta}_{ki_1}^*)\bar{\beta}_{ki_2}^*\bar{\beta}_{ki_3}^*)^2+(\sum_{k=1}^{K}\bar{\beta}_{ki_1}(\bar{\beta}_{ki_2}-\bar{\beta}_{ki_2}^*)\bar{\beta}_{ki_3}^*)^2\\
	&+(\sum_{k=1}^{K}\bar{\beta}_{ki_1}\bar{\beta}_{ki_2}(\bar{\beta}_{ki_3}-\bar{\beta}_{ki_3}^*))^2\Big].
	\end{split}
\end{equation*}
Equipped with Cauchy-Schwarz inequality, RHS can be further bounded by
\begin{eqnarray*}
	\text{RHS}
	&\leq& 3\sum_{i_1,i_2,i_3}\Big[\sum_{k=1}^{K}(\bar{\beta}_{ki_1}-\bar{\beta}_{ki_1}^*)^2\sum_{k=1}^K\bar{\beta}_{ki_2}^{*2}\bar{\beta}_{ki_3}^{*2}\\
	&&+\sum_{k=1}^{K}(\bar{\beta}_{ki_2}-\bar{\beta}_{ki_2}^*)^2\sum_{k=1}^K\bar{\beta}_{ki_1}^{2}\bar{\beta}_{ki_3}^{*2}\\
	&&+\sum_{k=1}^{K}(\bar{\beta}_{ki_3}-\bar{\beta}_{ki_3}^*)^2\sum_{k=1}^K\bar{\beta}_{ki_2}^{2}\ \bar{\beta}_{ki_1}^{2}\Big]
\end{eqnarray*}
At the same time, using $\eta_k\leq (1+c)\eta_k^*$ for $k\in[K]$,
\begin{eqnarray*}
	\Big\|\mT - \mT^*\Big\|_F^2&\leq& 3\Big[\sum_{i_1=1}^p\sum_{k=1}^{K}(\bar{\beta}_{ki_1}-\bar{\beta}_{ki_1}^*)^2(\sum_{i_2=1}^p\sum_{i_3=1}^p\sum_{k=1}^K\bar{\beta}_{ki_2}^{*2}\bar{\beta}_{ki_3}^{*2})\\
	&&+\sum_{i_2=1}^p\sum_{k=1}^{K}(\bar{\beta}_{ki_2}-\bar{\beta}_{ki_2}^*)^2(\sum_{i_1=1}^p\sum_{i_3=1}^p\sum_{k=1}^K\bar{\beta}_{ki_1}^{2}\bar{\beta}_{ki_3}^{*2})\\
	&&+\sum_{i_3=1}^p\sum_{k=1}^{K}(\bar{\beta}_{ki_3}-\bar{\beta}_{ki_3}^*)^2(\sum_{i_2=1}^p\sum_{i_1=1}^p\sum_{k=1}^K\bar{\beta}_{ki_2}^{2}\bar{\beta}_{ki_1}^{2})\Big]\\
	&=&3\Big(\sum_{k=1}^K\|\bar{\bbeta}_k-\bar{\bbeta}_k^*\|_2^2\Big)\Big(\sum_{k=1}^K(\sqrt[3]{\eta_k^*})^4+\sum_{k=1}^K(\sqrt[3]{\eta_k^*})^2(\sqrt[3]{\eta_k})^2+\sum_{k=1}^K(\sqrt[3]{\eta_k})^4\Big)\\
	&\leq& 9(1+c) \Big(\sum_{k=1}^K\|\bar{\bbeta}_k-\bar{\bbeta}_k^*\|_2^2\Big)\Big(\sum_{k=1}^K(\sqrt[3]{\eta_k^*})^4\Big).
\end{eqnarray*}

For the non-symmetric tensor estimation model, we have
\begin{eqnarray*}
	\Big\|\mT - \mT^*\Big\|_F^2 = \sum_{i_1, i_2, i_3}\Big(\sum_{k=1}^{K}\eta_k\beta_{1ki_1}\beta_{2ki_2}\beta_{3ki_3}-\sum_{k=1}^K\eta_k^*\beta_{1ki_1}^*\beta_{2ki_2}^*\beta_{3ki_3}^*\Big)^2.
\end{eqnarray*}
Following the same strategy above, we obtain
\begin{equation*}
\begin{split}
	\Big\|\mT - \mT^*\Big\|_F^2\leq& 3(1+c) \Big(\sum_{k=1}^K\|\bar{\bbeta}_{1k}-\bar{\bbeta}_{1k}^*\|_2^2+\sum_{k=1}^K\|\bar{\bbeta}_{2k}-\bar{\bbeta}_{2k}^*\|_2^2\\
	&+\sum_{k=1}^K\|\bar{\bbeta}_{3k}-\bar{\bbeta}_{3k}^*\|_2^2\Big)\Big(\sum_{k=1}^K(\sqrt[3]{\eta_k^*})^4\Big).
\end{split}
\end{equation*}
This ends the proof. \hfill $\blacksquare$\\

\subsection{Proof of Lemma \ref{lemma:oracle_update_sym}}\label{proof:one-step_sym}

First of all, we state a lemma to illustrate the effect of weight $\phi$. The proof of Lemma \ref{lemma:marginal_effect} is deferred to Section \ref{lemma:marginal_effect}. 
\begin{Lemma}\label{lemma:marginal_effect}
	Consider $\{y_i\}_{i=1}^n$ come from either non-symmetric tensor estimation model \eqref{eq:model_asymmetric} or symmetric tensor estimation model \eqref{eq:model}. Suppose Conditions \ref{con:incoherence_sym}-\ref{con:sample} hold. Then $\phi=\frac{1}{n} \sum_{i=1}^n y_i^2$ is upper and lower bounded by  
	\begin{equation*}
	(16-6\Gamma^3-9\Gamma)(\sum_{k=1}^K\eta_k^*)^2\leq \frac{1}{n}\sum_{i=1}^{n}y_i^2\leq (16+6\Gamma^3+9\Gamma)(\sum_{k=1}^K\eta_k^*)^2, 
	\end{equation*}
	with probability at least $1-(K^2+K+3)/n$, where $\Gamma$ is the incoherence parameter defined in Definition \ref{con:incoherence_sym}.
\end{Lemma}

According to Lemma \ref{lemma:marginal_effect}, $\frac{1}{n}\sum_{i=1}^ny_i^2$ approximates $(\sum_{k=1}^K\eta_k^*)^2$ up to some constants with high probability. Moreover, we know that from \eqref{con:initial_input}, $\max_{k}|\eta_k-\eta_k^*|\leq \varepsilon_0$ for some small $\varepsilon_0$. Based on those two facts described above, we replace $\eta_k$ by $\eta_k^*$ and $\phi$ by $(\sum_{k=1}^K\eta_k^*)^{2}$ for the sake of completeness. Note that this change could only result in some constant scale changes for final results. Similar simplification was used in matrix recovery scenario \cite{Tu2015}. Therefore, we define the weighted estimator and weighted true parameter as $\bar{\bbeta}_k = \sqrt[3]{\eta_k^*}\bbeta_k$, $\bar{\bbeta}_k^* = \sqrt[3]{\eta_k^*}\bbeta_k^*$. Now, $\eta_k^*\beta_k\circ\beta_k\circ\beta_k = \bar{\beta}_k\circ \bar{\beta}_k\circ  \bar\beta_k$. Recall $\mathcal{\cdot}$
is the loss function defined in \eqref{eq:L_s-T}. Correspondingly with a slight abuse of notation, define the gradient function $\nabla_k\cL(\bar{\bbeta}_k)$ on $F$ as
\begin{eqnarray*}
	\nabla_k\cL(\bar{\bbeta}_k)_{F}=\frac{6\sqrt[3]{\eta_k^*}}{n}\sum_{i=1}^n\Big(\sum_{k'=1}^K(\bx_{i_{F}}^{\top}\bar{\bbeta}_{k'})^3-y_i\Big)(\bx_{i_{F}}^{\top}\bar{\bbeta}_k)^2\bx_{i_{F}},
\end{eqnarray*}
and its noiseless version as
\begin{equation}
\label{def:noiseless_gradient}
\nabla_k\tilde{\cL}(\bar{\bbeta}_k)_{F}=\frac{6\sqrt[3]{\eta_k^*}}{n}\sum_{i=1}^n\Big(\sum_{k'=1}^K(\bx_{i_{F}}^{\top}\bar{\bbeta}_{k'})^3-\sum_{k'=1}^K(\bx_{i_{F}}^{\top}\bar{\bbeta}^*_{k'})^3\Big)(\bx_{i_F}^{\top}\bar{\bbeta}_k)^2\bx_{i_{F}}.
\end{equation}
According to the definition of thresholding function \eqref{eqn:thresholding}, $\tilde{\bbeta}_k^+$ can be written as 
\begin{equation*}
\tilde{\bbeta}_k^+=\bbeta_k-\frac{\mu}{\phi}\nabla_k\cL(\bar{\bbeta}_k)_{F} +\frac{\mu}{\phi}h(\bar{\bbeta}_k)\bgamma_k,
\end{equation*}
where $\bgamma_k\in \mathbb R^p$ satisfies $\supp(\bgamma_k)\subset F$, $\|\bgamma_k\|_{\infty}\leq 1$ and $h(\bar{\bbeta}_k)$ is defined as
\begin{equation}
\label{def:thresholded_function_sym}
h(\bar{\bbeta}_k) = \frac{\sqrt{4\log(np)}}{n}\sqrt{\sum_{i=1}^n\Big(\sum_{k=1}^K (\bx_{i_{F}}^{\top}\bar{\bbeta}_k)^3-y_i\Big)^2\eta_k^{*\tfrac{2}{3}}(\bx_{i_{F}}^{\top}\bar{\bbeta}_k)^2}.
\end{equation}
Moreover, we denote $\bz_k = \bar{\bbeta}_k-\bar{\bbeta}_k^*$.  With a little abuse of notations, we also drop the subscript $F$ in this section for notation simplicities. 

We expand and decompose the sum of square error by three parts as follows:
\begin{equation}\label{eqn:error_decom}
\begin{split}
&\sum_{k=1}^K\Big\|\sqrt[3]{\eta_k^*}\tilde{\bbeta}_k^+-\sqrt[3]{\eta_k^*}\bbeta^*_k\Big\|_2^2 \\
&= \sum_{k=1}^K\Big\|\bz_k - \frac{\mu\sqrt[3]{\eta_k^*}}{\phi}\nabla_k\cL(\bar{\bbeta}_k) +\frac{\mu\sqrt[3]{\eta_k^*}}{\phi}h(\bar{\bbeta}_k)\bgamma_k\Big\|_2^2\\
&=  \underbrace{\sum_{k=1}^K\Big\|\bz_k-\mu\frac{\sqrt[3]{\eta_k^*}}{\phi}\nabla_k\cL(\bar{\bbeta}_k)\Big\|_2^2}_{\text{A: gradient update effect}}+\underbrace{\sum_{k=1}^K\Big\|\frac{\mu\sqrt[3]{\eta_k^*}}{\phi}h(\bar{\bbeta}_k)\gamma_k\Big\|_2^2}_{\text{B: threshoding effect}}\\
&+\underbrace{\sum_{k=1}^K\Big\langle\bz_k-\mu\frac{\sqrt[3]{\eta_k^*}}{\phi}\nabla_k\cL(\bar{\bbeta}_k), \frac{\mu\sqrt[3]{\eta_k}}{\phi}h(\bar{\bbeta}_k)\gamma_k\Big\rangle}_{\text{C: cross term}}.
\end{split}
\end{equation}
In the following proof, we will bound three parts sequentially. 

\subsubsection{Bounding gradient update effect}
In order to separate the optimization error and statistical error, we use the noiseless gradient $\nabla_k \tilde{\cL}(\bar{\bbeta}_k)$ as a bridge such that $A$ can be decomposed as
\begin{equation}
\label{eqn:symmetric_A}
\begin{split}
A=&\sum_{k=1}^K\|\bz_k\|_2^2- 2\mu\sum_{k=1}^K\Big\langle\frac{\sqrt[3]{\eta_k^*}}{\phi}\nabla_k\cL(\bar{\bbeta}_k), \bz_k\Big\rangle+\mu^2\sum_{k=1}^K\Big\|\frac{\sqrt[3]{\eta_k^*}}{\phi}\nabla_k\cL(\bar{\bbeta}_k)\Big\|_2^2\\
\leq&\sum_{k=1}^K\|\bz_k\|_2^2-2\mu\underbrace{\sum_{k=1}^K\Big\langle\frac{\sqrt[3]{\eta_k^*}}{\phi}\nabla_k\tilde{\cL}(\bar{\bbeta}_k), \bz_k\Big\rangle}_{A_1}+2\mu^2\underbrace{\sum_{k=1}^K\Big\|\frac{\sqrt[3]{\eta_k^*}}{\phi}\nabla_k\tilde{\cL}(\bar{\bbeta}_k)\Big\|_2^2}_{A_2}\\
& + 2\mu^2\underbrace{\sum_{k=1}^K \Big\|\frac{\sqrt[3]{\eta_k^*}}{\phi}\Big(\nabla_k\tilde{\cL}(\bar{\bbeta}_k)-\nabla_k\cL(\bar{\bbeta}_k)\Big)\Big\|_2^2}_{A_3}\\
&+2\mu\underbrace{\sum_{k=1}^K\Big\langle\bz_k, \frac{\sqrt[3]{\eta_k^*}}{\phi}\Big(\nabla_k\tilde{\cL}(\bar{\bbeta}_k)-\nabla_k\cL(\bar{\bbeta}_k)\Big)\Big\rangle}_{A_4},
\end{split}
\end{equation}
where $A_1$ and $A_2$ quantify the optimization error, $A_3$ quantifies the statistical error, and $A_4$ is a cross term which can be negligible comparing with the rate of the statistical error. The lower bound for $A_1$ and upper bound for $A_2$ together coincide with the verification of regularity conditions in the matrix recovery case \cite{CLS15}. 

~\\
\noindent \textbf{Step One: Lower bound for $A_1$.} Plugging in $\phi = (\sum_{k=1}^K\eta_k^*)^{2}$, we have
\begin{eqnarray}\label{ineqn:weight}
K^{-2}R^{-\tfrac{2}{3}}\eta_{\max}^{*-\tfrac{4}{3}}\leq\frac{(\sqrt[3]{\eta_k^*})^2}{\phi} = \frac{(\sqrt[3]{\eta_k^*})^2}{(\sum_{k=1}^K\eta_k^*)^{2}}\leq K^{-2}R^{\tfrac{2}{3}}\eta_{\min}^{*-\tfrac{4}{3}}.
\end{eqnarray}
According to the definition of noiseless gradient $\nabla_k \tilde{\cL}(\bbeta_k)$ and $\bz_k$, $A_1$ can be expanded and decomposed sequentially by nine terms,
\begin{equation}\label{eqn:decomposition}
\begin{split}
A_1 \geq&K^{-2}R^{-\tfrac{2}{3}}\eta_{\max}^{*-\tfrac{4}{3}}\Big[\frac{6}{n}\sum_{i=1}^n\Big(\sum_{k'=1}^K 3(\bx_i^{\top}\bz_{k'})(\bx_i^{\top}\bar{\bbeta}_{k'})^2\sum_{k=1}^K(\bx_i^{\top}\bz_k)(\bx_i^{\top}\bar{\bbeta}_k^*)^2\Big)\Leftarrow \ A_{11}\\
&+\frac{6}{n}\sum_{i=1}^n\Big(\sum_{k'=1}^K 3(\bx_i^{\top}\bz_{k'})(\bx_i^{\top}\bar{\bbeta}_{k'})^2\sum_{k=1}^K2(\bx_i^{\top}\bz_k)^2(\bx_i^{\top}\bar{\bbeta}_k^*)\Big)\Leftarrow \ A_{12}\\
&+\frac{6}{n}\sum_{i=1}^n\Big(\sum_{k'=1}^K 3(\bx_i^{\top}\bz_{k'})(\bx_i^{\top}\bar{\bbeta}_{k'})^2\sum_{k=1}^K(\bx_i^{\top}\bz_k)^3\Big)\Leftarrow \ A_{13}\\
&+\frac{6}{n}\sum_{i=1}^n\Big(\sum_{k'=1}^K 3(\bx_i^{\top}\bz_{k'})^2(\bx_i^{\top}\bar{\bbeta}_{k'})\sum_{k=1}^K(\bx_i^{\top}\bz_k)(\bx_i^{\top}\bar{\bbeta}_k^*)^2\Big)\Leftarrow \ A_{14}\\
&+\frac{6}{n}\sum_{i=1}^n\Big(\sum_{k'=1}^K 3(\bx_i^{\top}\bz_{k'})^2(\bx_i^{\top}\bar{\bbeta}_{k'})\sum_{k=1}^K2(\bx_i^{\top}\bz_k)^2(\bx_i^{\top}\bar{\bbeta}_k^*)\Big)\Leftarrow \ A_{15}\\
&+\frac{6}{n}\sum_{i=1}^n\Big(\sum_{k'=1}^K 3(\bx_i^{\top}\bz_{k'})^2(\bx_i^{\top}\bar{\bbeta}_{k'})\sum_{k=1}^K(\bx_i^{\top}\bz_k)^2(\bx_i^{\top}\bar{\bbeta}_k^*)\Big)\Leftarrow \ A_{16}\\
&+\frac{6}{n}\sum_{i=1}^n\Big(\sum_{k'=1}^K 3(\bx_i^{\top}\bz_{k'})^3\sum_{k=1}^K(\bx_i^{\top}\bz_k)(\bx_i^{\top}\bar{\bbeta}_k^*)^2\Big)\Leftarrow \ A_{17}\\
&+\frac{6}{n}\sum_{i=1}^n\Big(\sum_{k'=1}^K 3(\bx_i^{\top}\bz_{k'})^3\sum_{k=1}^K2(\bx_i^{\top}\bz_k)^2(\bx_i^{\top}\bar{\bbeta}_k^*)\Big)\Leftarrow \ A_{18}\\
&+\frac{6}{n}\sum_{i=1}^n\Big(\sum_{k'=1}^K 3(\bx_i^{\top}\bz_{k'})^3\sum_{k=1}^K(\bx_i^{\top}\bz_k)^3\Big)\Big]\Leftarrow \ A_{19},
\end{split}
\end{equation}
where $A_{11}$ is the main term according to the order of $\bar{\bbeta}_k^*$, while $A_{12}$ to $A_{19}$ are remainder terms. The proof of lower bound for $A_{11}$ to $A_{19}$ follows two steps:
\begin{enumerate}
	\item Calculate and lower bound the expectation of each term through Lemma \ref{cor:expectation_gaussian_vector}: high-order Gaussian moment;
	\item Argue that the empirical version is concentrated around their expectation with high probability through Lemma \ref{cor:concentr_high_order}: high-order concentration inequality.
\end{enumerate}

\textbf{Bounding $A_{11}$.} Note that $A_{11}$ involves the product of dependent Gaussian vectors. This brings difficulties on both the calculation of expectations and the use of concentration inequality.  According to the high-order Gaussian moment results in Lemma \ref{cor:expectation_gaussian_vector}, the expectation of $A_{11}$ can be calculated explicitly as  
\begin{equation}\label{eqn:expectation_A11}
\begin{split}
\mathbb E(A_{11}) &= 36 \sum_{k=1}^K\sum_{k'=1}^K (\bar{\bbeta}_{k'}^{*\top}\bar{\bbeta}_k^*)^2(\bz_{k'}^{\top}\bz_k) \Leftarrow \ I_{1}\\
&+72 \sum_{k=1}^K\sum_{k'=1}^K (\bar{\bbeta}_{k'}^{*\top}\bar{\bbeta}_k^*)(\bz_{k'}^{\top}\bar{\bbeta}_k^*)(\bz_k^{\top}\bar{\bbeta}_{k'}^*) \Leftarrow \ I_{2}\\
&+108 \sum_{k=1}^K\sum_{k'=1}^K (\bar{\bbeta}_{k'}^{*\top}\bar{\bbeta}_k^*)(\bz_{k'}^{\top}\bar{\bbeta}_{k'}^*)(\bz_k^{\top}\bar{\bbeta}_k^*)\Leftarrow \ I_{3}\\
&+54 \sum_{k=1}^K\sum_{k'=1}^K (\bar{\bbeta}_{k'}^{*\top}\bar{\bbeta}_{k'}^*)(\bar{\bbeta}_{k}^{*\top}\bar{\bbeta}_{k}^*)(\bz_{k'}^{\top}\bz_k)\Leftarrow \ I_{4}.
\end{split}
\end{equation}
Note that $I_1$ to $I_4$ involve the summation of $K^2$ term. To use incoherence Condition  \ref{con:incoherence_sym},  we isolate $K$ terms with $k=k'$. Then, $I_1$ to $I_4$ could be lower bounded as 
\begin{equation*}
\begin{split}
I_1&\geq 36\eta_{\min}^{*4/3}\Big[\sum_{k=1}^K \|\bz_k\|_2^2- \Gamma^2\Big(\sum_{k=1}^K\|\bz_k\|_2\Big)^2\Big]\\
I_2&\geq 72\eta_{\min}^{*4/3}\Big[\sum_{k=1}^K (\bz_k^{\top}\bar{\bbeta}_k^*)^2- \Gamma\Big(\sum_{k=1}^K\|\bz_k\|_2\Big)^2\Big]\\
I_3&\geq 108\eta_{\min}^{*4/3}\Big[\sum_{k=1}^K (\bz_k^{\top}\bar{\bbeta}_k^*)^2- \Gamma\Big(\sum_{k=1}^K\|\bz_k\|_2\Big)^2\Big]\\
I_4&\geq 54\eta_{\min}^{*4/3}\Big\|\sum_{k=1}^{K}\bz_k\Big\|_2^2\geq 0,
\end{split}
\end{equation*}
where $\Gamma$ is the incoherence parameter. Putting the above four bounds together, they jointly provide
\begin{equation}\label{ineqn:S4}
\mathbb E(A_{11})\geq 36\eta_{\min}^{*4/3}\sum_{k=1}^K \|\bz_k\|_2^2 - \Big(36\eta_{\min}^{*4/3}\Gamma^2 +180\eta_{\min}^{*4/3}\Gamma\Big)\Big(\sum_{k=1}^K\|\bz_k\|_2\Big)^2.
\end{equation}

On the other hand, repeatedly using Lemma \ref{cor:concentr_high_order}, we obtain that with probability at least $1-1/n$,
\begin{equation*}
\begin{split}
	&\Big|\frac{1}{n}\sum_{i=1}^n\Big((\bx_i^{\top}\bz_{k'})(\bx_i^{\top}\bar{\bbeta}_{k'}^*)^2(\bx_i^{\top}\bz_k)(\bx_i^{\top}\bar{\bbeta}_k^*)^2-\mathbb E(\bx_i^{\top}\bz_{k'})(\bx_i^{\top}\bar{\bbeta}_{k'}^*)^2(\bx_i^{\top}\bz_k)(\bx_i^{\top}\bar{\bbeta}_k^*)^2\Big)\Big|\\
	&\leq C\frac{(\log n)^3}{\sqrt{n}}(\sqrt[3]{\eta_{\max}^*})^4\|\bz_{k'}\|_2\|\bz_k\|_2.
	\end{split}
\end{equation*}
Taking the summation over $k,k'\in[K]$, it could further imply that for some absolute constant $C$,
\begin{eqnarray} \label{ineqn:S5}
\Big|A_{11}-\mathbb E(A_{11})\Big|\leq 18C\frac{(\log n)^3}{\sqrt{n}}(\sqrt[3]{\eta_{\max}^*})^4\Big(\sum_{k=1}^{K}\|\bz_{k}\|_2\Big)^2,
\end{eqnarray}
with probability at least $1-K^2/n$. Combining \eqref{ineqn:S4} and \eqref{ineqn:S5}, we obtain with probability at least $1-K^2/n$,
\begin{equation}\label{eqn:bound_A11}
\begin{split}
&K^{-2}R^{-\tfrac{2}{3}}\eta_{\max}^{*-\tfrac{4}{3}}A_{11}\\
\geq &\Big[36K^{-2}R^{-\tfrac{8}{3}}-K^{-\tfrac{3}{2}} \Big(216R^{-\tfrac{8}{3}}\Gamma+18C\frac{(\log n)^3}{\sqrt{n}}\Big)\Big]\sum_{k=1}^K \|\bz_k\|_2^2,
\end{split}
\end{equation}
where $R = \eta_{\max}^*/\eta_{\min}^*$. Here, we use the fact $\Gamma \leq 1$ and $(\sum_{k=1}^K\|\bz_k\|_2)^2\leq K(\sum_{k=1}^K\|\bz_k\|_2^2)$.

\textbf{Bounding $A_{12}$ to $A_{19}$:} For remainder terms, we follow the same proof strategy. According to Lemma \ref{cor:expectation_gaussian_vector}, the expectation of $A_{12}$ can be calculated as 
\begin{eqnarray*}
	\mathbb E(A_{12}) &=& 36 \sum_{k=1}^K\sum_{k'=1}^K (\bz_k^{\top}\bar{\bbeta}_{k'}^{*\top})^2(\bz_{k'}^{\top}\bar{\bbeta}_k^*) \Leftarrow \ I_{1}\\
	&&+72 \sum_{k=1}^K\sum_{k'=1}^K (\bz_k^{\top}\bar{\bbeta}_{k'}^*)(\bar{\bbeta}_{k'}^{*\top}\bar{\bbeta}_k^*)(\bz_{k'}^{\top}\bz_k) \Leftarrow \ I_{2}\\
	&&+108 \sum_{k=1}^K\sum_{k'=1}^K (\bz_k^{\top}\bar{\bbeta}_{k'})(\bz_{k'}^{\top}\bar{\bbeta}_{k'}^*)(\bz_k^{\top}\bar{\bbeta}_k^*)\Leftarrow \ I_{3}\\
	&&+54 \sum_{k=1}^K\sum_{k'=1}^K (\bar{\bbeta}_{k'}^{*\top}\bar{\bbeta}_{k'}^*)(\bz_{k'}^{\top}\bar{\bbeta}_{k}^*)(\bz_{k}^{\top}\bz_k)\Leftarrow \ I_{4}.
\end{eqnarray*}
Let us analyze $I_1$ first. Under \eqref{con:initial_input}, $\|\bz_k\|_2\leq \varepsilon_0\sqrt[3]{\eta_k^*}$, it suffices to show that 
\begin{eqnarray*}
	\sum_{k=1}^K\sum_{k'=1}^K(\bz_k^{\top}\bar{\bbeta}_{k'})^2(\bz_{k'}^{\top}\bar{\bbeta}_k^*)&\geq& -\sum_{k=1}^K\sum_{k'=1}^{K}\|\bz_k\|_2^2\|\bar{\bbeta}_{k'}^*\|_2^2\|\bz_k'\|_2\|\bar{\bbeta}_{k}^*\|_2\\
	&\geq& -\eta_{\max}^{*\tfrac{4}{3}}\varepsilon_0\Big(\sum_{k=1}^K\|\bz_k\|_2\Big)^2.
\end{eqnarray*}
This immediately implies a lower bound for $\mathbb E(A_{12})$ after we bound similarly for $I_2,I_3$ and $I_4$,
\begin{equation}
\mathbb E(A_{12})\geq -270\eta_{\max}^{*\tfrac{4}{3}}\varepsilon_0\Big(\sum_{k=1}^K\|\bz_k\|_2\Big)^2.
\end{equation}
By Lemma \ref{cor:concentr_high_order}, we obtain for some absolute constant $C$,
\begin{equation}\label{eqn:bound_A12}
\begin{split}
&K^{-2}R^{-\tfrac{2}{3}}\eta_{\max}^{*-\tfrac{4}{3}}A_{12}\\
 \geq& K^{-2}R^{-\tfrac{2}{3}}\eta_{\max}^{*-\tfrac{4}{3}}\Big[\mathbb E(A_{12}) - 18C\eta_{\max}^{*\tfrac{4}{3}}\varepsilon_0\Big(\sum_{k=1}^K\|\bz_k\|_2\Big)^2\frac{(\log n)^3}{\sqrt{n}}\Big]\\
\geq&-K^{-1}R^{-\tfrac{2}{3}}\varepsilon_0 \Big(270+18C\frac{(\log n)^3}{\sqrt{n}}\Big)\Big(\sum_{k=1}^K\|\bz_k\|_2^2\Big),
\end{split}
\end{equation}
with probability at least $1-K^2/n$. The detail derivation is the same as in \eqref{eqn:bound_A11}, so we omit here.

Similarly, the lower bounds of $A_{13}$ to $A_{19}$ can be derived as follows 
\begin{equation}\label{eqn:boundrest}
\begin{split}
&K^{-\tfrac{1}{2}}\eta_{\max}^{*-\tfrac{4}{3}}A_{14}\geq -K^{\tfrac{1}{2}}\varepsilon_0 \Big(270+18C\frac{(\log n)^3}{\sqrt{n}}\Big)\Big(\sum_{k=1}^K\|\bz_k\|_2^2\Big)\\
&K^{-\tfrac{1}{2}}\eta_{\max}^{*-\tfrac{4}{3}}A_{13}, A_{15}, A_{17}\geq -K^{\tfrac{1}{2}}\varepsilon_0^2 \Big(270+18C\frac{(\log n)^3}{\sqrt{n}}\Big)\Big(\sum_{k=1}^K\|\bz_k\|_2^2\Big)\\
&K^{-\tfrac{1}{2}}\eta_{\max}^{*-\tfrac{4}{3}}A_{16},A_{18}\geq -K^{\tfrac{1}{2}}\varepsilon_0^3 \Big(270+18C\frac{(\log n)^3}{\sqrt{n}}\Big)\Big(\sum_{k=1}^K\|\bz_k\|_2^2\Big)\\
&K^{-\tfrac{1}{2}}\eta_{\max}^{*-\tfrac{4}{3}}A_{19}\geq -K^{\tfrac{1}{2}}\varepsilon_0^4 \Big(270+18C\frac{(\log n)^3}{\sqrt{n}}\Big)\Big(\sum_{k=1}^K\|\bz_k\|_2^2\Big).
\end{split}
\end{equation}
Putting  \eqref{eqn:bound_A11}, \eqref{eqn:bound_A12} and \eqref{eqn:boundrest} together, we have with probability at least $1-9K^2/n$,
\begin{eqnarray*}
	A_1&\geq& \Big[36K^{-2}R^{-\tfrac{8}{3}}-K^{-\tfrac{3}{2}}\Big(2160R^{-\tfrac{3}{3}}\Gamma+18C\frac{(\log n)^3}{\sqrt{n}}\Big) \\
	&&- 8\varepsilon_0K^{-1}R^{-\tfrac{2}{3}}\Big(270+18C\frac{(\log n)^3}{\sqrt{n}}\Big)\Big]\Big(\sum_{k=1}^K\|\bz_k\|_2^2\Big).
\end{eqnarray*}
For the above bound,
\begin{itemize}
	\item When the sample size satisfies $n\geq (18CK^{1/2}R^{8/3}(\log n)^3)^2$, we have
	\begin{equation*}
	\max\Big\{ 18K^{-\tfrac{3}{2}}C\frac{(\log n)^3}{\sqrt{n}},8\varepsilon_0 K^{-1}R^{-\tfrac{2}{3}}18C\frac{(\log n)^3}{\sqrt{n}}\Big\}\leq K^{-2}R^{-\tfrac{8}{3}}.
	\end{equation*}
	\item When $\varepsilon_0\leq K^{-1}R^{-2}/2160$, we have
	\begin{equation*}
	8\varepsilon_0 K^{-1}R^{-\tfrac{2}{3}}270\leq K^{-2}R^{-\tfrac{8}{3}}.
	\end{equation*}
	\item When the incoherence parameter satisfies $\Gamma\leq K^{-1/2}/216$, we have
	\begin{equation*}
	K^{-\tfrac{3}{2}}2160R^{-\tfrac{8}{3}}\Gamma\leq K^{-2}R^{-\tfrac{8}{3}}.
	\end{equation*}
\end{itemize} 
Note that those above conditions can be fulfilled by Conditions \ref{con:incoherence_sym}, \ref{con:sample} and \eqref{con:initial_input}. Thus,  we are able to simplify $A_1$ by
\begin{equation}\label{ineqn:bound_A1}
A_1\geq 32K^{-2}R^{-\tfrac{8}{3}}\Big(\sum_{k=1}^{K}\|\bz_k\|_2^2\Big),
\end{equation}
with probability at least $1-9K^2/n$.

~\\
\noindent  \textbf{Step Two: Upper bound for $A_2$.}
We observe the fact that
\begin{equation}\label{eqn:def:A2}
\begin{split}
A_2=&\sum_{k=1}^K\Big\|\frac{1}{\phi}\sqrt[3]{\eta_k^*}\nabla_k\tilde{\cL}(\bar{\bbeta}_k)\Big\|_2^2\\
=&\sup_{\bw\in\mathbb S^{Ks-1}}\Big|\Big\langle\sum_{k=1}^{K}\frac{\sqrt[3]{\eta_k^*}}{\phi}\nabla_k\tilde{\cL}(\bar{\bbeta}_k), \bw\Big\rangle\Big|^2,
\end{split}
\end{equation} 
where $\mathbb S$ is a unit sphere. It is equivalent to show for any $\bw\in\mathbb S^{Ks-1}$, $A'_2 =|\langle\sum_{k=1}^{K}\frac{\sqrt[3]{\eta_k^*}}{\phi}\nabla_k\tilde{\cL}(\bar{\bbeta}_k), \bw\rangle| $ is upper bounded. According to the definition of noiseless gradient \eqref{def:noiseless_gradient}, $A_2'$ is explicitly written as
\begin{equation*}
A_2' = \frac{6}{n}\sum_{i=1}^n\Big(\sum_{k'=1}^K(\bx_i^{\top}\bar{\bbeta}_{k'})^3-\sum_{k'=1}^{K}(\bx_i^{\top}\bar{\bbeta}_{k'}^*)^3\Big)\Big(\sum_{k=1}^K\frac{(\sqrt[3]{\eta_k^*})^2}{\phi}(\bx_i^{\top}\bar{\bbeta}_k)^2(\bx_i^{\top}\bw)\Big).
\end{equation*}
Following by \eqref{ineqn:weight} and \eqref{eqn:decomposition}, similar decomposition can be made for $A_2'$ as follows, where the only difference is that we replace one $\bx_i^{\top}\bz_k$ by $\bx_i^{\top}\bw$.
\begin{eqnarray*}
	A_2' &\leq&K^{-2}R^{\tfrac{2}{3}}\eta_{\min}^{*-\tfrac{4}{3}}\Big[\frac{6}{n}\sum_{i=1}^n\Big(\sum_{k'=1}^K 3(\bx_i^{\top}\bz_{k'})(\bx_i^{\top}\bar{\bbeta}_{k'})^2\sum_{k=1}^K(\bx_i^{\top}\bw)(\bx_i^{\top}\bar{\bbeta}_k^*)^2\Big)\Leftarrow \ A_{21}'\\
	&&+\frac{6}{n}\sum_{i=1}^n\Big(\sum_{k'=1}^K 3(\bx_i^{\top}\bz_{k'})(\bx_i^{\top}\bar{\bbeta}_{k'})^2\sum_{k=1}^K2(\bx_i^{\top}\bz_k)(\bx_i^{\top}\bw)(\bx_i^{\top}\bar{\bbeta}_k^*)\Big)\Leftarrow \ A_{22}'\\
	&&+\frac{6}{n}\sum_{i=1}^n\Big(\sum_{k'=1}^K 3(\bx_i^{\top}\bz_{k'})(\bx_i^{\top}\bar{\bbeta}_{k'})^2\sum_{k=1}^K(\bx_i^{\top}\bz_k)^2(\bx_i^{\top}\bw)\Big)\Leftarrow \ A_{23}'\\
	&&+\frac{6}{n}\sum_{i=1}^n\Big(\sum_{k'=1}^K 3(\bx_i^{\top}\bz_{k'})^2(\bx_i^{\top}\bar{\bbeta}_{k'})\sum_{k=1}^K(\bx_i^{\top}\bw)(\bx_i^{\top}\bar{\bbeta}_k^*)^2\Big)\Leftarrow \ A_{24}'\\
	&&+\frac{6}{n}\sum_{i=1}^n\Big(\sum_{k'=1}^K 3(\bx_i^{\top}\bz_{k'})^2(\bx_i^{\top}\bar{\bbeta}_{k'})\sum_{k=1}^K2(\bx_i^{\top}\bz_k)(\bx_i^{\top}\bw)(\bx_i^{\top}\bar{\bbeta}_k^*)\Big)\Leftarrow \ A_{25}'\\
	&&+\frac{6}{n}\sum_{i=1}^n\Big(\sum_{k'=1}^K 3(\bx_i^{\top}\bz_{k'})^2(\bx_i^{\top}\bar{\bbeta}_{k'})\sum_{k=1}^K(\bx_i^{\top}\bz_k)(\bx_i^{\top}\bw)(\bx_i^{\top}\bar{\bbeta}_k^*)\Big)\Leftarrow \ A_{26}'\\
	&&+\frac{6}{n}\sum_{i=1}^n\Big(\sum_{k'=1}^K 3(\bx_i^{\top}\bz_{k'})^3\sum_{k=1}^K(\bx_i^{\top}\bw)(\bx_i^{\top}\bar{\bbeta}_k^*)^2\Big)\Leftarrow \ A_{27}'\\
	&&+\frac{6}{n}\sum_{i=1}^n\Big(\sum_{k'=1}^K 3(\bx_i^{\top}\bz_{k'})^3\sum_{k=1}^K2(\bx_i^{\top}\bz_k)(\bx_i^{\top}\bw)(\bx_i^{\top}\bar{\bbeta}_k^*)\Big)\Leftarrow \ A_{28}'\\
	&&+\frac{6}{n}\sum_{i=1}^n\Big(\sum_{k'=1}^K 3(\bx_i^{\top}\bz_{k'})^3\sum_{k=1}^K(\bx_i^{\top}\bz_k)^2(\bx_i^{\top}\bw)\Big)\Big].\Leftarrow \ A_{29}'\\
\end{eqnarray*}

Let's bound $A_{21}'$ first. By using the same technique when calculating $\mathbb E(A_{11})$ in \eqref{eqn:expectation_A11}, we derive an upper bound for $\mathbb E(A_{21}')$,
\begin{eqnarray*}
	\mathbb E(A_{21}')&\leq& 36\eta_{\max}^{*\tfrac{4}{3}}\Big(\sum_{k=1}^{K}\|\bz_k\|_2+(K-1)\sum_{k=1}^{K}\Gamma \|\bz_k\|_2\Big)\\
	&+&180\eta_{\max}^{*\tfrac{4}{3}}\Big(\sum_{k=1}^{K}\|\bz_k\|_2+(K-1)\sum_{k=1}^{K}\Gamma\|\bz_k\|_2\Big)
	 +54\eta_{\max}^{*\tfrac{4}{3}}\Big(K\sum_{k=1}^{K}\|\bz_k\|_2\Big).
\end{eqnarray*}
Equipped with Lemma \ref{lemma:tensor_spectral} and the definition of tensor spectral norm \eqref{def:sparse_spectral}, it suffices to bound $A_{21}'$ by
\begin{eqnarray*}
	R^{\tfrac{2}{3}}\eta_{\min}^{*-\tfrac{4}{3}}K^{-\tfrac{1}{2}}A_{21}'\leq K^{-2}R^{2}\Big[216+54K+216K\Gamma+18CK\delta_{n,p,s}\Big]\Big(\sum_{k=1}^{K}\|\bz_k\|_2\Big)
\end{eqnarray*} 
with probability at least $1-10K^2/n^3$, where $\delta_{n,p,s}$ is defined in \eqref{eqn:initial_stat_rate}.

The upper bounds for $A_{22}'$ to $A_{29}'$ follow similar forms. Combining them together, we can derive an upper bound for $A_2'$ as follows
\begin{eqnarray*}
	A_2'&\leq& K^{-2}R^{2}\Big[216+270K+18CK\delta_{n,p,s}\Big]\Big(\sum_{k=1}^{K}\|\bz_k\|_2\Big)\\
	&\leq& K^{-2}R^{2}\Big[220+270K\Big]\Big(\sum_{k=1}^{K}\|\bz_k\|_2\Big),
\end{eqnarray*} 
with probability at least $1-90K^2/n^3$, where the second inequality utilizes Condition \ref{con:sample}. Therefore, the upper bound of $A_2$ is given as follows
\begin{equation}\label{ineqn:bound_A2}
A_2\leq K^{-1}R^{4}[220+270K]^2\Big(\sum_{k=1}^K\|\bz_k\|_2^2\Big),
\end{equation} 
with probability at least $1-90K^2/n^3$.

~\\
\noindent  \textbf{Step Three: Upper bound for $A_3$.}
By the definition of noisy gradient and noiseless gradient, $A_3$ is explicitly written as
\begin{eqnarray*}
	A_3 &=& \sum_{k=1}^K\Big\|\frac{(\sqrt[3]{\eta_k^*})^2}{\phi}\frac{6}{n}\sum_{i=1}^n\epsilon_i(\bx_i^{\top}\bar{\bbeta}_k)^2\bx_i\Big\|_2^2\\
	&\leq& K^{-4}R^{\tfrac{4}{3}}\eta_{\min}^{*-\tfrac{8}{3}}\sum_{k=1}^{K}\Big(\sqrt{Ks}\max_j\frac{6}{n}\sum_{i=1}^n\epsilon_i(\bx_i^{\top}\bar{\bbeta}_k)^2x_{ij}\Big)^2,
\end{eqnarray*} 
where the second inequality comes from \eqref{ineqn:weight}. For fixed $\{\epsilon_i\}_{i=1}^n$, applying Lemma \ref{cor:concentr_high_order}, we have
\begin{eqnarray*}
	\Big|\sum_{i=1}^n\epsilon_i(\bx_i^{\top}\bar{\bbeta}_k)^2x_{ij} - \mathbb E\Big(\sum_{i=1}^n\epsilon_i(\bx_i^{\top}\bar{\bbeta}_k)^2x_{ij}\Big)\Big|\leq C(\log n)^{\tfrac{3}{2}}\|\bepsilon\|_2\|\bar{\bbeta}_k\|_2^2,
\end{eqnarray*}
with probability at least $1-1/n$. Together with Lemma \ref{lemma:sub_exp}, we obtain for any $j\in[Ks]$, 
\begin{eqnarray*}
	\Big|\frac{6}{n}\sum_{i=1}^n \epsilon_i(\bx_i^{\top}\bar{\bbeta}_k)^2x_{ij}\Big|\leq 6CC_0\sigma \|\bar{\bbeta}_k\|_2^2\frac{(\log n)^{3/2}}{\sqrt{n}},
\end{eqnarray*}
with probability at least $1-4/n$, where $\sigma$ is the noise level. According to \eqref{con:initial_input}, 
\begin{eqnarray*}
	\Big\|\bar{\bbeta}_k-\bar{\bbeta}_k^*\Big\|_2^2\leq \sum_{k=1}^K\Big\|\bar{\bbeta}_k-\bar{\bbeta}_k^*\Big\|_2^2\leq K\eta_{\max}^{*\tfrac{2}{3}}\varepsilon_0^2,
\end{eqnarray*}
which further implies $\|\bar{\bbeta}_k\|_2^2\leq (1+K^{\tfrac{1}{2}}\varepsilon_0)^2\eta_{\max}^{*\tfrac{2}{3}}$. Equipped with union bound over $j\in[Ks]$,
\begin{eqnarray*}
	\max_{j\in[Ks]}\Big|\frac{6}{n}\sum_{i=1}^n \epsilon_i(\bx_i^{\top}\bar{\bbeta}_k)^2x_{ij}\Big|\leq 6CC_0\sigma (1+K^{\tfrac{1}{2}}\varepsilon_0)^2(\sqrt[3]{\eta_{\max}^*})^2\frac{(\log n)^{3/2}}{\sqrt{n}},
\end{eqnarray*}
with probability at least $1-4Ks/n$. Letting $C=6C_0(Ce)^{-2/3}(1+K^{\tfrac{1}{2}}\varepsilon_0)^2$,
\begin{eqnarray}\label{ineqn:bound_A3}
A_3\leq C\eta_{\min}^{*-\tfrac{4}{3}}R^{\tfrac{8}{3}}\sigma^2K^{-2}\frac{s(\log n)^3}{n},
\end{eqnarray}
with probability at least $1-4Ks/n$.

~\\
\noindent  \textbf{Step Four: Upper bound for $A_4$.} This cross term can be written as 
\begin{eqnarray*}
	A_4 = 2\sum_{k=1}^{K}\frac{\mu}{\phi}(\sqrt[3]{\eta_k^*})^2\Big(\frac{1}{n}\sum_{i=1}^n\epsilon_i(\bx_i^{\top}\bar{\bbeta}_k)^2(\bx_i^{\top}\bz_k)\Big).
\end{eqnarray*}
To bound this term, we take the same step in Step Three which fixes the noise term $\{\epsilon_i\}_{i=1}^n$ first. Similarly, we obtain with probability at least $1-4K/n$, 
\begin{eqnarray}\label{ineqn:bound_A4}
A_4\leq 2C\sigma \frac{(\log n)^{\tfrac{3}{2}}}{\sqrt{n}}K^{-1}R^{\tfrac{4}{3}}\eta_{\min}^{*-\tfrac{2}{3}}.
\end{eqnarray} 
This term is negligible in terms of the order when comparing with \eqref{ineqn:bound_A3}.

~\\
\noindent  \textbf{Summary.} Putting the bounds \eqref{ineqn:bound_A1}, \eqref{ineqn:bound_A2}, \eqref{ineqn:bound_A3} and \eqref{ineqn:bound_A4} together, we achieve an upper bound for gradient update effect as follows,
\begin{equation}\label{upper_bound:gradient}
\begin{split}
A\leq &\Big(1-64\mu K^{-2}R^{-\tfrac{8}{3}}+2\mu^2K^{-1}R^4[220+270K]^2\Big)\sum_{k=1}^K\|\bz_k\|_2^2\\
&+4\mu CK^{-2}\eta_{\min}^{*-\tfrac{4}{3}}R^{\tfrac{8}{3}}\frac{\sigma^2s(\log n)^3}{n},
\end{split}
\end{equation}
with probability at least $1-(18K^2+4K+4Ks)/n$.\hfill $\blacksquare$\\

\subsubsection{Bounding thresholding effect}\label{subsubsection_threshold}

The thresholding effect term in \eqref{eqn:error_decom} can also be decomposed into optimization error and statistical error. Recall that $B$ can be explicitly written as
\begin{eqnarray*}
	B = \sum_{k=1}^K\Big\|\mu \frac{\eta_k^{*\tfrac{2}{3}}}{\phi}\frac{4\sqrt{\log (np)}}{n}\sqrt{\sum_{i=1}^n\Big(\sum_{k'=1}^K(\bx_i^{\top}\bar{\bbeta}_{k'})^3-y_i\Big)^2(\bx_i^{\top}\bar{\bbeta}_k)^4}\gamma_k\Big\|_2^2,
\end{eqnarray*}
where $\supp(\gamma_k)\subset F_k$ and $\|\gamma_k\|_{\infty}\leq 1$. By using $(a+b)^2\leq 2(a^2+b^2)$, we have
\begin{equation*}
\begin{split}
	B\leq& \mu^2 \frac{64Ks\log p}{n}\Big[\underbrace{\frac{1}{n}\sum_{i=1}^n\Big(\sum_{k'=1}^K(\bx_i^{\top}\bar{\bbeta}_{k'})^3-\sum_{k'=1}^K(\bx_i^{\top}\bar{\bbeta}_{k'}^*)^3\Big)\Big(\sum_{k=1}^K\frac{\eta_k^{*\tfrac{4}{3}}}{\phi^2}(\bx_i^{\top}\bar{\bbeta}_k)^4\Big)}_{B_1:\text{optimization error}}\\
	&+\underbrace{\frac{1}{n}\sum_{i=1}^n\epsilon_i^2\sum_{k=1}^K\frac{\eta_k^{*\tfrac{4}{3}}}{\phi^2}(\bx_i^{\top}\bar{\bbeta}_k)^4}_{B_2:\text{statistical error}}\Big]. 
	\end{split}
\end{equation*}

~\\
\noindent  \textbf{Bounding $B_1$.} This optimization error term shares similar structure with \eqref{eqn:def:A2} but with higher order. Therefore, we follow the same idea as we did in bounding \eqref{eqn:def:A2}. Following by \eqref{ineqn:weight} and some basic expansions and inequalities, 
\begin{equation*}
\begin{split}
	B_1\leq& K^{-2}R^{\tfrac{4}{3}}\eta_{\min}^{*-\tfrac{8}{3}}\frac{1}{n}\Big(\sum_{k'=1}^K(\bx_i^{\top}\bar{\bbeta}_{k'})^3-\sum_{k'=1}^K(\bx_i^{\top}\bar{\bbeta}_{k'}^*)^3\Big)\Big(\sum_{k=1}^K(\bx_i^{\top}\bar{\bbeta}_k)^4\Big)\\
	\leq& K^{-2}R^{\tfrac{4}{3}}\eta_{\min}^{*-\tfrac{8}{3}}\Big[\frac{1}{n}\sum_{i=1}^n\Big(\sum_{k=1}^K3K(\bx_i^{\top}\bz_k)^6+9K(\bx_i^{\top}\bz_k)^4(\bx_i^{\top}\bar{\bbeta}_k^*)^2\\
	&+9K(\bx_i^{\top}\bz_k)^2(\bx_i^{\top}\bar{\bbeta}_k^*)^4\Big)\sum_{k'=1}^K(\bx_i^{\top}\bar{\bbeta}_{k'})^4\Big].
	\end{split}
\end{equation*}
The main term is $(\bx_i^{\top}\bz_k)^2(\bx_i^{\top}\bar{\bbeta}_k^*)^4$ according to the order of $\bar{\bbeta}_k^*$. We bound the main term first. Note that there exists some positive large constant $C$ such that 
\begin{eqnarray*}
	\mathbb E\Big(\frac{1}{n}\sum_{i=1}^n(\bx_i^{\top}\bz_k)^2(\bx_i^{\top}\bar{\bbeta}_{k'}^*)^4(\bx_i^{\top}\bar{\bbeta}_{k'})^4\Big)\leq C\|\bz_k\|_2^2\|\bar{\bbeta}_k^*\|_2^4\|\bar{\bbeta}_{k'}\|_2^4.
\end{eqnarray*}
Together with Lemma \ref{cor:concentr_high_order} and \eqref{con:initial_input}, we have
\begin{equation*}
\begin{split}
	&\sum_{k=1}^K\sum_{k'=1}^K\Big(\frac{1}{n}\sum_{i=1}^n(\bx_i^{\top}\bz_k)^2(\bx_i^{\top}\bar{\bbeta}_{k'}^*)^4(\bx_i^{\top}\bar{\bbeta}_{k'})^4\Big)\\
	\leq& C\Big(1+\frac{(\log n)^5}{\sqrt{n}}\Big)K^2\eta_{\max}^{*\tfrac{8}{3}}(1+K^{\tfrac{1}{2}}\varepsilon_0)^4\sum_{k=1}^K\|\bz_k\|_2^2.
	\end{split}
\end{equation*}
with probability at least $1-3K^2/n$. Overall, the upper bound of $B_1$ takes the form 
\begin{equation}\label{ineqn:bound_B1}
\begin{split}
B_1\leq& K^{-2}R^{\tfrac{4}{3}}\eta_{\min}^{*-\tfrac{8}{3}}\Big[18C\Big(1+\frac{(\log n)^5}{\sqrt{n}}\Big)K^2\eta_{\max}^{*\tfrac{8}{3}}(1+K^{\tfrac{1}{2}}\varepsilon_0)^4\sum_{k=1}^K\|\bz_k\|_2^2\Big]\\
\leq& R^{4}18C\Big(1+\frac{(\log n)^5}{\sqrt{n}}\Big)(1+K^{\tfrac{1}{2}}\varepsilon_0)^4\sum_{k=1}^K\|\bz_k\|_2^2,
\end{split}
\end{equation}
with probability at least $1-3K^2/n$.

~\\
\noindent  \textbf{Bounding $B_2$.} We rewrite $B_2$ by 
\begin{eqnarray*}
	B_2 = \sum_{k=1}^K\frac{\eta_k^{*\tfrac{4}{3}}}{\phi^2}\Big(\frac{1}{n}\sum_{i=1}^n\epsilon_i^2(\bx_i^{\top}\bar{\bbeta}_k)^4\Big).
\end{eqnarray*}
For fixed $\{\epsilon_i\}_{i=1}^n$, accordingly to Lemma \ref{cor:concentr_high_order}, we have
\begin{eqnarray*}
	\Big|\sum_{i=1}^n\epsilon_i^2(\bx_i^{\top}\bar{\bbeta}_k)^4-\mathbb E\Big(\sum_{i=1}^n\epsilon_i^2(\bx_i^{\top}\bar{\bbeta}_k)^4\Big)\Big|\leq C(\log n)^2\|\bepsilon^2\|_2\|\bar{\bbeta}_k\|_2^4.
\end{eqnarray*}
Note that $\mathbb E((\bx_i^{\top}\bar{\bbeta}_k)^4) = 3\|\bar{\bbeta}_k\|_2^4$. It will reduce to
\begin{eqnarray*}
	\frac{1}{n}\sum_{i=1}^n\epsilon_i^2(\bx_i^{\top}\bar{\bbeta}_k)^4\leq \Big(\frac{3}{n}\sum_{i=1}^n\epsilon_i^2+C\frac{(\log n)^2}{n}\|\bepsilon^2\|_2\Big)\|\bar{\bbeta}_k\|_2^4.
\end{eqnarray*}
From Lemma \ref{lemma:sub_exp}, with probability at least $1-3/n$, 
\begin{eqnarray*}
	|\frac{1}{n}\sum_{i=1}^n\epsilon_i^2|\leq C_0\sigma^2, \ \frac{1}{n}\|\bepsilon^2\|_2\leq C_0 \frac{\sigma^2}{\sqrt{n}}.
\end{eqnarray*}
Combining the above two inequalities, we obtain
\begin{eqnarray}
\Big|\frac{1}{n}\sum_{i=1}^n\epsilon_i^2(\bx_i^{\top}\bar{\bbeta}_k)^4\Big|\leq 6C_0\sigma^2\|\bar{\bbeta}_k\|_2^4,
\end{eqnarray}
with probability at least $1-7/n$.
Plugging in the definition of $\phi$ and \eqref{con:initial_input}, $B_2$ is upper bounded by
\begin{equation}\label{ineqn:bound_B2}
B_2\leq 6C_0\sigma^2(1+K^{\tfrac{1}{2}}\varepsilon_0)^4\eta_{\min}^{*-\tfrac{4}{3}}R^{\tfrac{8}{3}}K^{-3},
\end{equation}
with probability at least $1-7K/n$.  	

~\\
\noindent  \textbf{Summary.} Putting the bounds \eqref{ineqn:bound_B1} and \eqref{ineqn:bound_B2} together, we have similar upper bound for thresholded effect,
\begin{equation}\label{upper_bound:threhsholding}
B\leq C_2\mu^2 R^{4}\sum_{k=1}^K\|\bz_k\|_2^2 + C_3\mu^2\eta_{\min}^{*-\tfrac{4}{3}}R^{\tfrac{8}{3}}K^{-2}\frac{\sigma^2s\log p}{n},
\end{equation}
with probability at least $1-(3K^2+7K)/n$. \hfill $\blacksquare$\\

\subsubsection{Ensemble}
From the definition of $\gamma_k$, it's not hard to see actually the cross term $C$ is equal to zero. Combining the upper bound of gradient update effect \eqref{upper_bound:gradient} and thresholding effect \eqref{upper_bound:threhsholding} together, we obtain 
\begin{equation*}
\begin{split}
&\sum_{k=1}^K\Big\|\sqrt[3]{\eta_k}\tilde{\bbeta}_k^+-\sqrt[3]{\eta_k^*}\bbeta_k^*\Big\|_2^2\\
\leq&\Big(1-64\mu K^{-2}R^{-\tfrac{8}{3}}+3\mu^2K^{-1}R^4[220+270K]^2\Big)\Big(\sum_{k=1}^K\|\bz_k\|_2^2\Big)\\
&+2C_3\mu^2R^{\tfrac{8}{3}}\eta_{\min}^{*-\tfrac{4}{3}}\frac{\sigma^2K^{-2}s\log p}{n}.
\end{split}
\end{equation*}
As long as the step size $\mu$ satisfies
\begin{equation*}
0<\mu\leq \frac{32 R^{-20/3}}{3K[220+270K]^2},
\end{equation*}
we reach the conclusion
\begin{equation}
\begin{split}
&\sum_{k=1}^K\Big\|\sqrt[3]{\eta_k}\tilde{\bbeta}_k^+-\sqrt[3]{\eta_k^*}\bbeta_k^*\Big\|_2^2\\
\leq &\Big(1-32\mu K^{-2}R^{-\tfrac{8}{3}}\Big)\sum_{k=1}^K\Big\|\sqrt[3]{\eta_k}\bbeta_k-\sqrt[3]{\eta_k^*}\bbeta_k^*\Big\|_2^2\\
&+2C_3\mu^2R^{-\tfrac{8}{3}}\eta_{\min}^{*-\tfrac{4}{3}}\frac{\sigma^2K^{-2}s\log p}{n},
\end{split}
\end{equation}
with probability at least $1-4Ks/n$.\hfill $\blacksquare$

\subsection{Proof of Lemma \ref{lemma:oracle_equivalence_sym}}

Let us consider $k$-th component first. Without loss of generality, suppose $F\subset\{1,2,\ldots,Ks\}$. For $j = Ks+1, \ldots,p$,
\begin{equation}
\frac{\partial}{\partial \beta_{kj}}\cL(\bbeta_k) = \frac{2}{n}\sum_{i=1}^n\Big(\sum_{k=1}^K\eta_k(\bx_i^{\top}\bbeta_k)^3-y_i\Big)\eta_k(\bx_i^{\top}\bbeta_k)^2x_{ij},
\end{equation} 
and it's not hard to see the independence between $\{\bx_i^{\top}\bbeta_k, y_i\}$ and $x_{ij}$. Applying standard Hoeffding's inequality, we have with probability at least $1-\tfrac{1}{n^2p^2}$,
\begin{eqnarray*}
	\Big|\frac{\partial}{\partial \beta_{kj}}\cL(\bbeta_k)\Big|\leq \frac{\sqrt{4\log(np)}}{n}\sqrt{\sum_{i=1}^n(\sum_{k=1}^K\eta_k(\bx_i^{\top}\bbeta_k)^3-y_i)^2(\eta_k(\bx_i^{\top}\bbeta_k))^2}=h(\bbeta_k).
\end{eqnarray*}
Equipped with union bound, with probability at least $1-\tfrac{1}{n^2p}$,
\begin{equation*}
\max_{Ks+1\leq j\leq p} \Big|\frac{\partial}{\partial \beta_{kj}}\cL(\bbeta_k)\Big|\leq h(\bbeta_k).
\end{equation*}
Therefore, according to the definition of thresholding function $\varphi(\bx)$, we obtain the following equivalence, 
\begin{equation}\label{eqn:equivalence}
\varphi_{\frac{\mu}{\phi}h(\bbeta_k)}\Big(\bbeta_k-\frac{\mu}{\phi}\nabla_{\bbeta_k}\cL(\bbeta_k)\Big) = \varphi_{\frac{\mu}{\phi}h(\bbeta_k)}\Big(\bbeta_k-\frac{\mu}{\phi}\nabla_{\bbeta_k}\cL(\bbeta_k)_F\Big),
\end{equation}
holds for $k\in[K]$, with probability  at least $1-\tfrac{1}{n^2p}$. \eqref{eqn:equivalence} also provides that $\supp(\bbeta_k^+)\subset F$ for every $k\in[K]$, which further implies $F^+\subset F$. Now we end the proof. \hfill $\blacksquare$

\subsection{Proof of Lemma \ref{lemma:marginal_effect}}

First, we consider symmetric case. According to the definition of $\{y_i\}_{i=1}^n$ from symmetric tensor estimation model \eqref{eq:model}, we separate the random noise  $\epsilon_i$ by the following expansion,
\begin{eqnarray}
\label{eqn:y_n}
\frac{1}{n}\sum_{i=1}^n y_i^2 &=& \frac{1}{n}\sum_{i=1}^{n}\Big[\sum_{k=1}^{K}\eta_k^*(\bx_i^{\top}\bbeta_k^*)^3+\epsilon_i\Big]^2\nonumber\\
&=&\underbrace{\frac{1}{n}\sum_{i=1}^{n}(\sum_{k=1}^{K}\eta_k^*(\bx_i^{\top}\bbeta_k^*)^3)^2}_{I_1}+\underbrace{\frac{2}{n}\sum_{i=1}^{n}\epsilon_i\sum_{k=1}^{K}\eta_k^*(\bx_i^{\top}\bbeta_k^*)^3}_{I_2}+\underbrace{\frac{1}{n}\sum_{i=1}^{n}\epsilon_i^2}_{I_3}.
\end{eqnarray}

~\\
\noindent\textbf{Bounding $I_1$.} We expand $i$-th component of $I_1$ as follows
\begin{equation}
\label{eqn:marginal_1}
\begin{split}
&(\sum_{k=1}^{K}\eta_k^*(\bx_i^{\top}\bbeta_k^*)^3)^2 \\
=& \sum_{k=1}^{K}\eta_k^*(\bx_i^{\top}\bbeta_k^*)^6 + 2 \sum_{k_i< k_j}\eta_{k_i}^*\eta_{k_j}^*(\bx_i^{\top}\bbeta_{k_i}^*)^3(\bx_i^{\top}\bbeta_{k_j}^*)^3.
\end{split}
\end{equation}
As shown in Corollary \ref{cor:expectation_gaussian_vector}, the expectations of above two parts takes forms of 
\begin{eqnarray*}
	 && \mathbb E(\bx_i^{\top}\bbeta_{k_i}^*)^3(\bx_i^{\top}\bbeta_{k_j}^*)^3 = 6(\bbeta_{k_i}^{*\top}\bbeta_{k_j}^*)^3 + 9(\bbeta_{k_i}^{*\top}\bbeta_{k_j}^*)\|\bbeta_{k_i}^*\|_2^2\|\bbeta_{k_j}^*\|_2^2\\
	 &&\mathbb E(\bx_i^{\top}\bbeta_k^*)^6 = 15\|\bbeta_k^*\|_2^2.
\end{eqnarray*} 
Recall that $\|\bbeta_k^*\|_2=1$ for any $k\in[K]$ and Condition \ref{con:incoherence_sym} implies for any $k_i\neq k_j$,
$|\bbeta_{k_i}^{*\top}\bbeta_{k_j}^*|\leq \Gamma,$
where $\Gamma$ is the incoherence parameter. Thus, $\mathbb E(\bx_i^{\top}\bbeta_{k_i}^*)^3(\bx_i^{\top}\bbeta_{k_j}^*)^3$ is upper bounded by
\begin{eqnarray}\label{eqn:marginal_2}
\Big|\mathbb E(\bx_i^{\top}\bbeta_{k_i}^*)^3(\bx_i^{\top}\bbeta_{k_j}^*)^3\Big| \leq 6\Gamma^3 +9\Gamma, \ \text{for any} \ k_i\neq k_j.
\end{eqnarray}
By using the concentration result in Lemma \ref{cor:concentr_high_order}, we have with probability at least $1-1/n$
\begin{equation}
\label{eqn:marginal_3}
\begin{split}
&\Big|\frac{1}{n}\sum_{i=1}^n(\bx_i^{\top}\bbeta_k^*)^6-\mathbb E(\frac{1}{n}\sum_{i=1}^n(\bx_i^{\top}\bbeta_k^*)^6)\Big|\leq C_1\frac{(\log n)^3}{\sqrt{n}},\\
&\Big|\frac{1}{n}\sum_{i=1}^n(\bx_i^{\top}\bbeta_{k_i}^*)^3(\bx_i^{\top}\bbeta_{k_j}^*)^3 -\mathbb E(\frac{1}{n}\sum_{i=1}^n(\bx_i^{\top}\bbeta_{k_i}^*)^3(\bx_i^{\top}\bbeta_{k_j}^*)^3)\Big|\leq C_1\frac{(\log n)^3}{\sqrt{n}}.
\end{split}
\end{equation}
Putting \eqref{eqn:marginal_1},\eqref{eqn:marginal_2} and \eqref{eqn:marginal_3} together, this essentially provides an upper bound for $I_1$, namely
\begin{equation}
\label{eqn:marginal_I1}
\frac{1}{n}\sum_{i=1}^{n}(\sum_{k=1}^{K}\eta_k^*(\bx_i^{\top}\bbeta_k^*)^3)^2\leq \Big(15+6\Gamma^3+9\Gamma+2 C_1\frac{(\log n)^3}{\sqrt{n}}\Big) (\sum_{k=1}^{K}\eta_k^*)^2,
\end{equation}
with probability at least $1-K^2/n$.

~\\
\noindent\textbf{Bounding $I_2$.} Since the random noise $\{\epsilon_i\}_{i=1}^n$ is of mean zero and independent of $\{\bx_i\}$, we have
$$\mathbb E(\epsilon_i\sum_{k=1}^{K}\eta_k^*(\bx_i^{\top}\bbeta_k^*)^3) = 0.$$
By using the independence and Corollary \ref{cor:concentr_high_order}, we have
\begin{eqnarray*}
	&&\mathbb P\Big(\frac{1}{n}\sum_{i=1}^{n}\epsilon_i(\bx_i^{\top}\bbeta_k^*)^3\geq C_2\frac{(\log n)^{\tfrac{3}{2}}}{n}\sqrt{n} \sigma\Big)\\
	&&\leq \mathbb P\Big(\frac{1}{n}\sum_{i=1}^{n}\epsilon_i(\bx_i^{\top}\bbeta_k^*)^3\geq C_2\frac{(\log n)^{\tfrac{3}{2}}}{n}\sqrt{n} \sigma\Big|\|\bepsilon\|_2\leq C_0\sigma \sqrt{n}\Big) +\mathbb P\Big(\|\bepsilon\|_2\geq C_0\sqrt{n}\sigma\Big)\\
	&&\leq \frac{1}{n}+\frac{3}{n} = \frac{4}{n}.
\end{eqnarray*}
This further implies that 
\begin{eqnarray}
\label{eqn:marginal_I2}
\frac{1}{n}\sum_{i=1}^n\sum_{k=1}^{K}\eta_k^*(\bx_i^{\top}\bbeta_k^*)^3\epsilon_i\leq (\sum_{k=1}^K\eta_k^*)C_2\frac{(\log n)^{\tfrac{3}{2}}}{\sqrt{n}}\sigma,
\end{eqnarray}
with probability at least $1-4K/n$. 

~\\
\noindent\textbf{Bounding $I_3$.} As shown in Lemma \ref{lemma:sub_exp}, the random noise $\epsilon_i$ with sub-exponential tail satisfies 
\begin{eqnarray}
\label{eqn:marginal_I3}
\frac{1}{n}\sum_{i=1}^n\epsilon_i^2 \leq C_3 \sigma^2 .
\end{eqnarray}
with probability at least $1-3/n$.

Overall, putting \eqref{eqn:marginal_I1}, \eqref{eqn:marginal_I2} and \eqref{eqn:marginal_I3} together, we have with probability at least $1-(K^2+4K+3)/n$,
\begin{equation*}
\frac{\frac{1}{n}\sum_{i=1}^{n}y_i^2}{(\sum_{k=1}^K\eta_k^*)^2}
\leq 15 + 6\Gamma^3+9\Gamma+2 C_1\frac{(\log n)^3}{\sqrt{n}} + \frac{2 C_2\sigma}{(\sum_{k=1}^K\eta_k^*)}\frac{(\log n)^{\tfrac{3}{2}}}{\sqrt{n}} + \frac{C_3 \sigma^2 }{(\sum_{k=1}^K\eta_k^*)^2}.
\end{equation*}
Under Conditions \ref{con:noise} \& \ref{con:sample}, the above bound reduces to
\begin{equation*}
\frac{1}{n}\sum_{i=1}^{n}y_i^2\leq (16+6\Gamma^3+9\Gamma)(\sum_{k=1}^K\eta_k^*)^2, 
\end{equation*}
with probability at least $1-(K^2+4K+3)/n$. The proof of lower bound is similar, and hence is omitted here.

Similar results will also hold for non-symmetric tensor estimation model. Throughout the proof, the only difference is that $$\mathbb E(\bu_i^{\top}\bbeta_{1k}^*)^2(\bv_i^{\top}\bbeta_{2k}^*)^2(\bw_i^{\top}\bbeta_{3k}^*)^2=1.$$ \hfill $\blacksquare$\\


\section{Non-symmetric Tensor Estimation}\label{sec:non_symmetric}

\subsection{Conditions and Algorithm}\label{con_alg_nonsym}
In this subsection, we provide several essential conditions for Theorem \ref{thm:main_nym} and the detail algorithm for non-symmetric tensor estimation.
\begin{Condition}[(Uniqueness of CP-decomposition)]\label{con:ident_nym}
	The CP-decomposition form \eqref{eq:low_rank_nonsym} is unique in the sense that if there exists another CP-decomposition $\mT^* = \sum_{k=1}^{K'} \eta_k^{*'}\bbeta_{1k}^{*'}\circ\bbeta_{2k}^{*'}\circ\bbeta_{3k}^{*'}$, it must have $K=K'$ and be invariant up to a permutation of $\{1,\ldots, K\}$.
\end{Condition}
\begin{Condition}[(Parameter space)]\label{con:parameter_nym}
	The CP-decomposition of $\mT^* = \sum_{k=1}^K \eta_k^*\bbeta_{1k}^*\circ\bbeta_{2k}^*\circ\bbeta_{3k}^*$ satisfies
	$$\|\mT^*\|_{op}\leq C_1\eta_{\max}^*, \quad K=\cO(s), \quad \text{and}\quad R =\eta_{\max}^*/\eta_{\min}^*\leq C_2
	$$
	for some absolute constants $C_1,C_2$.
\end{Condition}
\begin{Condition}[(Parameter incoherence)]
	\label{con:incoherence_nym}
	The true tensor components are incoherent such that
	\begin{eqnarray*}
		\Gamma:=\max_{k_i\neq k_j}\Big\{|\langle\bbeta_{1k_i}^*, \bbeta_{1k_j}^*\rangle|, |\langle\bbeta_{2k_i}^*, \bbeta_{2k_j}^*\rangle|,|\langle\bbeta_{3k_i}^*, \bbeta_{3k_j}^*\rangle|\Big\}\leq C\min\{K^{-\tfrac{3}{4}}R^{-1},s^{-\tfrac{1}{2}}\}.
	\end{eqnarray*}
\end{Condition}
\begin{Condition}[(Random noise)]\label{con:noise_non}
	We assume the random noise $\{\epsilon_i\}_{i=1}^n$ follows a sub-exponential tail with parameter $\sigma$ satisfying $0<\sigma<C\sum_{k=1}^K\eta_k^*$. 
\end{Condition}

\begin{algorithm}\label{alg:overall_nonsymmetric}
	\small
	\caption{Non-symmetric tensor estimation via cubic sketchings}
	\begin{algorithmic}[1]
		\REQUIRE response $\{y_i\}_{i=1}^n$, sketching vector $\{\bu_i,\bv_i,\bw_i\}_{i=1}^n$, truncation level $d$, step size $\mu$, rank $K$, stopping error $\epsilon = 10^{-4}$.
		\STATE  \textbf{Step 1:} Calculate the moment-based tensor $\cT$ as \eqref{eq:moment_asy} and do sparse tensor decomposition on $\cT$ to get a warm-start $\{\boldsymbol{\eta}^{(0)}, \bB_1^{(0)}, \bB_2^{(0)}, \bB_3^{(0)}\}$.
		\STATE \textbf{Step 2:} Let $t=0$. \\
		\STATE \ \ \textbf{Repeat} block-wise thresholded gradient update
		\STATE \begin{itemize}
			\item Compute threshold level $\bh(\bB_k)$ as defined in Step Two.
			\item Calculated block-wise thresholded gradient descent update
			\begin{equation*}
				\begin{split}
					&\tv(\bB_1^{(t+1)}) = \varphi_{\frac{\mu\bh(\bB_1)}{\phi}}\Big(\tv(\bB_1^{(t)}) -\frac{\mu}{\phi}\nabla_{\bB_1}\cL(\bB_1^{(t)}, \bB_2^{(t)}, \bB_3^{(t)})\Big)\\
					&\tv(\bB_2^{(t+1)}) = \varphi_{\frac{\mu\bh(\bB_2)}{\phi}}\Big(\tv(\bB_2^{(t)}) -\frac{\mu}{\phi}\nabla_{\bB_2}\cL(\bB_1^{(t)}, \bB_2^{(t)}, \bB_3^{(t)})\Big)\\
					&\tv(\bB_3^{(t+1)}) = \varphi_{\frac{\mu\bh(\bB_3)}{\phi}}\Big(\tv(\bB_3^{(t)}) -\frac{\mu}{\phi}\nabla_{\bB_3}\cL(\bB_1^{(t)}, \bB_2^{(t)}, \bB_3^{(t)})\Big),
				\end{split}
			\end{equation*} 
			where $\phi=\frac{1}{n}\sum_{i=1}^ny_i^2$. The detail form of $\nabla_{\bB_1}\cL, \nabla_{\bB_2}\cL, \nabla_{\bB_3}\cL$  can refer \eqref{eqn:gradient_nonsym}
		\end{itemize}
		\STATE \ \  \textbf{Until} $\max\{\|\bB_j^{(T+1)} - \bB_j^{(T)}\|_F\}\leq \epsilon$. \\
		
		\STATE \textbf{Step 3:} Do column-wise normalization as 
		\begin{eqnarray*}
			&&\hat{\bB}_1 = \Big(\frac{\bbeta_{11}^{(T)}}{\|\bbeta_{11}^{(T)}\|_2}, \ldots, \frac{\bbeta_{1K}^{(T)}}{\|\bbeta_{1K}^{(T)}\|_2}\Big), \\
			&&\hat{\bB}_2 = \Big(\frac{\bbeta_{21}^{(T)}}{\|\bbeta_{21}^{(T)}\|_2}, \ldots, \frac{\bbeta_{2K}^{(T)}}{\|\bbeta_{2K}^{(T)}\|_2}\Big), \\
			&&\hat{\bB}_3 = \Big(\frac{\bbeta_{31}^{(T)}}{\|\bbeta_{31}^{(T)}\|_2}, \ldots, \frac{\bbeta_{3K}^{(T)}}{\|\bbeta_{3K}^{(T)}\|_2}\Big).
		\end{eqnarray*} 
		And update the weight by 
		$$\hat{\boldsymbol{\eta}} = \boldsymbol{\eta}^{(0)} * (\|\bbeta_{11}^{(T)}\|_2\|\bbeta_{21}^{(T)}\|_2\|\bbeta_{31}^{(T)}\|_2, \ldots, \|\bbeta_{1K}^{(T)}\|_2\|\bbeta_{3K}^{(T)}\|_2\|\bbeta_{3K}^{(T)}\|_2)^{\top}.
		$$ The final estimator is $\hat{\mT} = \sum_{k=1}^{K}\hat{\eta}_{k}\hat{\bbeta}_{1k}\circ \hat{\bbeta}_{2k}\circ \hat{\bbeta}_{3k}$.\\
		\RETURN non-symmetric tensor estimator $\hat{\mT}$.
	\end{algorithmic}
\end{algorithm}

\subsection{Proof of Theorem \ref{thm:main_nym}}\label{sec:proof_asym}

The main distinguished part of the proof for non-symmetric update is Lemma \ref{lemma:onestep_update_nym}: one-step oracle estimator, which is parallel to Lemma \ref{lemma:onestep_update_sym}. For the sake of completeness, we limit our attention to rank-one case and only provide the theoretical development for one-step oracle estimator in this subsection. The generalization to general rank case follows the exact same idea in the proof of symmetric update by incorporating the incoherence condition \eqref{con:incoherence_nym}. 

For rank-one non-symmetric tensor estimation, the model \eqref{eq:model_asymmetric} reduces to
\begin{equation*}
	y_i = \langle\eta^*\bbeta_1^*\circ \bbeta_2^*\circ \bbeta_3^*, \bu_i\circ \bv_i \circ \bw_i\rangle+\epsilon_i, \ \text{for} \ i=1,\ldots,n.
\end{equation*}
Suppose $|\supp(\bbeta_1^*)|=s_1$, $|\supp(\bbeta_2^*)|=s_2$, $|\supp(\bbeta_3^*)|=s_3$ and denote $s=\max\{s_1,s_2,s_3\}$. Define $F_j^{(t)} = \supp(\bbeta_j^*)\cup\supp(\bbeta_j^{(t)})$, $F^{(t)}=\cup_{j=1}^3F_j^{(t)}$ and the oracle estimator as
\begin{equation*}
	\tilde{\bbeta}_1^{(t+1)} = \varphi_{\frac{\mu}{\phi}h(\bbeta_1^{(t)})}\Big(\bbeta_j^{(t)}-\frac{\mu}{\phi}\nabla_1\cL(\bbeta_1^{(t)}, \bbeta_2^{(t)}, \bbeta_3^{(t)})_{F^{(t)}}\Big),
\end{equation*}
where $h(\bbeta_1^{(t)})$ has the form of 
\begin{equation}\label{def:threhsold_level_nym}
\frac{\sqrt{4\log np}}{n}\sqrt{\sum_{i=1}^n\Big(\eta(\bu_i^{\top}\bbeta_1^{(t)})(\bv_i^{\top}\bbeta_2^{(t)})(\bw_i^{\top}\bbeta_3^{(t)})-y_i\Big)^2\eta^{\tfrac{2}{3}}(\bv_i^{\top}\bbeta_2^{(t)})^2(\bw_i^{\top}\bbeta_3^{(t)})^2}.
\end{equation}
The definitions of $\tilde{\bbeta}_2^{(t+1)}$ and $\tilde{\bbeta}_3^{(t+1)}$ are similar.
\begin{Lemma}
	\label{lemma:onestep_update_nym}
	Let $t\geq 0$ be an integer. Suppose Conditions \ref{con:ident_nym}-\ref{con:noise_non} hold and $\{\bbeta_j^{(t)}, \eta\}$ satisfies the following upper bound
	\begin{equation}\label{con:initial_nym}
		\max_{j=1,2,3}\|\sqrt[3]{\eta}\bbeta_j^{(t)}-\sqrt[3]{\eta^*}\bbeta_j^*\|_2\leq \sqrt[3]{\eta^*}\varepsilon_0, \ |\eta-\eta^*|\leq \varepsilon_0
	\end{equation}
	with probability at least $1-CO(1/n)$. Assume the step size $\mu$ satisfies $0<\mu<\mu_0$ for some small absolute constant $\mu_0$ and $s\leq d\leq Cs$. Then $\{\tilde{\bbeta}_j^{(t+1)}\}$ can be upper bounded as
	\begin{eqnarray*}
		&&\max_{j=1,2,3}\Big\|\sqrt[3]{\eta}\tilde{\bbeta}_j^{(t+1)}-\sqrt[3]{\eta^*}\bbeta_j^*\Big\|_2\\
		&\leq& (1-\frac{\mu}{12})\max_{j=1,2,3}\Big\|\sqrt[3]{\eta}\bbeta_j^{(t)}-\sqrt[3]{\eta^*}\bbeta_j^*\Big\|_2+\mu\frac{3\sigma}{(\sqrt[3]{\eta^*})^2}\sqrt{\frac{3s \log p}{n}},
	\end{eqnarray*}
	with probability at least $1-12s/n$.
\end{Lemma}

\emph{Proof.} We focus on $j=1$ first. To simplify the notation, we drop the superscript of iteration index $t$, and denote iteration index $t+1$ by $+$. Moreover, denote $\bar{\bbeta}_j = \sqrt[3]{\eta}\bbeta_j$, $\bar{\bbeta}_j^+=\sqrt[3]{\eta}\bbeta_j$, $\bar{\bbeta}_j^{*}=\sqrt[3]{\eta^*}\bbeta_j^*$ for $j=1,2,3$. Then, the gradient function is rewritten as
\begin{eqnarray*}
	\nabla_1\cL(\bar{\bbeta}_1,\bar{\bbeta}_2,\bar{\bbeta}_3) =\sqrt[3]{\eta} \frac{2}{n}\sum_{i=1}^n\Big((\bu_i^{\top}\bar{\bbeta}_1)(\bv_i^{\top}\bar{\bbeta}_2)(\bw_i^{\top}\bar{\bbeta}_3)\Big)(\bv_i^{\top}\bar{\bbeta}_2)(\bw_i^{\top}\bar{\bbeta}_3)\bu_i.
\end{eqnarray*}

According to the definition of thresholded function, $\tilde{\bbeta}_1^+$ can be explicitly written by
\begin{eqnarray*}
	\tilde{\bbeta}_1^+ &=& \varphi_{\frac{\mu}{\phi} h(\bar{\bbeta}_1)}\Big(\bbeta_1-\frac{\mu}{\phi}\nabla_{1}\cL(\bar{\bbeta}_1, \bar{\bbeta}_2, \bar{\bbeta}_3)_{F}\Big)\\
	&=& \bbeta_1-\frac{\mu}{\phi}\nabla_{1}\cL(\bar{\bbeta}_1, \bar{\bbeta}_2, \bar{\bbeta}_3)_{F}+\frac{\mu}{\phi} h(\bar{\bbeta}_1)\bgamma,
\end{eqnarray*}
where $\bgamma\in \mathbb R^p, \supp(\bgamma)\subset F$ and $\|\bgamma\|_{\infty}\leq 1$. Then the oracle estimation error $\|\sqrt[3]{\eta}\tilde{\bbeta}_1^+-\sqrt[3]{\eta^*}\bbeta_1^*\|_2$ can be decomposed by the gradient update effect and the thresholded effect,
\begin{equation}
	\label{eqn:main_decomposition}
	\begin{split}
	&\Big\|\sqrt[3]{\eta}\tilde{\bbeta}_1^+-\sqrt[3]{\eta^*}\bbeta_1^*\Big\|_2\\
	 = &\underbrace{\Big\|\bar{\bbeta}_1-\bar{\bbeta}_1^*-\mu\frac{\sqrt[3]{\eta}}{\phi}\nabla_{1}\cL(\bar{\bbeta}_1, \bar{\bbeta}_2, \bar{\bbeta}_3)_{F}\Big\|_2}_{\text{gradient update effect}} + \underbrace{\mu\frac{\sqrt[3]{\eta}}{\phi}\Big|h(\bar{\bbeta}_1)\Big|\sqrt{3s}}_{\text{thresholded effect}}.
	 \end{split}
\end{equation}
By using the tri-convex structure of $\cL(\bar{\bbeta}_1, \bar{\bbeta}_2, \bar{\bbeta}_3)$, we borrow the analysis tool for vanilla gradient descent \cite{Bubeck15} given sufficient good initial. Following this proof strategy, we decompose the gradient update effect in \eqref{eqn:main_decomposition} by three parts,
\begin{eqnarray*}
	\Big\|\sqrt[3]{\eta}\tilde{\bbeta}_1^+-\sqrt[3]{\eta^*}\bbeta_1^*\Big\|_2&\leq&\underbrace{\Big\|\bar{\bbeta}_1-\bar{\bbeta}_1^*-\mu \frac{\sqrt[3]{\eta}}{\phi}\nabla_{1} \tilde{\cL}(\bar{\bbeta}_1,\bar{\bbeta}_2^*,\bar{\bbeta}_3^*)_{F}\Big\|_2}_{I_1}\\
	&+&\mu \underbrace{\frac{\sqrt[3]{\eta}}{\phi}\Big\|\nabla_{1} \tilde{\cL}(\bar{\bbeta}_1,\bar{\bbeta}_2^*,\bar{\bbeta}_3^*)_{F}-\nabla_{1}\tilde{\cL}(\bar{\bbeta}_1,\bar{\bbeta}_2,\bar{\bbeta}_3)_{F}\Big\|_2}_{I_2}\\
	&+&\mu \underbrace{\frac{\sqrt[3]{\eta}}{\phi}\Big\|\nabla_{1}\tilde{\cL}(\bar{\bbeta}_1,\bar{\bbeta}_2,\bar{\bbeta}_3)_{F}-\nabla_{1}\cL(\bar{\bbeta}_1,\bar{\bbeta}_2,\bar{\bbeta}_3)_{F}\Big\|_2}_{I_3}\\
	&+&\underbrace{\mu\frac{\sqrt[3]{\eta}}{\phi}\Big|h(\bar{\bbeta}_1)\Big|\sqrt{3s}}_{I_4},
\end{eqnarray*}
where $\nabla_1\tilde{\cL}$ is the noiseless gradient as we defined in \eqref{def:noiseless_gradient}.
We will bound $I_1,I_2, I_3, I_4$ successively in the following four subsections. For simplicity, during the following proof, we drop the index subscript $F$ as we did in Section \ref{proof:one-step_sym}. And $\phi = \sum_{i=1}^ny_i^2$ approximates $\eta^{*2}$ up to constant due to Lemma \ref{lemma:marginal_effect}.

\subsubsection{Bounding $I_1$}\label{proof_gradient_effect}

In this section, let us denote 
$$
\sqrt[3]{\eta}\tilde{\cL}(\bar{\bbeta}_1,\bar{\bbeta}_2^*, \bar{\bbeta}_3^*)/\phi = f(\bar{\bbeta}_1), \ \sqrt[3]{\eta}\nabla_{1}\tilde{\cL}(\bar{\bbeta}_1,\bar{\bbeta}_2^*,\bar{\bbeta}_3^*)/\phi=\nabla f(\bar{\bbeta}_1),
$$
 where $\supp(\nabla f(\bar{\bbeta}_1))=F$. When $\bbeta_2$ and $\bbeta_3$ are fixed, the update can be treated as a vanilla gradient descent update. The following proof follows three steps. The first two steps show that $f(\bar{\bbeta}_1)$ is Lipshitz differentiable and strongly convex on the constraint set $F$, and the last step utilizes the classical convex gradient analysis.

\textbf{Step One:} Verify $f(\bar{\bbeta}_1)$ is $L$-Lipschitz differentiable. For any $\bar{\bbeta}_1^{(1)}$ and $\bar{\bbeta}_1^{(2)}$ whose support belong to $F$, 
\begin{eqnarray*}
	\nabla f(\bar{\bbeta}_1^{(1)})-\nabla f(\bar{\bbeta}_1^{(2)}) = \frac{(\sqrt[3]{\eta})^2}{\phi}\frac{2}{n}\sum_{i=1}^n\Big(\bu_i^{\top}(\bar{\bbeta}_1^{(1)}-\bar{\bbeta}_1^{(2)})(\bv_i^{\top}\bar{\bbeta}_2^*)^2(\bw_i^{\top}\bar{\bbeta}_3^*)^2\Big)\bu_{i}.
\end{eqnarray*}
Then, there exist $\pi\in\mathbb S^{s-1}$ such that
\begin{eqnarray*}
	&&\Big\|\nabla f(\bar{\bbeta}_1^{(1)})-\nabla f(\bar{\bbeta}_1^{(2)})\Big\|_2 \\
	&=& \frac{(\sqrt[3]{\eta})^2}{\phi}\Big|\frac{1}{n}\sum_{i=1}^n\Big(\bu_i^{\top}(\bar{\bbeta}_1^{(1)}-\bar{\bbeta}_1^{(2)})(\bv_i^{\top}\bar{\bbeta}_2^*)^2(\bw_i^{\top}\bar{\bbeta}_3^*)^2\Big)\bu_{i}^{\top}\pi\Big|.
\end{eqnarray*}
Applying Lemma \ref{lemma:tensor_spectral} with multiplying $(\bar{\bbeta}_1^{(1)}-\bar{\bbeta}_1^{(2)}) \circ\bar{\bbeta}_2^*\circ \bar{\bbeta}_3^* $, it shows
\begin{eqnarray*}
	&&\Big|\sum_{i=1}^n\Big[(\bu_i^{\top}(\bar{\bbeta}_1^{(1)}-\bar{\bbeta}_1^{(2)})(\bu_i^{\top}\pi)(\bv_i^{\top}\bar{\bbeta}_2^*)^2(\bw_i^{\top}\bar{\bbeta}_3^*)^2)\Big]\Big|\\
	&\leq& \Big(1+\delta_{n,p,s}\Big)\Big\|\bar{\bbeta}_1^{(1)}-\bar{\bbeta}_1^{(2)}\Big\|_2\eta^{*\tfrac{4}{3}},
\end{eqnarray*}
with probability at least $1-10/n^3$, where $\delta_{n,p,s}$ is defined in \eqref{eqn:initial_stat_rate}. Under Condition \eqref{con:sample} with some constant adjustments, we obtain
\begin{eqnarray}
	\Big\|\nabla f(\bar{\bbeta}_1^{(1)})-\nabla f(\bar{\bbeta}_1^{(2)})\Big\|_2 \leq\frac{57}{16}\Big\|\bar{\bbeta}_1^{(1)}-\bar{\bbeta}_1^{(2)}\Big\|_2.
\end{eqnarray}
with probability at least $1-10/n^3$. Therefore, $f(\bar{\bbeta}_1)$ is Lipschitz differentiable with Lipschitz constant $L=\frac{57}{8}$.

\textbf{Step Two:} Verify $f(\bar{\bbeta}_1)$ is $\alpha$-strongly convex. It is equivalent to prove that $\nabla^2f(\bar{\bbeta}_1)\succeq m \mI_p$. Based on the inequality (3.3.19) in \cite{HJ88}, it shows that
\begin{equation}
	\label{minimum_eigen}
	\begin{split}
	\lambda_{\min}\big(\nabla^2(f(\bar{\bbeta}_1))\big)\geq \lambda_{\min}\big(\mathbb E(\nabla^2 f(\bar{\bbeta}_1))\big)-\lambda_{\max}\big(\nabla^2 f(\bar{\bbeta}_1)-\mathbb E(\nabla^2 f(\bar{\bbeta}_1)\big).
	\end{split}
\end{equation}
The lower bound of $\lambda_{\min}(\nabla^2(f(\bar{\bbeta}_1)))$ breaks into two parts: an lower bound for $\lambda_{\min}(\mathbb E(\nabla^2 f(\bar{\bbeta}_1)))$, and an upper bound for $\lambda_{\max}(\nabla^2 f(\bar{\bbeta}_1)-\mathbb E(\nabla^2 f(\bar{\bbeta}_1))$.
The Hessian matrix of $f(\bar{\bbeta}_1)$ is given by
\begin{eqnarray*}
	\nabla^2f(\bar{\bbeta}_1) = \frac{(\sqrt[3]{\eta})^2}{\phi}\frac{2}{n}\sum_{i=1}^n(\bv_1^{\top}\bar{\bbeta}_2^*)^2(\bw_i^{\top}\bar{\bbeta}_3^*)^2\bu_i\bu_i^{\top}.
\end{eqnarray*}
Since $\bu_i, \bv_i,\bw_i$ are independent with each other, we have $\mathbb E(\nabla^2f(\bar{\bbeta}_1)) = 2\mI$, which implies $\lambda_{\min}\big(\mathbb E(\nabla^2f(\bar{\bbeta}_1))\big)\geq2$.

 On the other hand,
\begin{eqnarray*}
	&&\lambda_{\max}\Big(\nabla^2f(\bar{\bbeta}_1)-\mathbb E(\nabla^2f(\bar{\bbeta}_1))\Big) = \Big\|\nabla^2f(\bar{\bbeta}_1)-\mathbb E(\nabla^2f(\bar{\bbeta}_1))\Big\|_2\\
	&&\leq \ba^{\top}\Big(\nabla^2f(\bar{\bbeta}_1)-\mathbb E(\nabla^2f(\bar{\bbeta}_1))\Big)\bb=\frac{2}{n}\sum_{i=1}^n (\bv_i^{\top}\bar{\bbeta}_2^*)^2(\bw_i^{\top}\bar{\bbeta}_3^*)^2(\bu_i^{\top}\ba)(\bu_i^{\top}\bb)\\
	&& - \mathbb E\Big(\sum_{i=1}^n (\bv_i^{\top}\bar{\bbeta}_2^*)^2(\bw_i^{\top}\bar{\bbeta}_3^*)^2(\bu_i^{\top}\ba)(\bu_i^{\top}\bb)\Big)\eta^{*-\tfrac{4}{3}}.
\end{eqnarray*}
where $\ba,\bb\in\mathbb S^{s-1}$. Equipped with Lemma \ref{lemma:tensor_spectral}, it yields that with probability at least $1-10/n^3$,
\begin{eqnarray*}
	\lambda_{\max}\Big(\nabla^2f(\bar{\bbeta}_1)-\mathbb E(\nabla^2f(\bar{\bbeta}_1))\Big)\leq 2\delta_{n,s,p}.
\end{eqnarray*}
 Together with the lower bound of $\lambda_{\min}(\mathbb E(\nabla^2 f(\bar{\bbeta}_1)))$, we have 
$$\lambda_{\min}(\nabla^2f(\bar{\bbeta}_1))\geq 2-2\delta_{n,p,s},
$$
Under Condition \ref{con:sample}, the minimum eigenvalue of Hessian matrix $\nabla^2f(\bar{\bbeta}_1)$ is lower bounded by $\frac{19}{10}$ with probability at least $1-10/n^3$. This guarantees that $f(\bar{\bbeta}_1)$ is strongly-convex with $\alpha= \frac{19}{10}$.

\textbf{Step Three:}
Combining the Lipschitz condition, strongly-convexity and Lemma 3.11 in \cite{Bubeck15}, it shows that
\begin{eqnarray*}
	&&\big(\nabla f(\bar{\bbeta}_1)-\nabla f(\bar{\bbeta}_1^*)^{\top}\big)\big(\bar{\bbeta}_1-\bar{\bbeta}^*\big)\\
	&\geq& \frac{\alpha L}{\alpha+L}\Big\|\bar{\bbeta}_1-\bar{\bbeta}_1^*\Big\|_2^2+\frac{1}{\alpha+L}\Big\|\nabla f(\bar{\bbeta}_1)-\nabla f(\bar{\bbeta}_1^*)\Big\|_2^2.
\end{eqnarray*}
Since the gradient vanishes at the optimal point, the above inequality times $2\mu$ simplifies to
\begin{equation}
	\label{inequality_convex}
	\begin{split}
	-2\mu \nabla f(\bar{\bbeta}_1)^{\top} (\bar{\bbeta}_1-\bar{\bbeta}_1^*)
	\leq -\frac{2\mu \alpha L}{\alpha+L}\Big\|\bar{\bbeta}_1-\bar{\bbeta}_1^*\Big\|_2^2-\frac{2\mu}{\alpha+L}\Big\|\nabla f(\bar{\bbeta}_1)\Big\|_2^2.
	\end{split}
\end{equation}
Now it's sufficient to bound $\|\bar{\bbeta}_1-\bar{\bbeta}_1^*-\mu\nabla f(\bar{\bbeta}_1)\|_2$ as follows
\begin{eqnarray*}
	&&\Big\|\bar{\bbeta}_1-\bar{\bbeta}_1^*-\mu\nabla f(\bar{\bbeta}_1)\Big\|_2^2\\
	 &=& \Big\|\bar{\bbeta}_1^{t}-\bar{\bbeta}_1^*\Big\|_2^2 + \mu^2\Big\|\nabla f(\bar{\bbeta}_1)\Big\|_2^2-2\mu\nabla f(\bar{\bbeta}_1)^{\top}(\bar{\bbeta}_1-\bar{\bbeta}^*)\\
	&\leq& \Big(1-2\mu \frac{\alpha L}{\alpha + L}\Big)\Big\|\bar{\bbeta}_1-\bar{\bbeta}_1^*\Big\|_2^2+\mu\big(\mu-\frac{2}{\alpha+L}\big)\Big\|\nabla f(\bar{\bbeta}_1)\Big\|_2^2.
\end{eqnarray*}
where $L,\alpha$ are Lipschitz constant and strongly convexity parameter, respectively. If $\mu<\frac{80}{361}$, the last term can be neglected and we obtain the desired upper bound,
\begin{equation}\label{bound:I1_nym}
	\Big\|\bar{\bbeta}_1-\bar{\bbeta}_1^*-\mu\frac{\sqrt[3]{\eta}}{\phi}\nabla_{1} \tilde{\cL}\big(\bar{\bbeta}_1, \bar{\bbeta}_2^*, \bar{\bbeta}_3^*\big)\Big\|_2 \leq \Big(1-3\mu\Big)\Big\|\bar{\bbeta}_1-\bar{\bbeta}_1^*\Big\|_2,
\end{equation}
with probability $1-20/n^3$. This ends the proof. \hfill $\blacksquare$\\

\subsubsection{Bounding $I_2$}

For simplicity, we write $\bz_1 = \bar{\bbeta}_1-\bar{\bbeta}_1^*$, $\bz_2 = \bar{\bbeta}_2 - \bar{\bbeta}_2^*$, $\bz_3 = \bar{\bbeta}_3 - \bar{\bbeta}_2^*$. By the definition of noiseless gradient, it suffices to decompose $I_2$ by
\begin{eqnarray*}
	&&\eta^{-\tfrac{1}{3}}\Big\|\nabla_{1}\tilde{\cL}(\bar{\bbeta}_1, \bar{\bbeta}_2^*, \bar{\bbeta}_3^*) - \nabla_{1}\tilde{\cL}(\bar{\bbeta}_1, \bar{\bbeta}_2, \bar{\bbeta}_3)\Big\|_2\\
	&\leq& \Big\|\frac{1}{n}\sum_{i=1}^n(\bu_i^{\top}\bar{\bbeta}_1)(\bv_i^{\top}\bar{\bbeta}_2)(\bv_i^{\top}\bz_2)(\bw_i^{\top}\bar{\bbeta}_3)(\bw_i^{\top}\bz_3)\bu_i\Big\|_2\\
	&&+\Big\|\frac{1}{n}\sum_{i=1}^n(\bu_i^{\top}\bar{\bbeta}_1)(\bv_i^{\top}\bar{\bbeta}_2)(\bv_i^{\top}\bz_2)(\bw_i^{\top}\bar{\bbeta}_3)(\bw_i^{\top}\bar{\bbeta}_3^*)\bu_i\Big\|_2\nonumber\\
	&&+\Big\|\frac{1}{n}\sum_{i=1}^n(\bu_i^{\top}\bar{\bbeta}_1)(\bv_i^{\top}\bar{\bbeta}_2^*)(\bv_i^{\top}\bar{\bbeta}_2)(\bw_i^{\top}\bar{\bbeta}_3)(\bw_i^{\top}\bz_3)\bu_i\Big\|_2\\
	&&+\Big\|\frac{1}{n}\sum_{i=1}^n(\bu_i^{\top}\bz_1)(\bv_i^{\top}\bar{\bbeta}_2^*)(\bv_i^{\top}\bz_2)(\bw_i^{\top}\bar{\bbeta}_3^*)(\bw_i^{\top}\bz_3)\bu_i\Big\|_2\nonumber\\
	&& +\Big\|\frac{1}{n}\sum_{i=1}^n(\bu_i^{\top}\bz_1)(\bv_i^{\top}\bar{\bbeta}_2^*)(\bv_i^{\top}\bz_2)(\bw_i^{\top}\bar{\bbeta}_3^*)^2\bu_i\Big\|_2\\
	&&+\Big\|\frac{1}{n}\sum_{i=1}^n(\bu_i^{\top}\bz_1)(\bv_i^{\top}\bar{\bbeta}_2^*)^2(\bw_i^{\top}\bar{\bbeta}_3^*)(\bw_i^{\top}\bz_3)\bu_i\Big\|_2.\nonumber
\end{eqnarray*}
Repeatedly using Lemma \ref{lemma:tensor_spectral}, we obtain
\begin{equation*}
\begin{split}
	&\eta^{-\tfrac{1}{3}}\Big\|\nabla_{1}\tilde{\cL}(\bar{\bbeta}_1, \bar{\bbeta}_2^*, \bar{\bbeta}_3^*) - \nabla_{1}\tilde{\cL}(\bar{\bbeta}_1, \bar{\bbeta}_2, \bar{\bbeta}_3)\Big\|_2\\
	&\leq \Big(1+\delta_{n,p,s}\Big)\Big[(1+\varepsilon_0)^3\varepsilon_0+(1+\varepsilon_0)^3+(1+\varepsilon_0)^3+\varepsilon_0^2+2\varepsilon_0\Big]\eta^{*\tfrac{4}{3}}\max_j\|\bz_j\|_2\\
	&\leq \frac{5}{2}\Big(1+\delta_{n,p,s}\Big)\eta^{*\tfrac{4}{3}}\max_j\|\bz_j\|_2,
	\end{split}
\end{equation*}
for sufficiently small $\varepsilon_0$ with probability at least $1-60/n^3$. Under Condition \ref{con:sample}, it suffices to get 
\begin{equation}\label{bound:I2_nym}
	\frac{\sqrt[3]{\eta}}{\phi}\Big\|\nabla_{1}\tilde{\cL}(\bar{\bbeta}_1, \bar{\bbeta}_2, \bar{\bbeta}_3) - \nabla_{1}\tilde{\cL}(\bar{\bbeta}_1, \bar{\bbeta}_2^*, \bar{\bbeta}_3^*)\Big\|_2\leq \frac{8}{3}\max_{j=1,2,3}\Big\|\bar{\bbeta}_j-\bar{\bbeta}_j^*\Big\|_2,
\end{equation}
with probability at least $1-6/n$. \hfill $\blacksquare$\\

\subsubsection{Bounding $I_3$}\label{proof_statistical_error}

$I_3$ quantifies the statistical error. By the definition of noiseless gradient and noisy gradient, we have
\begin{eqnarray*}
	&&\frac{\sqrt[3]{\eta}}{\phi}\Big\|\nabla_{1}\tilde{\cL}(\bar{\bbeta}_1,\bar{\bbeta}_2,\bar{\bbeta}_3)-\nabla_{1}\cL(\bar{\bbeta}_1,\bar{\bbeta}_2,\bar{\bbeta}_3)\Big\|_2\\	&=&\frac{(\sqrt[3]{\eta})^2}{\phi}\Big\|\frac{2}{n}\sum_{i=1}^n\epsilon_i(\bv_i^{\top}\bar{\bbeta}_2)(\bw_i^{\top}\bar{\bbeta}_3)\bu_i\Big\|_2.
\end{eqnarray*}

The proof of this part essentially coincides with the proof for symmetric tensor estimation.  Combining Lemmas \ref{cor:concentr_high_order} and \ref{lemma:sub_exp}, we have
\begin{eqnarray*}
	\Big|\frac{2}{n}\sum_{i=1}^n\epsilon_i(\bv_i^{\top}\bar{\bbeta}_2)(\bw_i^{\top}\bar{\bbeta}_3)u_{ij}\Big|\leq C(1+\varepsilon_0)^2\eta^{*\tfrac{2}{3}}\sigma\frac{(\log n)^{\tfrac{3}{2}}}{\sqrt{n}},
\end{eqnarray*}
with probability at least $1-4/n$. Applying union bound over $3s$ coordinates, it suffices to get
\begin{eqnarray*}
	\mathbb P\Big(\max_{j\in[3s]}\Big|\frac{1}{n}\sum_{i=1}^n\epsilon_i(\bv_i^{\top}\bar{\bbeta}_2)(\bw_i^{\top}\bar{\bbeta}_3)u_{ij}\Big|\geq C(1+\varepsilon_0)^2\eta^{*-\tfrac{2}{3}}\sigma\frac{(\log n)^{\tfrac{3}{2}}}{\sqrt{n}}\Big)\leq \frac{12s}{n}.
\end{eqnarray*}
Therefore, we reach
\begin{equation*}
	\frac{\sqrt[3]{\eta}}{\phi}\Big\|\nabla_{1}\tilde{\cL}(\bar{\bbeta}_1,\bar{\bbeta}_2,\bar{\bbeta}_3)-\nabla_{1}\cL(\bar{\bbeta}_1,\bar{\bbeta}_2,\bar{\bbeta}_3)\Big\|_2 \leq 2C\eta^{*-\tfrac{2}{3}}\sigma\sqrt{\frac{3s(\log n)^{3}}{n}},
\end{equation*}
with probability at least $1-12s/n$. \hfill $\blacksquare$\\

\subsubsection{Bounding $I_4$}

According to the definition of thresholding level $h(\bbeta_1)$ in \eqref{def:threhsold_level_nym}, we can bound the square as follows,
\begin{eqnarray*}
	&&\frac{(\sqrt[3]{\eta})^2}{\phi^2}h^2(\bar{\bbeta}_1)= \frac{(\sqrt[3]{\eta})^4}{\phi^2}\frac{4 \log np}{n^2}\sum_{i=1}^n\Big((\bu_i^{\top}\bar{\bbeta}_1)(\bv_i^{\top}\bar{\bbeta}_2)(\bw_i^{\top}\bar{\bbeta}_3)\\
	&&-(\bu_i^{\top}\bar{\bbeta}_1^*)(\bv_i^{\top}\bar{\bbeta}_2^*)(\bw_i^{\top}\bar{\bbeta}_3^*)-\epsilon_i\Big)^2(\bv_i^{\top}\bar{\bbeta}_2)^2(\bw_i^{\top}\bar{\bbeta}_3)^2
\end{eqnarray*}
Based on the basic inequality $(a+b)^2\leq 2(a^2+b^2)$, we have
\begin{eqnarray*}
	&&\Big((\bu_i^{\top}\bar{\bbeta}_1)(\bv_i^{\top}\bar{\bbeta}_2)(\bw_i^{\top}\bar{\bbeta}_3)-(\bu_i^{\top}\bar{\bbeta}_1^*)(\bv_i^{\top}\bar{\bbeta}_2^*)(\bw_i^{\top}\bar{\bbeta}_3^*)-\epsilon_i\Big)^2\\
	&\leq& 2\Big((\bu_i^{\top}\bar{\bbeta}_1)(\bv_i^{\top}\bar{\bbeta}_2)(\bw_i^{\top}\bar{\bbeta}_3)-(\bu_i^{\top}\bar{\bbeta}_1^*)(\bv_i^{\top}\bar{\bbeta}_2^*)(\bw_i^{\top}\bar{\bbeta}_3^*)\Big)^2+2\epsilon_i^2.
\end{eqnarray*}
Denote $I_1$ and $I_2$ corresponding to optimization error and statistical error,
\begin{eqnarray*}
I_1 &=& \frac{(\sqrt[3]{\eta})^4}{\phi^2}\frac{4 \log np}{n^2}\sum_{i=1}^{n}\Big((\bu_i^{\top}\bar{\bbeta}_1)(\bv_i^{\top}\bar{\bbeta}_2)(\bw_i^{\top}\bar{\bbeta}_3)\\
	&&-(\bu_i^{\top}\bar{\bbeta}_1^*)(\bv_i^{\top}\bar{\bbeta}_2^*)(\bw_i^{\top}\bar{\bbeta}_3^*)\Big)^2(\bv_i^{\top}\bar{\bbeta}_2)^2(\bw_i^{\top}\bar{\bbeta}_3)^2\\
	I_2 &=& \frac{(\sqrt[3]{\eta})^4}{\phi^2}\frac{4 \log np}{n^2}\sum_{i=1}^{n}\epsilon_i^2(\bv_i^{\top}\bar{\bbeta}_2)^2(\bw_i^{\top}\bar{\bbeta}_3)^2.
\end{eqnarray*}
Next, $I_1$ is decomposed by some high-order polynomials as follows
\begin{equation}
	\label{eqn:thresholded_I1}
	\begin{split}
	I_1 =& \frac{(\sqrt[3]{\eta})^4}{\phi^2}\frac{4 \log np}{n^2}\Big(\sum_{i=1}^n(\bu_i^{\top}\bz_1)^2(\bv_i^{\top}\bz_2)^2(\bw_i^{\top}\bz_3)^2(\bv_i^{\top}\bar{\bbeta}_2)^2(\bw_i^{\top}\bar{\bbeta}_3)^2\\
	&+\sum_{i=1}^n(\bu_i^{\top}\bz_1)^2(\bv_i^{\top}\bz_2)^2(\bw_i^{\top}\bar{\bbeta}_3^*)^2(\bv_i^{\top}\bar{\bbeta}_2)^2(\bw_i^{\top}\bar{\bbeta}_3)^2\\
	&+\sum_{i=1}^n(\bu_i^{\top}\bz_1)^2(\bv_i^{\top}\bar{\bbeta}_2^*)^2(\bw_i^{\top}\bz_3)^2(\bv_i^{\top}\bar{\bbeta}_2)^2(\bw_i^{\top}\bar{\bbeta}_3)^2\\
	&+\sum_{i=1}^n(\bu_i^{\top}\bz_1)^2(\bv_i^{\top}\bar{\bbeta}_2^*)^2(\bw_i^{\top}\bar{\bbeta}_3^*)^2(\bv_i^{\top}\bar{\bbeta}_2)^2(\bw_i^{\top}\bar{\bbeta}_3)^2\\
	&+\sum_{i=1}^n(\bu_i^{\top}\bar{\bbeta}_1^*)^2(\bv_i^{\top}\bar{\bbeta}_2^*)^2(\bw_i^{\top}\bz_3)^2(\bv_i^{\top}\bar{\bbeta}_2)^2(\bw_i^{\top}\bar{\bbeta}_3)^2\\
	&+\sum_{i=1}^n(\bu_i^{\top}\bar{\bbeta}_1^*)^2(\bv_i^{\top}\bz_2)^2(\bw_i^{\top}\bar{\bbeta}_3^*)^2(\bv_i^{\top}\bar{\bbeta}_2)^2(\bw_i^{\top}\bar{\bbeta}_3)^2\\
	&+\sum_{i=1}^n(\bu_i^{\top}\bar{\bbeta}_1^*)^2(\bv_i^{\top}\bz_2)^2(\bw_i^{\top}\bz_3)^2(\bv_i^{\top}\bar{\bbeta}_2)^2(\bw_i^{\top}\bar{\bbeta}_3)^2\Big).
	\end{split}
\end{equation}
Each term contains the product of Gaussian random vectors form up to power ten. For the first term, by using Lemma \ref{cor:concentr_high_order},
\begin{eqnarray*}
	&&\frac{1}{n}\sum_{i=1}^n(\bu_i^{\top}\bz_1)^2(\bv_i^{\top}\bz_2)^2(\bw_i^{\top}\bz_3)^2(\bv_i^{\top}\bar{\bbeta}_2)^2(\bw_i^{\top}\bar{\bbeta}_3)^2\\
	&\leq& (1+\varepsilon_0)^4\varepsilon_0^4\Big(1+C\frac{(\log n)^5}{\sqrt{n}}\Big)\eta^{*\tfrac{8}{3}}\max_{j=1,2,3}\|\bz_j\|_2^2,
\end{eqnarray*}
with probability at least $1-1/n$. Similar bounds holds for other terms. As long as $n\geq C\log^{10} n$, we have with probability at least $1-7/n$, 
\begin{equation}
	\label{ineqn:bound_I1}
	I_1\leq \frac{7 \log p}{n}\max_{j=1,2,3}\Big\|\bar{\bbeta}_j-\bar{\bbeta}_j^*\Big\|_2^2.
\end{equation}

Now we turn to bound $I_2$. For fixed $\{\epsilon_i\}$, we have,
\begin{eqnarray*}
	&&\Big|\sum_{i=1}^n\epsilon_i^2(\bv_i^{\top}\bar{\bbeta}_2)^2(\bw_i^{\top}\bar{\bbeta}_3)^2- \sum_{i=1}^n\epsilon_i^2\|\bar{\bbeta}_2\|_2^2\|\bar{\bbeta}_3\|_2^2\Big|\\
	&\leq& C(\log n)^2\|\bepsilon^2\|_2\|\bar{\bbeta}_2\|_2^2\|\bar{\bbeta}_3\|_2^2,
\end{eqnarray*}
with probability at least $1-n^{-1}$. Combining with Lemma \ref{lemma:sub_exp},
\begin{equation}\label{ineqn:bound_I2}
	I_2\leq 4\sigma^2\eta^{*\tfrac{4}{3}}\frac{\log p}{n}.
\end{equation}
Putting \eqref{ineqn:bound_I1} and \eqref{ineqn:bound_I2} together, the thresholded effect can be bound by
\begin{equation}\label{bound:I3_nym}
	\frac{\sqrt[3]{\eta}}{\phi}|h(\bbeta_1)|\leq\sqrt{\frac{7 \log np}{n}}\max_{j=1,2,3}\Big\|\bar{\bbeta}_j-\bar{\bbeta}_j^*\Big\|_2+\frac{2\sigma}{(\sqrt[3]{\eta^*})^2}\sqrt{\frac{ \log np}{n}},
\end{equation}
with probability at least $1-8/n$, provided $n\gtrsim(\log n)^{10}$. \hfill $\blacksquare$\\

\subsubsection{Summary}

Putting the upper bounds \eqref{bound:I1_nym}, \eqref{bound:I2_nym} and \eqref{bound:I3_nym} together, we obtain that if step size $\mu$ satisfies $0<\mu<\mu_0$ for some small $\mu_0$,
\begin{eqnarray*}
	\Big\|\sqrt[3]{\eta}\tilde{\bbeta}_1^+-\sqrt[3]{\eta^*}\bbeta_1^*\Big\|_2\leq (1-\frac{\mu}{12})\max_{j=1,2,3}\Big\|\bar{\bbeta}_j-\bar{\bbeta}_j^*\Big\|_2+\mu\frac{3\sigma}{(\sqrt[3]{\eta^*})^2}\sqrt{\frac{3s \log p}{n}},
\end{eqnarray*}
with probability at least $1-12s/n$. This finishes our proof. \hfill $\blacksquare$\\

\section{Matrix Form Gradient and Stochastic Gradient descent}\label{supp_sec:sgd}

\subsection{Matrix Formulation of Gradient}

In this section, we provide detail derivations for \eqref{eqn:symmetric_gradient} and \eqref{eqn:gradient_nonsym}.
\begin{lemma}
	\label{lemma:gradient_equivalence_sym}
	Let $\boldsymbol{\eta} = (\eta_1,\ldots,\eta_K)\in \mathbb R^{K\times 1}, \bX = (\bx_1,\ldots, \bx_n)\in \mathbb R^{p\times n}$ and $\bB = (\bbeta_1, \ldots, \bbeta_K)\in \mathbb R^{p\times K}$. The gradient of symmetric tensor estimation empirical risk function \eqref{eqn:risk_function_sym} can be written in a matrix form as follows 
	\begin{eqnarray*}
		\nabla_{\bB} \cL(\bB, \boldsymbol{\eta}) = \frac{6}{n}[((\bB^{\top}\bX)^{\top})^3\boldsymbol{\eta} - \by]^{\top}[(((\bB^{\top}\bX)^{\top})^2\odot\boldsymbol{\eta}^{\top})^{\top} \odot \bX]^{\top}.
	\end{eqnarray*}
\end{lemma}

\emph{Proof.} First let's have a look at the gradient for $k$-th component,
\begin{eqnarray*}
	\nabla \cL_k(\bbeta_k) = \frac{6}{n}(\sum_{k=1}^K\eta_k(\bx_i^{\top}\bbeta_k)^3-y_i)\eta_k(\bx_i^{\top}\bbeta_k)\bx_i\in \mathbb R^{p\times 1}, \ \ \text{for} \ k=1,\ldots, K.
\end{eqnarray*}
Correspondingly, each part can be written as a matrix form,
\begin{eqnarray*}
	&&((\underbrace{\bB^{\top}\bX}_{K\times n})^{\top})^3\boldsymbol{\eta} - \by\in\mathbb R^{n\times 1}\\
	&&(((\bB^{\top}\bX)^{\top})^2\odot\boldsymbol{\eta}^{\top})^{\top} \odot \bX \in \mathbb R^{pK\times n}.
\end{eqnarray*}
This implies that $[((\bB^{\top}\bX)^{\top})^3\eta - \by]^{\top}[(((\bB^{\top}\bX)^{\top})^2\odot\eta^{\top})^{\top} \odot \bX]^{\top}\in \mathbb R^{1\times pK}$. Note that $\nabla_{\bB} \cL(\bB, \boldsymbol{\eta}) = (\nabla \cL_1(\bbeta_1)^{\top}, \ldots, \nabla \cL_K(\bbeta_K)^{\top})\in \mathbb R^{1\times pK}$. The conclusion can be easily derived. \hfill $\blacksquare$\\


\begin{Lemma}
	\label{lemma:gradient_equivalence_aym}
	Let $\boldsymbol{\eta} = (\eta_1,\ldots,\eta_K)\in \mathbb R^{K\times 1}, \bU = (\bu_1,\ldots, \bu_n)\in \mathbb R^{p_1\times n}, \bV = (\bv_1,\ldots, \bv_n)\in \mathbb R^{p_2\times n}, \bW = (\bw_1,\ldots, \bw_n)\in \mathbb R^{p_3\times n}$ and $\bB_1 = (\bbeta_{11}, \ldots, \bbeta_{1K})\in \mathbb R^{p_1\times K}, \bB_2 = (\bbeta_{21}, \ldots, \bbeta_{2K})\in \mathbb R^{p_2\times K}, \bB_3 = (\bbeta_{31}, \ldots, \bbeta_{3K})\in \mathbb R^{p_3\times K}$. The gradient of non-symmetric tensor estimation empirical risk function \eqref{def:risk_function_nym} can be written in a matrix form as follows 
	\begin{eqnarray*}
		\nabla_{\bB_1}\cL(\bB_1, \bB_2, \bB_3, \boldsymbol{\eta})  = \bD^{\top}(\bC_1^{\top}\odot\bU)^{\top},
	\end{eqnarray*}
	where $\bD = (\bB_1^{\top}\bU)^{\top}*(\bB_2^{\top}\bV)^{\top}*(\bB_3^{\top}\bW)^{\top}\boldsymbol{\eta} - \by$ and $\bC_1 = (\bB_2^{\top}\bV)^{\top}*(\bB_3^{\top}\bW)^{\top}\odot\boldsymbol{\eta}^{\top}$.
\end{Lemma}

\emph{Proof.} Recall that $\{*,\odot\}$ represent Hadamard product and Khatri-Rao product respectively. Then the dimensionality of $\bD, \bC_1, \bC_1\odot \bU$ can be calculated as follows
\begin{eqnarray*}
	&&\bD=\underbrace{(\bB_1^{\top}\bU)^{\top}}_{n\times K}*\underbrace{(\bB_2^{\top}\bV)^{\top}}_{n\times K}*\underbrace{(\bB_3^{\top}\bW)^{\top}}_{n\times K}\boldsymbol{\eta}-\by\in \mathbb R^{n\times 1},\\
	&&\bC_1 = (\bB_2^{\top}\bV)^{\top}*(\bB_3^{\top}\bW)^{\top}\odot\boldsymbol{\eta}^{\top}\in \mathbb R^{n\times K}, \bC_1^{\top}\odot \bU\in\mathbb R^{Kp_1\times n}.
\end{eqnarray*}
Therefore, $$\nabla_{\bB_1}\cL(\bB_1, \bB_2, \bB_3, \boldsymbol{\eta})  = \bD^{\top}(\bC_1^{\top}\odot\bU)^{\top} = (\nabla_1\cL(\bbeta_1)^{\top}, \ldots, \nabla_K\cL(\bbeta_K)^{\top}).$$\hfill $\blacksquare$\\

\subsection{Stochastic Gradient descent}\label{subsec:stochas}

Stochastic thresholded gradient descent is a stochastic approximation of the gradient descent optimization method. Note that the empirical risk function (\ref{eqn:risk_function_sym}) that can be written as a sum of differentiable functions. Followed by \eqref{eqn:symmetric_gradient},  the gradient of (\ref{eqn:risk_function_sym}) evaluated at $i$-th sketching $\{y_i, \bx_i\}$ can be written as
\begin{eqnarray*}
	\nabla_{\bB}\cL_i(\bB, \boldsymbol{\eta})=[((\bB^{\top}\bx_i)^{\top})^3\boldsymbol{\eta} - y_i][(((\bB^{\top}\bx_i)^{\top})^2\odot \boldsymbol{\eta}^{\top})^{\top}\odot \bx_i]^{\top}\in \mathbb R^{1\times pK},
\end{eqnarray*}
Thus, the overall gradient $\nabla_{\bB}\cL_i(\bB, \boldsymbol{\eta})$ defined in \eqref{eqn:symmetric_gradient} can be expressed as a summand of $\nabla_{\bB}\cL_i(\bB, \boldsymbol{\eta})$,
\begin{eqnarray*}
	\nabla_{\bB}\cL_i(\bB, \boldsymbol{\eta})=\frac{1}{n}\sum_{i=1}^n \nabla_{\bB}\cL_i(\bB, \boldsymbol{\eta}).
\end{eqnarray*}
The thresholded step remains the same as Step 3 in Algorithm\ref{alg:overall_symmetric}. Then the symmetric update of stochastic thresholded gradient descent within one iteration is summarized by
\begin{eqnarray*}
	\tv(\bB^{(t+1)}) = \varphi_{\frac{\mu_{SGD}}{\phi}\bh(\bB^{(t)})}\Big(\tv(\bB^{(t)}) -\frac{\mu_{SGD}}{\phi}\nabla_{\bB}\cL_i(\bB^{(t)})\Big).
\end{eqnarray*}

\section{Technical Lemmas}\label{sec:tech_lemma}

\begin{Lemma}
	\label{lemma:tensor_expectation}
	Suppose $\bx\in\mathbb R^p$ is a standard Gaussian random vector. For any non-random vector $\ba, \bb, \bc\in\mathbb R^p$, we have the following tensor expectation calculation,
	\begin{equation}
	\label{eqn:tensor_expectation}
	\begin{split}
	&\mathbb E\Big((\ba^{\top}\bx)(\bb^{\top}\bx)(\bc^{\top}\bx)\bx \circ\bx \circ \bx\Big)\\
	&= \Big(\ba\circ \bb \circ \bc + \ba\circ \bc \circ \bb + \bb \circ \ba\circ\bc+\bb\circ \bc \circ\ba+\bc\circ\bb\circ\ba+\bc\circ\ba\circ\bb\Big)\\
	&+3\sum_{m=1}^{p}\Big(\ba \circ \be_m\circ \be_m(\bb^{\top}\bc)+\be_m\circ\bb\circ\be_m(\ba^{\top}\bc)+\be_m\circ\be_m\circ \bc(\ba^{\top}\bb)\Big),
	\end{split}
	\end{equation}
	where $\be_m$ is a canonical vector in $\mathbb R^p$. 
\end{Lemma}
\emph{Proof.} Recall that for a standard Gaussian random variable $x$, its odd moments are zero and even moments are $\mathbb E(x^6) = 15, \mathbb E(x^4) = 4$. Expanding the LHS of \eqref{eqn:tensor_expectation} and comparing LHS and RHS, we will reach the conclusion. Details are omitted here. \hfill $\blacksquare$\\

\begin{Lemma}
	\label{lemma:tensor_expectation_asy}
	Suppose $\bu\in\mathbb R^{p_1}, \bv\in\mathbb R^{p_2}, \bw\in\mathbb R^{p_3} $ are independent standard Gaussian random vectors. For any non-random vector $\ba\in\mathbb R^{p_1}, \bb\in\mathbb R^{p_2}, \bc\in\mathbb R^{p_3}$, we have the following tensor expectation calculation
	\begin{equation}
		\label{eqn:tensor_expectation_asy}
		\begin{split}
			&\mathbb E\Big((\ba^{\top}\bu)(\bb^{\top}\bv)(\bc^{\top}\bw)\bu \circ\bv \circ \bw\Big)= \ba\circ \bb \circ \bc.
		\end{split}
	\end{equation}
\end{Lemma}
\emph{Proof.} Due to the independence among $\bu, \bv, \bw$, the conclusion is easy to obtain by using the moment of standard Gaussian random variable. \hfill $\blacksquare$\\

Note that in the left side of \eqref{eqn:tensor_expectation}, it involves an expectation of rank-one tensor. When multiplying any non-random rank-one tensor with same dimensionality, i.e., $\ba_1 \circ \bb_1 \circ \bc_1$, on both sides, it will facilitate us to calculate the expectation of product of Gaussian vectors, see next Lemma for details.

\begin{lemma}
	\label{cor:expectation_gaussian_vector}
	Suppose $\bx\in\mathbb R^p$ is a standard Gaussian random vector. For any non-random vector $\ba, \bb,\bc,\bd \in\mathbb R^p$, we have the following expectation calculation
	\begin{eqnarray*}
	&&\mathbb E(\bx^{\top}\ba)^6 = 15\|\ba\|_2^6,\\
	&&\mathbb E(\bx^{\top}\ba)^5(\bx^{\top}\bb) = 15\|\ba\|_2^4(\ba^{\top}\bb),\\
	&&\mathbb E(\bx^{\top}\ba)^4(\bx^{\top}\bb)^2 = 12\|\ba\|_2^2(\ba^{\top}\bb)^2+3\|\ba\|_2^4\|\bb\|_2^2,\\
	&&\mathbb E(\bx^{\top}\ba)^3(\bx^{\top}\bb)^3 = 6(\ba^{\top}\bb)^3+9(\ba^{\top}\bb)\|\ba\|_2^2\|\bb\|_2^2,\\
	&&\mathbb E(\bx^{\top}\ba)^3(\bx^{\top}\bb)^2(\bx^{\top}\bc) = 6(\ba^{\top}\bb)^2(\ba^{\top}\bc)+6(\ba^{\top}\bb)(\bb^{\top}\bc)(\ba^{\top}\ba)\\
	&&+3(\ba^{\top}\bc)(\bb^{\top}\bb)(\ba^{\top}\ba),\\
	&&\mathbb E(\bx^{\top}\ba)^2(\bx^{\top}\bb)(\bx^{\top}\bc)^2(\bx^{\top}\bd) = 2(\ba^{\top}\bc)^2(\bb^{\top}\bd) + 4(\ba^{\top}\bc)(\bb^{\top}\bc)(\ba^{\top}\bd)\\
	&&+6(\ba^{\top}\bc)(\ba^{\top}\bb)(\bc^{\top}\bd)+3(\bc^{\top}\bx)(\bb^{\top}\bd)(\ba^{\top}\ba).
	\end{eqnarray*}
\end{lemma}

\emph{Proof.} Note that $\mathbb E((\bx^{\top}\ba)^3(\bx^{\top}\bb)^3) = \mathbb E ((\bx^{\top}\ba)^3\langle\bx\circ \bx\circ \bx, \bb \circ \bb \circ \bb \rangle)$. Then we can apply the general result in Lemma \ref{lemma:tensor_expectation}. Comparing both sides, we will obtain the conclusion. Others part follows the similar strategy. \hfill $\blacksquare$\\

Next lemma provides a probabilistic concentration bound for non-symmetric rank-one tensor under tensor spectral norm.
\begin{Lemma} \label{lemma:tensor_concentration}
	Suppose $\bX =(\bx_1^{\top}, \cdots, \bx_n^{\top})^{\top}, \bY =(\by_1^{\top}, \cdots, \by_n^{\top})^{\top}, \bZ =(\bz_1^{\top}, \cdots, \bz_n^{\top})^{\top}$ are three $n\times p$ random matrices. The $\psi_2$-norm of each entry is bounded, s.t. $\|X_{ij}\|_{\psi_2}= K_x, \|Y_{ij}\|_{\psi_2}= K_y, \|Z_{ij}\|_{\psi_2}=K_z$. We assume the row of $\bX,\bY,\bZ$ are independent. There exists an absolute constant $C$ such that,
	\begin{equation*}
	\begin{split}
&\mathbb P\Big(\Big\|\frac{1}{n}\sum_{i=1}^n\Big[\bx_i\circ\by_i\circ\bz_i-\mathbb E(\bx_i\circ\by_i\circ\bz_i)\Big]\Big\|_{s}\geq C K_xK_yK_z\delta_{n,p,s}\Big)\leq p^{-1}.\\
&\mathbb P\Big(\Big\|\frac{1}{n}\sum_{i=1}^n\Big[\bx_i\circ\bx_i\circ\bx_i-\mathbb E(\bx_i\circ\bx_i\circ\bx_i)\Big]\Big\|_{s}\geq C K_x^3\delta_{n,p,s}\Big)\leq p^{-1}.
\end{split}
	\end{equation*}
Here, $\|\cdot\|_{s}$ is the sparse tensor spectral norm defined in \eqref{def:sparse_spectral} and $\delta_{n,p,s}= \sqrt{s\log (ep/s)/n}+\sqrt{s^3\log(ep/s)^3/n^2}$.
\end{Lemma}

\emph{Proof.}
Bounding spectral norm always relies on the construction of the $\epsilon$-net. Since we will bound a sparse tensor spectral norm, our strategy is to discrete the sparse set and construct the $\epsilon$-net on each one. Let us define a sparse set $\cB_0=\{\bx\in \mathbb R^p, \|\bx\|_2=1, \|\bx\|_0\leq s \}$. And let $\cB_{0,s}$ be the $s$-dimensional set defined by $\cB_{0,s}=\{\bx\in\mathbb R^s,\|\bx\|_2=1\}$. Note that $\cB_0$ is corresponding to $s$-sparse unit vector set which can be expressed as a union of subsets of dimension $s$ by expanding some zeros, namely $\cB_0=\cup \ \cB_{0,s}$. There should be at most ${p\choose s}\leq (\frac{ep}{s})^s$ such set $\cB_{0,s}$.

Recalling the definition of sparse tensor spectral norm in \eqref{def:sparse_spectral}, we have
\begin{eqnarray*}
	&&A=\Big\|\frac{1}{n}\sum_{i=1}^n\Big[\bx_i\circ\by_i\circ\bz_i-\mathbb E(\bx_i\circ\by_i\circ\bz_i)\Big]\Big\|_{s}\\
	&=& \sup_{\bchi_1,\bchi_2,\bchi_3\in \cB_0}\Big|\frac{1}{n}\sum_{i=1}^n \Big[\langle\bx_i,\bchi_1\rangle\langle\by_i,\bchi_2\rangle\langle\bz_i,\bchi_3\rangle-\mathbb E(\langle\bx_i,\bchi_1\rangle\langle\by_i,\bchi_2\rangle\langle\bz_i,\bchi_3\rangle)\Big]\Big|.
\end{eqnarray*}
Instead of constructing the $\epsilon$-net on $\cB_0$, we will construct an $\epsilon$-net for each of subsets $\cB_{0,s}$. Define $\cN_{\cB_{0,s}}$ as the $1/2$-set of $\cB_{0,s}$. From Lemma 3.18 in \cite{ledoux2005concentration}, the cardinality of $\cN_{0,s}$ is bounded by $5^s$. By Lemma \ref{lemma:varepsilon_net}, we obtain
\begin{equation}\label{ineqn:sup}
\begin{split}
&\sup_{\bchi_1,\bchi_2,\bchi_3\in \cB_{0,s}}\Big|\frac{1}{n}\sum_{i=1}^n \Big[\langle\bx_i,\bchi_1\rangle\langle\by_i,\bchi_2\rangle\langle\bz_i,\bchi_3\rangle-\mathbb E(\langle\bx_i,\bchi_1\rangle\langle\by_i,\bchi_2\rangle\langle\bz_i,\bchi_3\rangle)\Big]\Big|\\
\leq &2^3\sup_{\bchi_1,\bchi_2,\bchi_3\in \cN_{\cB_{0,s}}}\Big|\frac{1}{n}\sum_{i=1}^n \Big[\langle\bx_i,\bchi_1\rangle\langle\by_i,\bchi_2\rangle\langle\bz_i,\bchi_3\rangle-\mathbb E(\langle\bx_i,\bchi_1\rangle\langle\by_i,\bchi_2\rangle\langle\bz_i,\bchi_3\rangle)\Big]\Big|.
\end{split}
\end{equation}
By rotation invariance of sub-Gaussian random variable, $\langle\bx_i,\bchi_1\rangle$, $\langle\by_i,\bchi_2\rangle$, $\langle\bz_i,\bchi_3\rangle$ are still sub-Gaussian random variables with $\psi_2$-norm bounded by $K_x, K_y, K_z$, respectively.  Applying Lemma \ref{cor:concentr_high_order} and union bound over $\cN_{\cB_{0,s}}$, the right hand side of \eqref{ineqn:sup} can be bounded by
\begin{eqnarray*}
\mathbb P\Big(\text{RHS}\geq 8 K_xK_yK_zC\Big(\sqrt{\frac{\log\delta^{-1}}{n}}+\sqrt{\frac{(\log \delta^{-1})^{3}}{n^2}}\Big)\Big)\leq (5^s)^3\delta,
\end{eqnarray*}
for any $0<\delta<1$.

Lastly, taking the union bound over all possible subsets $\cB_{0,s}$ yields that
\begin{eqnarray*}
&&\mathbb P\Big(A\geq 8 K_xK_yK_zC\Big(\sqrt{\frac{\log\delta^{-1}}{n}}+\sqrt{\frac{(\log \delta^{-1})^{3}}{n^2}}\Big)\Big)\\
&\leq &(\frac{ep}{s})^s(5^s)^3\delta = (\frac{125ep}{s})^s\delta.
\end{eqnarray*}
Letting $p^{-1} = (\frac{125ep}{s})^s\delta$, we obtain with probability at least $1-1/p$
\begin{eqnarray*}
A
\leq C K_xK_yK_z\Big(\sqrt{\frac{s\log (p/s)}{n}}+\sqrt{\frac{s^3\log^3(p/s)}{n^2}}\Big),
\end{eqnarray*}
with some adjustments on constant C. The proof for symmetric case is similar to non-symmetric case so we omit here. 
\hfill $\blacksquare$\\

\begin{Lemma}[(Tensor Covering Number(Lemma 4 in  \cite{nguyen2015tensor} ))]\label{lemma:varepsilon_net}
	Let $\mathbb N$ be an $\epsilon$-net for a set $\bB$ associated with a norm $\|\cdot\|$. Then, the spectral norm of a $d$-mode tensor $\cA$ is bounded by 
	\begin{equation*}
		\begin{split}
	&\sup_{\bx_1,\ldots,\bx_{d-1}\in \bB}\|\cA \times_1\bx_1\ldots \times_{d-1}\bx_{d-1}\|_2\\
	&\leq \Big(\frac{1}{1-\varepsilon}\Big)^{d-1}\sup_{\bx_1\cdots\bx_{d-1}\in\mathbb N}\|\cA \times_1 \bx_1\cdots\times_{d-1}\bx_{d-1}\|_2.
	\end{split}
	\end{equation*}
	This immediately implies that the spectral norm of a $d$-mode tensor $\cA$ is bounded by
	\begin{eqnarray*}
	\|\cA\|_2 \leq (\frac{1}{1-\epsilon})^{d-1}\sup_{\bx_1\ldots\bx_{d-1}\in \cN}\|\cA\times_1\bx_1\cdots\times_{d-1}\bx_{d-1}\|_2,
	\end{eqnarray*}
	where $\mathbb N$ is the $\epsilon$-net for the unit sphere $\mathbb S^{n-1}$ in $\mathbb R^n$.
\end{Lemma}

\begin{Lemma}[(Sub-Gaussianess of the Product of Random Variables)]
	\label{lemma:bounded_gaussian}
	Suppose $X_1$ is a bounded random variable with $|X_1|\leq K_1$ almost surely for some $K_1$ and $X_2$ is a sub-Gaussian random variable with Orlicz norm $\|X_2\|_{\psi_2}K_2$. Then $X_1X_2$ is still a sub-Gaussian random variable with Orlicz norm $\|X_1X_2\|_{\psi_2} = K_1K_2$.
\end{Lemma}

\emph{Proof:} Following the definition of sub-Gaussian random variable, we have
\begin{eqnarray*}
  \mathbb P\Big(\big|X_1X_2\big|>t\Big)=\mathbb P\Big(\big|X_2\big|>\frac{t}{\big|X_1\big|}\Big)\leq \mathbb P\Big(\big|X_2\big|>\frac{t}{\big|K_1\big|}\Big)\leq \exp\Big(1-t^2/K_1^2K_2^2\Big),
\end{eqnarray*}
holds for all $t\geq 0$. This ends the proof.
\hfill $\blacksquare$\\

\begin{Lemma}[(Tail Probability for the Sum of Sub-exponential Random Variables (Lemma A.7 in \cite{CLM16}))]
	\label{lemma:sub_exp}
	Suppose $\epsilon_1,\ldots,\epsilon_n$ are independent centered sub-exponential random variables with
	\begin{eqnarray*}
	\sigma:=\max_{1\leq i \leq n}\|\epsilon_i\|_{\psi_1}.
	\end{eqnarray*}
	Then with probability at least $1-3/n$, we have
	\begin{equation*}
	\begin{split}
	&\Big|\frac{1}{n}\sum_{i=1}^n\epsilon_i\Big|\leq C_0\sigma\sqrt{\frac{\log n}{n}}, \   \big\|\bepsilon\big\|_{\infty}\leq C_0\sigma\log n,\\ &  \Big|\frac{1}{n}\sum_{i=1}^n\epsilon_i^2\Big|\leq C_0 \sigma^2, \Big|\frac{1}{n}\sum_{i=1}^n\epsilon_i^4\Big| \leq C_0 \sigma^4,
	\end{split}
		\end{equation*}
	for some constant $C_0$.
\end{Lemma}

\begin{Lemma}[(Tail Probability for the Sum of Weibull Distributions (Lemma 3.6 in \cite{ALPT2011}))]\label{lemma:alpha_1}
    Let $\alpha\in[1,2]$ and $Y_1,\ldots, Y_n$ be independent symmetric random variables satisfying $\mathbb P(|Y_i|\geq t)=\exp (-t^{\alpha})$. Then for every vector $\ba = (a_1,\ldots, a_n)\in \mathbb R^n$ and every $t\geq 0$,
    \begin{equation*}
        \mathbb P\Big(|\sum_{i=1}^n a_iY_i|\geq t\Big) \leq 2\exp\Big(-c\min\Big(\frac{t^2}{\|\ba\|_2^2}, \frac{t^{\alpha}}{\|\ba\|_{\alpha^*}^{\alpha}}\Big)\Big)
    \end{equation*}
\end{Lemma}
\emph{Proof.} It is a combination of Corollaries 2.9 and 2.10 in \cite{talagrand1994supremum}.

\begin{Lemma}[(Moments for the Sum of Weibull Distributions  (Corollary 1.2 in \cite{bogucki2015suprema}))]\label{lemma:alpha2}
    Let $X_1, X_2,\ldots, X_n$ be a sequence of independent symmetric random variables satisfying $\mathbb P(|Y_i|\geq t)=\exp (-t^{\alpha})$, where $0<\alpha<1$. Then, for $p\geq 2$ and some constant $C(\alpha)$ which depends only on $\alpha$,
    \begin{equation*}
        \left\|\sum_{i=1}^n a_iX_i\right\|_p\leq C(\alpha)(\sqrt{p} \|\ba\|_2+ p^{1/\alpha}\|\ba\|_{\infty}).
    \end{equation*}
\end{Lemma}

\begin{Lemma}[(Stein's Lemma \cite{stein2004use})]\label{lemma:stein_lemma}
    Let $\bx\in\mathbb R^{d}$ be a random vector with joint density function $p(\bx)$. Suppose the score function $\nabla_{\bx}\log p(\bx)$ exists. Consider any continuously differentiable function $G(\bx):\mathbb R^{d_x}\rightarrow \mathbb R$. Then, we have 
    \begin{equation*}
        \mathbb E\Big[G(\bx)\cdot \nabla_{\bx}\log p(\bx)\Big] = -\mathbb E\Big[\nabla_{\bx}G(\bx)\Big].
    \end{equation*}
\end{Lemma}

\begin{Lemma}[(Khinchin-Kahane Inequality (Theorem 1.3.1 in \cite{de2012decoupling}))]\label{lemma:Khi}
    Let $\{a_i\}_{i=1}^n$ a finite non-random sequence, $\{\varepsilon_i\}_{i=1}^n$ be a sequence of independent Rademacher variables and $1<p<q<\infty$. Then 
    \begin{equation*}
        \Big\|\sum_{i=1}^n \varepsilon_ia_i\Big\|_q \leq \Big(\frac{q-1}{p-1}\Big)^{1/2}\Big\|\sum_{i=1}^n\varepsilon_i a_i\Big\|_p.
    \end{equation*}
\end{Lemma}

\begin{Lemma}\label{lemma:proof_incoherence}
    Suppose each non-zero element of $\{\bx_k\}_{k=1}^K$ is drawn from standard Gaussian distribution and $\|\bx_k\|_0\leq s$ for $k\in[K]$. Then we have for any $0<\delta\leq 1$,
    \begin{equation*}
        \mathbb P\Big(\max_{1\leq k_1<k_2\leq K}|\langle \bx_{k_1}, \bx_{k_2} \rangle|\leq C\sqrt{s}\sqrt{\log K+\log 1/\delta}\Big)\geq 1-\delta,
    \end{equation*}
    where $C$ is some constant.
\end{Lemma}
\emph{Proof.} Let us denote $\cS_{k_1k_2}\subset [1,2,\ldots, p]$ as an index set such that for any $i, j\in \cS_{k_1k_2}$, we have $x_{k_1i}\neq 0$ and $x_{k_2j}\neq 0$. From the definition of $\cS_{k_1k_2}$, we know that $|\cS_{k_1k_2}|\leq s$ and $\bx_{k_1}^{\top}\bx_{k_2} = \sum_{j=1}^px_{k_1j}x_{k_2j} = \sum_{j\in\cS_{k_1k_2}}x_{k_1j}x_{k_2j}$. We apply standard Hoeffding's concentration inequality,
\begin{eqnarray*}
 \mathbb P\Big(|\langle \bx_{k_1}, \bx_{k_2} \rangle|\geq t\Big)=\mathbb P\Big(|\sum_{j\in\cS_{k_1k_2}}x_{k_1j}x_{k_2j}|\geq t\Big)\leq e\exp\Big(-\frac{ct^2}{s}\Big).
\end{eqnarray*}
Letting $ct^2/s = \log (1/\delta)$, we reach the conclusion.

\end{document}